\documentclass[11pt,reqno]{amsart}

\usepackage{tikz}
\usepackage{diagbox,multirow}
\usepackage{graphicx}
\usepackage{subcaption}
\usepackage{epstopdf}
\usepackage{epsfig}
\usepackage{math}
\graphicspath{{./figures/},{./figures050125/}} 
\usepackage{mathtools}
\usepackage{tikz}
\usetikzlibrary{shapes,arrows}
\usepackage{algorithm}
\usepackage{algorithmic}
\tikzstyle{decision} = [diamond, draw, fill=blue!20, 
    text width=4.5em, text badly centered, node distance=3cm, inner sep=0pt]
\tikzstyle{block} = [rectangle, draw, fill=blue!20, 
    text width=5em, text centered, rounded corners, minimum height=4em]
\tikzstyle{line} = [draw, -latex']
\tikzstyle{cloud} = [draw, ellipse,fill=red!20, node distance=3cm,
    minimum height=2em]   
\usetikzlibrary{positioning}
\tikzset{main node/.style={circle,fill=blue!20,draw,minimum size=1cm,inner sep=0pt},  }

\newcommand{\mw}{\mathrm{W}}

\newcommand{\wt}{\widetilde}

\newcommand{\ts}{\mathsf{T}}

\newcommand{\mb}{\mathrm{b}}

\newcommand{\mJ}{\mathrm{J}}

\newcommand{\mN}{\mathrm{T}}
\newcommand{\mR}{\mathrm{R}}
\newcommand{\cR}{\mathcal{R}}
\newcommand{\nm}{\mathrm{s}}
\newcommand{\mS}{\mathrm{S}}
\newcommand{\ms}{\nm}
\newcommand{\eps}{\varepsilon}

\newcommand{\cY}{\mathcal{P}}

 \newcommand{\cE}{\mathcal{E}}
 
 \newcommand{\cM}{\mathcal{M}}
  \newcommand{\cA}{\mathcal{A}}
   \newcommand{\cG}{\mathcal{G}}

\newcommand{\En}{\mathbb{E}}

\newcommand{\cS}{\mathcal{S}}
\newcommand{\cH}{\mathcal{H}}
\def\argmin{\mathop{\rm argmin}}

\renewcommand{\R}{\mathbb{R}}

\def\argmax{\mathop{\rm argmax}}

\renewcommand{\H}{\mathbb{H}}

\newcommand{\be}{\begin{equation}}
\newcommand{\ee}{\end{equation}}

\newcommand{\mnote}[1]{\marginpar{\scriptsize \textcolor{red}{WD: #1}}}
\newcommand\eref[1]{(\ref{#1})}

\newcommand{\cth}{\stackrel{\mbox{\tiny{$\circ$}}}{\theta}}
\renewcommand{\e}{\varepsilon}

\newcommand{\cI}{\mathcal{I}}
\newcommand{\nN}{\mathrm{N}}
\newcommand{\cT}{\mathcal{T}}
\renewcommand{\T}{\mathfrak{T}}
\newcommand{\Dt}{\Delta\theta}
\renewcommand{\nN}{\mathrm{f}}

\renewcommand{\X}{\mathbb{X}}
\renewcommand{\V}{\mathbb{V}}
\renewcommand{\U}{\mathbb{U}}
\renewcommand{\Y}{\mathbb{Y}}
\newcommand{\tanp}{\mathfrak{T}}

\newcommand \id[1]  { \,\textrm d{#1}                       }   

\begin{document}
\title[Expansive Natural Neural Gradient Flows]{Expansive Natural Neural Gradient Flows for Energy Minimization}
\author[Dahmen]{Wolfgang Dahmen}
\email{wolfgang.anton.dahmen@googlemail.com}
\address{Department of Mathematics, University of South Carolina}
\author[Li]{Wuchen Li}
\email{wuchen@mailbox.sc.edu}
\address{Department of Mathematics, University of South Carolina}
\author[Teng]{Yuankai Teng}
\author[Wang]{Zhu Wang}
\email{wangzhu@math.sc.edu.}
\address{Department of Mathematics, University of South Carolina}
\begin{abstract}
This paper develops expansive gradient dynamics in deep neural network-induced mapping spaces. Specifically, we generate  tools and concepts for minimizing a class of energy functionals in an abstract Hilbert space setting covering a wide scope of applications 
such as   PDEs-based inverse problems and supervised learning. The approach hinges on a Hilbert space metric in the full diffeomorphism mapping space, which could be 
viewed as  a generalized Wasserstein-2 metric. We then study a projection gradient descent method within deep neural network parameterized sets. More importantly, we develop an adaptation and expanding strategy to step-by-step  enlarge the deep neural network structures. In particular, the expansion mechanism aims to enhance the alignment of the
neural manifold induced natural gradient direction as well as possible with
the ideal Hilbert space gradient descent direction leveraging the fact that we can evaluate projections of the Hilbert space gradient. We demonstrate the efficacy of the proposed
strategy for several simple model problems for energies arising in the context of supervised learning, model reduction, or inverse problems. In particular, we highlight the importance
of assembling the neural flow matrix based on the inner product for the ambient Hilbert space. The actual algorithms are the simplest specifications of a broader spectrum 
based on a correspondingly wider discussion, postponing a detailed analysis to forthcoming work.
\end{abstract}
 \thanks{This work was supported in part by the NSF under an FRG grant DMS-2245097,   DMS-2012469, and DMS-2513112. The work also benefited from the activities organized under the NSF RTG program (grant DMS-2038080). The authors also acknowledge funding by the Deutsche Forschungsgemeinschaft (DFG, German Research Foundation)---project number 442047500---through the Collaborative Research Center ``Sparsity and Singular Structures'' (SFB 1481).
 W. Li's work is supported by AFOSR MURI FP 9550-18-1-502, AFOSR YIP award No. FA9550-23-1-008.} 

\maketitle

\section{Introduction}\label{sec:intro}
\subsection{Motivation}
The vibrant development of {\em Deep Learning} (DL) methodologies is a core foundation of an ever-increasing impact  of  {\em Artificial Intelligence} (AI) in virtually all branches of science, medicine, and social life. An attractive 
perspective is the fact that an enormous scope of tasks, including statistical estimation, supervised and unsupervised learning, classification, state- and parameter estimation, and data assimilation, can be recast in a general  common framework as an {\em optimization problem}    over {\em Deep Neural Networks} (DNNs) as hypothesis classes.
This latter preference is largely due to their remarkable universal expressive power,
see e.g. \cite{Yar,DHP}.

However, on the downside an intrinsic  concern
is a significant {\em uncertainty of optimization success} 
because  
the encountered objective functionals are highly {\em non-convex}. As a result, one may fail to (even nearly) exploit the theoretically possible expressive power and it is usually not possible to estimate the computational work required to warrant
the desired outcome quality. On the other hand, accuracy certification is 
extremely important in technological error averse applications striving for prediction capability.
One of the  
strategies to overcome this hurdle   is - in contrast to traditional estimation postulate - to resort to {\em over-parametrization}. In combination with suitable regularization
(stochastic gradient descent (SGD) in effect is also a regularizer) accuracy is achieved
by driving the loss to zero while maintaining an acceptable generalization error.
However,
in science or technology oriented applications,  DL often aims at generating DNN based 
{\em reduced models} for low-complexity representations which should avoid    over-parametrization. This calls for new
adapted optimization techniques. While in SGD gradients are computed with regard
to Cartesian coordinates, an interesting alternative is to follow gradient flows
with respect to a model dependent metric induced by the network architecture.
Roughly speaking,  denoting by $\Theta$ the architecture and weight budget for a class
of neural networks,
one views the corresponding {\em neural manifold} $\cM(\Theta)$, comprised of all functions
generated when traversing $\Theta$, as (nearly) a  
  Riemannian manifold. One then aims at expediting  the underlying optimization task by following a {\em natural gradient flow}, see e.g. \cite{AGS,IG2,AG,Cayci,JVZ,Peher1,Peher2} for a non-exhaustive list of contributions.
Strictly speaking, reference to a metric, however, constrains the type of admissible
network architectures so that we prefer to adopt a somewhat wider interpretation
giving rise to flows that we refer to in what follows as {\em Neural Gradient Flows} (NGF).

Of course,   
 the realization of such natural or neural gradient flows   comes at a price
of (significantly) higher computational cost. A justification is therefore based
on the hope that corresponding descent steps are significantly more effective or help escaping from
``mediocre terrains'' of the objective landscape,

A central role in related works, as well as in this paper, is played by the so-called {\em neural flow matrix}
$G(\theta)$ which is in essence the Gramian of a spanning system of tangential directions
at a point $\nN(\theta)\in \cM(\Theta)$.

\subsection{Our Contributions}\label{ssec:contrib}
We briefly indicate some points that delineate the present developments from previous works.
\\[1.5mm]
\noindent
{\em (a) Scope of Applications:} Most works in this area consider loss functions in the form
of discrete or continuous mean square risks, i.e., (for convenience) the underlying
Hilbert space is either a Euclidean space or a space $L_2$ of square integrable functions. Instead, we outline in Section \ref{sec:energy}
a more general setting revolving around minimizing a convex (often even quadratic) energy in some ambient
(finite or infinite dimensional) Hilbert space $\H$ and describe some examples.
For instance, $\H$ could be a Sobolev space of positive order. The flow matrix $G(\theta)$
has then to be assembled with respect to the underlying inner product in $\H$ which
raises a number of practical issues. The purpose of some of the later experiments is
to highlight related effects. \\[1.5mm]
{\em (b) Regularization:} In most previous works $G(\theta)$ is assumed to be non-singular.
Schemes based on this assumption can then be viewed as explicit Euler schemes for the underlying
natural gradient flow where $G(\theta)$ determines the Riemannian structure of $\cM(\theta)$. In most realistic scenarios, however, especially when dealing with deep networks,
$G(\theta)$ is singular and even significantly rank-deficient. In this paper we 
consistently work under this assumption. Note that even when $G(\theta)$ happens to be regular, one generally cannot control its condition number or just bound the 
smallest eigenvalue away from zero, which heavily impedes descent properties. Hence the dynamical system under consideration
is either very stiff or even a differential algebraic system of positive degree.
Explicit time-stepping is therefore extremely delicate and necessitates very small learning rates. Singularity of $G(\theta)$ is, for instance, also addressed in \cite{Peher1,Peher2}.
Their strategy is ``sketching'', i.e., the gradient is  projected to a randomly chosen subspace of the tangent space at the current NN approximation $\nN(\theta)$,
hoping that the  Gramian (again in an $L_2$-setting) with respect to the reduced spanning system is reasonably well-conditioned.
In the present paper we also resort at times to projecting to a subspace of the full tangent space which is, however,
generated in a different way. As in several previous works our approach can be viewed as
a Levenberg-Marquardt method - a regularized {\em Gau{\ss}-Newton} scheme, see e.g. \cite{Cayci}. While the regularization controls small singular values, we have observed in numerous cases that NGF turns 
unstable even causing divergence when the largest eigenvalue of $G(\theta)$
is very large. This can only happen when the largest diagonal entry of $G(\theta)$ is very large, indicating a very sensitive dependence on this particular weight.
As described later in Section \ref{ssec:issues}, our regularization is depends  on the size of the largest diagonal entry of $G(\theta)$. In particular, an option is to exclude
hyper-sensitive weights and corresponding tangent directions to reduce the projection space.\\[1.5mm]
{\em (c) Extension and Alignment:} NGF in combination with whatever regularization by no means excludes that the objective may cease to decrease. One reason for this to happen is that
at a current point $\nN(\theta)\in \cM(\Theta)$ the tangent plane to $\cM(\Theta)$ is nearly perpendicular 
(with respect to the inner product of the ambient Hilbert space) to the ideal (Hilbert space gradient) descent direction for the underlying energy, posed over the whole Hilbert space. As will be seen later in the examples,  relevant energies are very benign and 
a fictitious gradient flow in the ambient Hilbert space would converge rapidly.
In such a situation of ``near-perpendiculiarity'' we propose to {\em extend} the network
architecture $\Theta$ in a way that the natural gradient direction is better {\em aligned}
with the ideal Hilbert space gradient direction. Of course, this latter direction cannot be computed but we can compute its (Hilbert space) inner product with the $\theta$-gradient
of $\nN(\theta)$. In other words we can compute its projection to the tangent plane. We discuss in Section \ref{ssec:align} how to exploit this fact so as to enhance
alignment with the ideal descent direction. \\[1.5mm]
\noindent
{\em Objective and layout:} We consider, first on an abstract level, a general setting
centered on minimizing a given strictly convex (often even quadratic)  energy in a given (finite or infinite dimensional) Hilbert space. Corresponding examples
are given in Section \ref{sec:energy}. Ensuing NGF concepts are first formulated in this generality but, so as to bring  out the common mechanisms. In particular, we consider energies that do not just 
live in an $L_2$-setting and emphasize the importance of respecting the correct metric in
the ambient Hilbert space. We illustrate in Section \ref{experiments} key effects by experiments for a representative selection of very simple test problems that nevertheless exhibit the targeted relevant features. In particular, to ensure an acceptable scalability for larger scale problems we employ a hybrid strategy, where NGF is confined to the last two layers while using Adam to update in an interlacing manner all trainable weights. The rationale behind the proposed algorithms
is explained in Section \ref{sec:NGF}. The scope of the discussion, especially in Section \ref{ssec:issues}, is somewhat broader than needed for the specific implementation used in the experiments. We emphasize that these implementations are just the simplest specifications, leaving room for more elaborate (and more expensive) versions, motivated by the reasoning in
Section \ref{ssec:issues}, once principal efficacy has been established. 
While focussing here on rationale and experimental exploration, we postpone a more in depth analysis - beyond the present scope - to a forthcoming paper.

\section{Whenever there is an Energy ...}\label{sec:energy}
On an abstract level, the general problem discussed in this paper has the following
form. Let $\H$ be a Hilbert space with inner product $(\cdot,\cdot)_\H$ and
norm $\|\cdot\|_\H^2= (\cdot,\cdot)_\H$. $\H$ could be finite or infinite dimensional and will be referred to as the ``ambient Hilbert space''. Given 
 a (strictly) convex functional $\cE: \H\to \R$, we wish to solve
\be
\label{1}
\min_{v\in \H} \cE(v).
\ee
A practical method results from restricting the minimization to some finitely parametrized
{\em hypothesis class} $\cH=\cH(\Theta)$. While problem \eqref{1}, when posed over $\H$, is 
convex and - from a theoretical viewpoint - a very benign optimization problem, it becomes
non-convex, when posed over a {\em nonlinear} hypothesis class $\cH(\Theta)$, i.e., over
the finite dimensional coordinate system induced by the trainable wights in $\Theta$.
In what follows we will exploit the favorable properties of the continuous (background) problem \eqref{1}, posed over $\H$.

We briefly indicate next the envisaged scope of application scenarios by a few examples.

\subsection{Supervised Learning}\label{ssec:supervised}
The standard model in supervised learning is to consider a product probability 
space $Z:=X\times Y$ with an unknown measure $\rho$ subject to weak assumptions like 
almost sure bounded support in $Y$. We view the components $y\in Y$ as 
labels associated with $x\in X$. Regression then asks for the ``expected label''
for any given $x$, which is the conditional expectation, viz. or 
regression function,  $f_\rho(x) = \En_\rho\big[y|x]= \int_Y y \id\rho(y|x)$. 
Denoting by $\rho_X$ the marginal of $\rho$ on $X$ and defining
 $\H:= L_2(X,\rho_X)$ with norm $\|g\|_H^2=
\int_X \|g\|^2 \id\rho_X$, 
 it is easy to see
that $f_\rho$  minimizes  
\begin{align}
\label{Elearning}
\frac 12 \int_{Z}\|y- v(x)\|^2 \id\rho(z) & = \frac 12\int_{X\times Y}\|y\|^2 \id\rho 
+ \frac 12 \int_X \|v(x)\|^2 \id\rho_X(x)\nonumber\\
&\qquad\qquad\qquad - \int_X f_\rho(x)\cdot v(x) \id\rho_X(x).
\end{align}
over $v\in \H$.
This, in turn, is equivalent to minimizing the {\em energy}
\be
\label{E2}
\cE(v):= \frac 12 \|v\|_\H^2 - (f_\rho,v)_\H
\ee
over $\H$ - a  first example  of an ``idealized'' energy minimization
in a Hilbert space $\H$.

A common strategy for approximating the unique minimizer of \eqref{E2} is 
to choose
 a hypothesis class $\cH= \{ \nN(\cdot;\theta): \theta \in \Theta\}$ whose
elements are parametrized by a finite number $D$ of ``trainable parameters'' $\theta\in \Theta \subset \R^D$. 
Then the best estimator from $\cH$ would be a minimizer
\be
\label{ideal0}
\min_{\theta\in \Theta}\frac 12 \|f_\rho- \nN(\cdot;\theta)\|^2_\H = 
\min_{\theta\in \Theta} \frac 12\En_{x\sim \rho_X}\big[\|f_\rho(\cdot)- \nN(\cdot;\theta)\|^2\big].
\ee
Since neither is $f_\rho$ known nor can one compute the expectation exactly, 
classical regression in supervised learning tries to estimate $f_\rho$
from {\em independent identically distributed}  (i.i.d.)
samples $\cS=\{z_i=(x_i,y_i):i=1,\ldots,M\}\subset X\times Y$ 
  by minimizing an {\em empirical risk} - a discrete counterpart to \eqref{ideal0} -
\be
\label{ER}
\hat f_\cS := \argmin_{\theta \in \Theta} \frac 1M \sum_{i=1}^M\| y_i- \nN(x_i;\theta)\|^2
\ee
over $\cH$. Here and below 
 $\|\cdot\|$ is the corresponding Euclidean norm when $y$ is vector valued.
This simply results from approximating $\rho$ by the discrete 
measure
\be
\label{rhodiscr} 
\rho_\cS := \frac 1M \sum_{i=1}^M \delta_{z_i},
\ee
which 
allows us to rewrite \eref{ER} as
\be
\label{ER2}
\hat f_\cS := \argmin_{\theta \in \Theta} \int_{X\times Y}\|y - \nN(x;\theta)\|^2 \id\rho_\cS(x,y) = \argmin_{\theta \in \Theta}\En_{x\sim \rho_{\cS,X}}
\big[\|f_\rho(\cdot)- \nN(\cdot;\theta)\|^2\big].
\ee
This again gives rise to an energy minimization problem analogous to \eref{E2}.

In regression 
one is primarily interested in the expected deviation between $f_\rho$ and 
$\hat f_\cS$ or $\cE(\cdot)$ and $\cE_\cS(\cdot)$ as $M=\#\cS \to \infty$,
typically for $M=\#\cS \sim N(\log N)^a$, $N= {\rm dim}\Theta$, or more generally,
$N$ representing the VC-dimension of $\cH$. For corresponding results in the scalar-valued case $Y\subset\R$, see e.g. \cite{BCDDT,BCDD,BaCDD}.  Here our main concern is how to actually
perform the energy minimization when $\cH$ is a nonlinear class.

\subsection{Elliptic PDEs}\label{ssec:ellip}
Solving certain partial differential equations (PDEs) is a seemingly 
different matter. To explain the connection, suppose the symmetric bilinear form $a(\cdot,\cdot): \H\times \H\to \R$ satisfies
\be
\label{ellip}
|a(v,w)|\le C_a\|v\|_\H \|w\|_\H,\quad a(v,v)\ge c_a\|v\|_\H^2,\quad \forall\, v,w\in \H.
\ee
In other words, $\|\cdot\|_a:= a(v,v)^{1/2}$ defines an equivalent inner product
on $\H$. By Lax-Milgram (or Riesz Representation Theorem), the variational problem:
given $\ell\in \H'$ (the dual of $\H$), find $u\in\H$ such that
\be
\label{varprob}
a(u,v)  = \ell(v),\quad \forall\, v\in \H,
\ee
has a unique solution.  Note that this solution is the {\em Riesz lift} of $\ell\in \H'$ to
$\H$, endowed with the norm $\|\cdot\|_a$. Moreover, $u$ is equivalently determined
by the energy minimization
\be
\label{enmin}
u= \argmin_{v\in\H}\Big\{\frac 12 a(v,v)-\ell(v)\Big\}.
\ee
Here the energy 
$$
\cE(v) = \frac 12 a(v,v)-\ell(v)
$$
is again quadratic.
The best known and perhaps simplest example is
the weak formulation of Poisson's equation $-\Delta u=\ell$ on $\Omega\subset\R^d$,
$u|_{\partial\Omega}=0$. In this case the relevant Hilbert space is $\H = H^1_0(\Omega)$
(the Sobolev space of $L_2$ functions that have first order weak derivatives in $L_2$)
whose traces on $\partial\Omega$ vanish.

When trying to approximately solve \eqref{varprob} over a {\em nonlinear} ansatz space
(such as DNNs), \eqref{enmin} offers an optimization-based approach. For moderate spatial dimension $d$
one might argue that there are better ways to solve \eqref{varprob}, but for larger $d$
or when $a(\cdot,\cdot)$ depends on parameters (see the next section), the optimization approach might be a 
competitive alternative.

\subsection{Model Reduction}\label{ssec:modelred}
A further   scenario concerns {\em model reduction}
in Uncertainty Quantification. Here one is typically given a {\em physical 
model} in terms of a {\em parameter dependent} family of PDEs (or more generally
operator equations) written in abstract residual form as
\be
\label{par}
F(u,p)=0,\quad p\in \cY\subset \R^{d_p},
\ee
where $\cY$ is a (most often) high-dimensional {\em parameter domain}. So for
each $p$ one asks for the state $u=u(p)$ that satisfies the PDE \eref{par} for
this $p$ and the collection 
$$
\mathcal{M}:= \{u(p): p\in \cY\}
$$
is often referred to as {\em solution manifold}. It represents the set of viable 
states one expects to encounter under a certain physical law represented
by the PDEs. The uncertainty represented by $\cY$ accounts for unknown constitutive
entities or model deficiencies, while a large parameter dimensionality 
reflects attempts to achieve reasonable model accuracy.

Thus, each state $u(p)\in \cM$ can be viewed as a function of (a small number) $d_x$
of spatial variables and of (a  potentially large number) $d_p$ of parametric variables.
An efficient exploration of $\cM$ calls for reduced modelling, e.g. by constructing
a surrogate for the parameter-to-solution map $\cS:p\mapsto u(p)$. Again, a learning approach
is a way to construct approximations to $\cS$ which requires contriving a suitable 
training loss. 

To that end, suppose $\U$ is a Hilbert space so that posing \eqref{par} over $\U$ is
well-posed for every $p\in \cY$ - in brief $\U$ is {\em model compliant} (see below for more details in this regard). A natural option is
then  to view $p$ as a ``stochastic variable'' and to consider, as a Hilbert space-valued counterpart to \eqref{Elearning}, the quadratic energy
\be
\label{exp}
2\cE(v):=\En\big[ \|u- v \|_\U^2\big]
= \int_\cY \|u(p)- v\|_\U^2 \id\mu(p)= \|u- v\|^2_{L_2(\cY;\U)}
\ee
over $v\in \U$. 
Here the underlying Hilbert space is the {\em Bochner space} $\H = L_2(\cY;\U)$.
This is a purely data-driven approach (see \cite{Stuartetal})
which yields automatically approximations in the ``right'' compliant norm $\H$
(which may have to be approximated by quadrature in the spatial variables).
A disadvantage is that instead of minimizing \eqref{exp}, one minimizes an empirical risk
\be
\label{uy}
\theta^*\in \argmin_{\theta\in \Theta} \frac 1M \sum_{i=1}^M \|u(p_i)- \nN(\theta)\|^2_\U,
\ee
which is equivalent to energy minimization in the ``discrete'' Bochner space 
$\H:= \ell_2(\widehat P;\U)$, where $\widehat P=\{p_i\}_{i=1}^M$ is a finite training set.
This requires the (potentially very costly) 
computation of a typically large number of high-fidelity approximate solutions
$\tilde u(p_i)$ as synthetic training data.
 
A remedy to this latter fact is to contrive  {\em residual losses} that do not
involve the unknown solution $u$. 
To that end, an instructive example is the pressure model for porous media flow 
$F(u,p)= f+ {\rm div}(a(p)\nabla u)=0$, $u|_{\partial\Omega}=0$, which is a parameter dependent version of Poisson's equation. A corresponding high-dimensional counterpart to
\eqref{varprob} reads: find $u\in \H:= L_2(\cY;H^1_0(\Omega))$ such that for $\ell\in \H'$
(i.e., $\U= H^1_0(\Omega)$)
\be
\label{liftpoisson}
b(u,v):= \int_\cY a(p)\nabla u\cdot\nabla v \id x\id \mu(p)= \int_\cY \ell(v)\id\mu(p),
\quad \forall\, v\in \H.
\ee
This is equivalent to the energy minimization
\be
\label{enlift}
u=\argmin_{v\in \H} \cE(v),\quad \cE(v):= \frac 12 b(v,v)- \ell(v).
\ee

\subsection{PINN, GAN, and Saddle Point Problems}\label{ssec:PG}
\eqref{enlift} takes advantage of the fact that each fiber problem $f+ {\rm div}(a(p)\nabla u)=0$
is {\em elliptic}. For more general PDE models, such as time dependent problems, transport
dominated problems or dispersive models, one can no longer resort to an energy minimization
in a straightforward way. Nevertheless, as shown in \cite{BDO,DW} an appropriate {\em minimum-residual} formulation can still be cooked up for a much wider scope of
PDE models. Physics Informed Neural Networks (PINNs) employ residual based
losses as well. But simply measuring residuals in $L_2$ (or $\ell_2$) often 
(as in the case of Poisson's) incurs a variational crime that may result in an 
ill-posed minimization problem. Therefore,  it is suggested in \cite{BDO,DW}
to employ {\em variationally correct} residual losses, which means that the loss
is always {\em uniformly proportional} to the error in a model compliant norm.
For linear PDEs this is ensured by employing a {\em stable variational formulation}.
This requires finding for each $p\in \cY$  a pair $\U_p,\V_p$ of trial and test spaces such that the weak
formulation: find $u=u(p)\in\U_p$ such that 
$$
F(u,p)(v)= b(u,v;p)- \ell(v)=0,\quad \forall\, v\in \V_p
$$
is well-posed uniformly in $p\in\cY$. This is equivalent to $b(\cdot,\cdot;p)$
satisfying the conditions in the Babu$\check{s}$ka-Ne$\check{c}$as-Theorem (continuity, inf-sup, and subjectivity condition). This   in turn, is equivalent to the {\em error-residual relation}
\be
\label{err-res}
\|u(p)- w\|_{\U_\p}\eqsim \|F(w,p)\|_{\V'_p},\quad w\in \U_p,\, p\in \cY.
\ee
This family of fiber problems can be lifted to a single variational problem 
$b(u,v)=\ell(v)$ over $\X:=L_2(\cY;(\U_p)_{p\in\cY})\times \Y:= L_2(\cY;(\V_p)_{p\in\cY})$
in analogy to \eqref{liftpoisson}.
Under mild measurability conditions on the bilinear forms $b(\cdot,\cdot;p)$,   $\X, \Y$ are Hilbert spaces as so called direct integrals, see \cite{BDO}. With
$F(w)(v):= \int_\cY F(w,p)(v) d\mu(p)$, we arive at the {\em variationally correct}
minimum residual problem
\be
\label{var}
\min_{w\in \X}\|F(w)\|_{\Y'}. 
\ee
 
The obvious caveat  here is that (unless $\Y$ is an $L_2$-type space and hence self-dual) the dual norm $\|F(w)\|_{\Y'}:= \sup_{v\in \Y}\frac{F(w)(v)}{\|v\|_\Y}$ cannot be evaluated in a straightforward way but involves 
a supremization over $\Y$. So, quite in the spirit of GAN-concepts (Generative Adversarial Networks), choosing a ``generative parameter budget'' $\Theta$
as well as an ``adversarial budget'' $\Psi$,  \eref{var} boils down to solving a min-max
problem
\be
\label{var2}
\min_{\theta\in \Theta}\max_{\psi\in\Psi}\frac{F(\cG(\theta))(\cA(\psi))}{\|\cA(\psi)\|_\Y},  
\ee
see e.g. \cite{Bertol,BYZZ,Canuto,Friedrichs,Pardo,ZBYZ}. These results reflect the delicacy of this task, in many cases strongly limiting the scope of applicability, e.g. to a low number
of variables. In fact,
it is shown in \cite{DW} that this will indeed approximate the true solution in the right norm
only if $\cA(\psi)$ approximates the true maximizer $v^* = \argmin_{v\in \Y}
F(\cG(\theta))(v)/\|v\|_\Y$, within relative accuracy better than one.
Unfortunately, it is by no means clear that any a priori choice of $\Theta$ and $\Psi$ permits the realization of a sufficiently accurate approximation to $v^*$ 
from the adversarial budget $\Psi$. 
\begin{remark}
This hints at a central point addressed
in this paper, namely the ability to {\em adapt}  network approximations to specific
accuracy requirements, see (c) in Section \ref{ssec:contrib}.
\end{remark}
 To see how this ties into the above framework, one observes first that the true maximizer $v^*$
for a given $\cG(\theta)$, is the Riesz-lift of the functional $F(\cG(\theta))\in \Y'$ to $\Y$, which is actually the solution of an {\em elliptic} variational problem
\be
\label{Riesz}
(v^*,v)_\Y = F(\cG(\theta))(v),\quad \forall\,v\in \Y,
\ee
and hence the minimizer of a quadratic energy
\be
\label{vstar}
v^*=\argmin_{v\in \Y}\frac 12 \|v\|_\Y^2 - F(\cG(\theta))(v),
\ee
see \eqref{varprob}, \eqref{enmin} with $a(\cdot,\cdot)$ replaced by the
$\Y$-inner product $(\cdot,\cdot)_\Y$.
As shown in \cite{DW}, \eqref{var2} can be equivalently formulated as an affine-quadratic
saddle point problem over $\X\times \Y$. For such problems one can formulate {\em proximal 
primal-dual schemes} over $\X\times \Y$ whose convergence rate can be quantified. This, in turn, 
leads to numerical schemes that can be viewed as perturbations of such fictitious 
Hilbert space schemes
that still converge provided that the involved proximal maps are approximately realized
within  judiciously chosen target tolerances, derived from the Hilbert space setting. 
Each proximal map, in turn, requires again minimizing a convex energy where, depending
on the stepsize, the minimum is relatively close to the current iterate. The schemes 
developed in this paper lead to an effective energy decay especially during an initial 
descent phase which one can leverage to meet target tolerances for the proximal maps.
So, any efficient scheme for \eqref{1} can serve as a tool in a primal-dual context that will be addressed
in detail in \cite{DW}.

In summary, we have sketched several different application scenarios leading
to \eqref{1} as the core problem for different Hilbert spaces $\H$. 
Moreover, in all these cases it is crucial to respect the particular Hilbert space norm
in a numerical method.
\section{Natural Gradient Flows}\label{sec:NGF}
\subsection{Hilbert Space Gradient}\label{ssec:genprob}
Recall that the gradient $\nabla \cE\in \H$ is defined  by
\be
\label{2}
\lim_{t\to 0}\frac 1t \big\{\cE(v+th)- \cE(v)\big\} = (\nabla \cE(v),h)_\H,\quad 
\forall\,\, h\in \H.
\ee
Whenever we wish to emphasize that the gradient is taken in $\H$, we write $\nabla_\H$ but drop the subscript if there is no risk of confusion. Obviously,
the path $\{v(t): t>0, v(0)=v_0\}$, given by  
\be
\label{4}
(\dot{v}(t),h)_\H= - (\nabla \cE(v(t)),h)_\H,\quad \forall\, h\in \H,
\ee
``flows'' into the unique minimum and defines our ``ideal'' natural gradient flow.

In particular, note that in all preceding examples   $\cE$ can be cast in the form
\be
\label{Egeneral}
\cE(v) = \frac 12 \|v\|^2_\H - g(v),
\ee
where $g\in \H'$ is a bounded linear functional. The   gradient flow can then be determined explicitly. It follows from \eref{4} that
\be
\label{J0}
(\nabla \cE(v),h)_\H = (v,h)_\H - g(h) = (v- \cR g,h)_\H,
\ee
where $\cR g\in H$ is the Riesz-lift of $g\in \H'$, i.e., $\cR g$ solves
the variational problem
\be
\label{RL}
(\cR g,v)_\H = g(v),\quad \forall\, v\in \H,
\ee 
and is the unique minimizer of $\cE(\cdot)$ over all of $\H$.
From \eref{J0} and \eref{4} one easily verifies that
\be
\label{ODEV2}
\dot{v}(t) = \cR g - v(t) \quad \mbox{hence}\quad v(t)= e^{-t}v_0 +(1-e^{-t})\cR g,\quad t\geq 0.
\ee
Thus, this flow exhibits a very rapid convergence to the minimizer $\cR g$.
Although, when $\H$ is infinite-dimensional, it is not possible to realize such a flow exactly, this will serve as an important orientation when trying
to realize this flow {\em approximately}.
Since the following simple observation will prove quite useful it is worth 
being exposed.
\begin{remark}
\label{rem:proj}
In case \eref{Egeneral} one can compute projections of the Hilbert space gradient 
$\nabla \cE(v)$ through \eqref{J0}
$$
( \nabla \cE(v),h)_\H = ( v,h)_\H - g(h)
$$
with any precision permitted by quadrature on $(\cdot,\cdot)_\H$ and the
application of the functional $g$.
\end{remark}

\subsection{Flow Matrix}
Of course, we cannot compute ideal Hilbert space gradient flows but need to resort
to approximations.
We adopt the rationale and notation of the previous section and choose some {\em hypothesis class} $\cH = \{\nN(\cdot;\theta): 
\theta \in \Theta\}$ for a  practical approach. We defer a description of specific classes
to the next section.

We denote by $\#\theta=D$ the number
of entries in the parameter-array $\theta$. In this sense the parameter space
$\Theta$ can be viewed as a subset (or all) of $\R^{\#\theta}$. A standard approach to approximately solving \eref{1}  is to find
\be
\label{ET}
\theta^*\in \argmin_{\theta\in \Theta}\cE(\nN(\theta)).
\ee
Whenever $\nN(\theta)$ depends {\em nonlinearly} on $\theta\in \Theta$ the resulting objective $E(\theta):= \cE(\nN(\theta))$ is in general no longer convex as a function of $\theta$ (so the $``\argmin''$ is in general set-valued). Thus, the conceptually simple optimization problem in the possibly infinite dimensional space $\H$
is traded against a potentially very difficult optimization problem in finite dimensions over $\R^{\#\theta}$. This is in particular the case 
when employing highly nonlinear families like DNNs as hypothesis classes. 

We briefly recapitulate next some essentially known facts from an angle that 
suits the developments in this paper.
Our strategy to extend  the classical work horse SGD by NGF concepts is
based on the following simple and routinely used
consequence of the chain rule 
\be
\label{grads}
\nabla_\theta E(\theta) = \big(\nabla_\H \cE(\nN(\theta)),
\nabla_\theta \nN(\theta)\big)_\H,
\ee
which is also the basis of the so-called {\em Neural Galerkin} method in \cite{Peher1,Peher2} and in many other related works.
 \begin{remark}
 \label{rem:specproj}
 When $\cE$ is quadratic \eqref{J0} allows one to observe $\nabla_\H\cE(\nN(\theta))$
 through projections on arbitrary $h\in \H$.  
 Even when $\cE$ is not quadratic, and hence can no longer be cast in the form \eqref{Egeneral}, one can still compute projections of $\nabla_\H\cE(\nN(\theta))$ to the tangent
 plane
 $$
 \nN(\theta) +  \tanp(\theta),\quad \tanp(\theta) := {\rm span}\,\big\{
 \partial_{\theta_i}\nN(\theta): i=1,\ldots,D\}
 $$
 of $\cM(\Theta)$ at $\nN(\theta)$ through differentiation in parameter space.
 \end{remark}
 
 In summary, replacing the ideal flow $v(t)$ in \eqref{4} by  $\nN(\theta(t))$, where
 $\theta(t)$ is a vector of parameters, ranging over some budget $\Theta$, the minimization can be viewed as realizing a ``favorable'' parameter flow
 $\dot{\theta}(t)$ which, upon using the chain rule, reads 
\be
\label{parflow}
\big(\nabla_\theta \nN( {\theta}(t))^\top\dot{\theta}(t),h\big)_\H = - \big(\nabla\cE(\nN(\theta(t))),h)_\H,\quad h\in \H.
\ee
Of course, for a finite parameter field $\theta(t)\in \R^D$ this flow does not exist because \eqref{parflow} requires validity under testing with {\em all} elements 
$h$ in the full Hilbert space. The defect
\be
\label{defect}
\nabla_\theta \nN( {\theta}(t))^\top\dot{\theta}(t)+ \nabla\cE(\nN(\theta(t)))
\ee
in $\H$ reflects the inability of the functions $\nN(\theta)$ to follow the ideal
descent path, a viewpoint to be taken up later again.

A natural response is to settle with   a 
{\em projected} parameter flow obtained by testing only with the functions $h= \partial_{\theta_i}\nN(\theta(t))$, $i=1,\ldots,D$, which means to project the ideal flow 
  to the tangent space
  \be
  \label{tan}
  \tanp(\theta):= {\rm span}\,\{\partial_{\theta_i}\nN(\theta): i=1,\ldots,\#\theta\}
  \ee
   at the point $\nN(\theta)$ in the
neural manifold $\cM(\Theta)$.
In view of \eqref{grads}, this yields the dynamical system
 \begin{equation}\label{NGD2}
G_{\nN, \Theta }(\theta)\frac{d\theta}{dt}=
- \nabla_\theta E(\theta),
\ee
%
where 
\be
\label{G}
G_{\nN,\Theta }(\theta) := \big(\nabla_\theta f(\theta),\nabla_\theta \nN(\theta)\big)_\H := \Big(\big(\partial_{\theta_i}\nN(\theta),\partial_{\theta_j}\nN(\theta)\big)_\H\Big)_{i,j=1}^D\in \R^{D\times D}.
\ee
  $G_{\nN,\Theta }(\theta)$ is sometimes referred to   as {\em neural flow matrix}
when $\nN(\theta)$ is a neural network.

 As an example of $G_{\nN,\Theta }(\theta)$,
  consider the case \eref{ER} (or \eref{ER2}) for labels  $y_i
\in \R^d$ so that $\nN(\cdot;\theta)$ is vector-valued $\nN(\cdot;\theta)= (\nN_1(\cdot;\theta),\ldots,\nN_d(\cdot;\theta))^\top$. One then has
$$
G_{\nN,\Theta }(\theta)_{i,j}= \sum_{m=1}^d  \En_{x\sim \rho_X}\big[\partial_{\theta_i}\nN_m(x;\theta)\partial_{\theta_j}\nN_m(x;\theta)\big].
$$ 
%

When $G_{\nN,\Theta }(\theta)$ is nonsingular   it  gives rise to a   Riemannian metric on $\cM(\Theta)$
and, in view of the preceding comments, the system
of ODEs
\be
\label{ODE}
\frac{d}{dt}\theta =- G_{\nN,\Theta}(\theta)^{-1}\nabla_\theta E(\theta) 
\ee
describes the natural gradient flow in the Riemannian manifold $\cM(\Theta)$.

Clearly, $G_{\nN,\Theta }(\theta)$, as a Gramian matrix, is  
symmetric positive semi-definite for each $\theta\in \Theta$, but unfortunately, cannot be expected  in general to be definite, especially when $\nN$ is a DNN. Hence, the frequently used
approach to consider \eref{ODE} in this form
is rather problematic. In fact, a singular or very ill-conditioned $G_{\nN,\Theta}(\theta)$
gives rise to a very stiff problem or even to an {\em algebro-differential} system of positive index so that one encounters {\em stiffness} in either version which is not well in line
with using explicit time stepping.
\subsection{A Least Squares Approach}\label{ssec:common}
\renewcommand{\mN}{\nN}
Recall the following common approach.
 Since the evolution \eqref{ODEV2} is not stiff in $\H$,
 a simple and justifiable approximation of an ideal gradient descent
step \eqref{ODEV2} in $\H$ at time $t$ for a given step size $\tau$ is
$$
v(t+\tau)= v(t)- \tau \nabla \cE(v(t)).
$$
 Thus, for $v(t)= \mN(\theta(t))$
we seek a parameter update $\theta + \Delta\theta$ such that $\mN(\theta +\Delta
\theta)$ is close to $v(t+\tau)$, i.e., $\Delta\theta$ should satisfy
\be
\label{wish}
\mN(\theta + \Delta\theta) - \mN(\theta) \approx -\tau\nabla\cE(\mN(\theta)).
\ee
Upon linearization
$\nN(\theta + \Delta\theta) = \nN(\theta) + \nabla_\theta \nN(\theta) \cdot \Delta\theta
+ o(\|\Delta\theta\|_2)$,  we obtain, in analogy to \eqref{parflow},  
\be
\label{2.2}
\big(\nabla_\theta \mN(\theta) \cdot \Delta\theta  +\tau \nabla
\cE(\mN(\theta)),h\big)_\H =0,\quad \forall\, h\in \H.
\ee
Just as \eqref{parflow}, $\Delta\theta$ is, however, way over-determined by  this relation which therefore can, in general, not be satisfied for
  all $h\in \H$, except for the trivial case $\tau=0$, $\Delta\theta =0$.
  
 If, as commonly done, one projects $\nabla\cE(\nN(\theta))$ to the tangent plane $\tanp(\theta)$ at $\nN(\theta)$, the system
\be
\label{normaleqs}
(\nabla_\theta \mN(\theta),\nabla_\theta \mN(\theta))_\H \Delta\theta =
- \tau\big(\nabla_\theta \mN(\theta),\nabla
\cE(\mN(\theta))\big)_\H = -\tau \nabla_\theta E(\theta), 
\ee
(see \eqref{grads}) always has at least one solution because it gives rise to the 
normal equations of the least squares problem
 \be
 \label{ls}
\min_{\Delta\theta\in \R^D} \big\| \nabla_\theta \mN(\theta) \cdot \Delta\theta + \tau\nabla
\cE(\mN(\theta))\big\|_\H^2.
\ee
 In view of
  \eqref{G}, \eqref{normaleqs} can be restated as
\be
\label{NGD3}
G_{\mN,\Theta}(\theta)\Delta\theta = - \tau\nabla_\theta E(\theta),
\ee
and thus can be viewed as an explicit Euler step   for \eqref{NGD2}.
%
\begin{remark}
\label{rem:hom}
One could divide \eqref{ls} by $\tau$ resulting in an equivalent least squares
problem with stepsize $\tau=1$ adjusting the length of $\Dt$. In practice,
the choice of $\tau$ turns out to matter.
\end{remark} 
 
Since the components $\partial_{\theta_i}\nN(\theta)$ of $\nabla_\theta \nN(\theta)$
are in general  linearly dependent, \eqref{ls} needs to be regularized.
One option is a Gau{\ss}-Newton ansatz, appending to $\nabla_\theta \nN(\theta)$
a fraction of an identity $\lambda I_D$ which results in
\be
\label{GN}
\min_{\Dt\in \R^D} \Big\{ \big\| \nabla_\theta \mN(\theta) \cdot \Delta\theta + \tau\nabla
\cE(\mN(\theta))\big\|_\H^2 +  \lambda^2\|\Dt\|^2\Big\}
\ee
and amounts to solving (see e.g. \cite{Cayci})
\be
\label{GN2}
 \big(G_{\nN,\Theta} +\lambda I\big)\Dt =-\tau \nabla_\theta E(\theta),
\ee
followed by the  update $\nN(\theta)\to \nN(\theta+\Dt)$.

\subsection{Rationale and Orientation}\label{ssec:issues}
\newcommand{\bV}{\mathbf{V}}
\newcommand{\bS}{\mathbf{S}}
\newcommand{\bE}{\mathbf{E}}
\newcommand{\cJ}{\mathcal{J}}
\newcommand{\bg}{\mathbf{g}}
\renewcommand{\W}{\mathbb{W}}
The subsequent discussion not only provides the rationale for the proposed 
algorithm but prepares also the ground for possible extensions or alternatives 
to later specifications.

Of course, we wish to decay the energy in $\cM(\Theta)$ which, due to its complex structure,  is difficult to  analyze. Instead, energy decay within the tangent plane $\nN(\theta)+\tanp(\theta)$ is  somewhat better understood and therefore the central objective 
in this section. Nearby $\nN(\theta)\in \cM(\Theta)$ a decay in $\nN(\theta)+\tanp(\theta)$
is indeed a strong indication of the desired decay in $\cM(\Theta)$ but generally {\em not} 
an {\em implication}, unless one has second order information. 

Keeping this in mind, we collect first some (essentially known)  facts about gradient descent for convex optimization
that  will serve later as an orientation for
specifying various algorithmic structure decisions.  
We assume below that the energy $\cE$
has the following properties which in particular hold for strictly convex quadratic energies:\\[1.2mm]
{\em $\alpha$-strict convexity:} there exists an $\alpha>0$ such that
\be
\label{strict}
\cE(v)\ge \cE(u) +(\nabla\cE(u),v-u)_\H +\frac{\alpha}2 \|u-v\|^2_\H,
\quad \forall\, u,v\in\H.
\ee 
{\em $\beta$-smoothness:} there exists a $\beta>0$ such that
\be
\label{smooth}
\|\nabla \cE(u')- \nabla \cE(u)\|_\H \le \beta\|u-u'\|_\H,\quad \forall\,u, u'\in \H.
\ee

In what follows we always denote by $u^*=\argmin_{u\in\H}\cE(u)$ the exact minimizer in $\H$.
 \medskip

 
\begin{proposition}
\label{prop:descent}
(a) If $\cE$ is $\alpha$-strictly convex and $\beta$-smooth,  one has for
  $0<\tau\le   \frac{1}{\beta}$
\be
\label{decay}
\cE(u-\tau\nabla\cE(u))-\cE(u^*)\le \big(1-\alpha\tau(2- \tau \beta)\big)\big(\cE(u)-\cE(u^*)\big).
\ee
(b) Moreover, assume that $\bg(u)\in\H$ satisfies
\be
\label{gbounds}
c_1 \|\bg(u)\|_\H \le \|\nabla\cE(u)\|_\H \le c_2\|\bg(u)\|_\H
\ee
for some $0< c_1\le c_2 <\infty$ as well as
\be
\label{align}
( \nabla\cE(u),\bg(u))_\H\ge c_3\|\bg(u)\|_\H\|\nabla\cE(u)\|_\H
\ee
for some $c_3>0$.

Then one has
\be
\label{rednew}
\cE(u-\tau\bg(u))-\cE(u^*)\le \Big[1-2\alpha\tau\Big(\frac{c_3}{c_2}- \frac{\tau\beta}{2c_1^2}\Big) \Big]\big(\cE(u)-\cE(u^*) \big).
\ee
and hence a strict energy decay whenever
\be
\label{taucond}
\tau < \frac{2c_3c_1^2}{c_2\beta}.
\ee
\end{proposition}

The claim under (a) follows from standard arguments; (b) can be derived
along similar lines. We refer to   \cite{DW} for a proof.

For instance for 
$$
\cE(v)= \frac 12 \|v\|^2_\H- \ell(v)
$$
one has $\alpha = \beta =1$.

We first quantify how the minimization \eqref{GN} relates to a perturbed gradient descent scheme. To that end, denote for any closed subspace $V$ of $\H$ by   $P_V$ the $\H$-orthogonal projection onto $V$.  We can then write (see \eqref{tan})
\begin{align}
\label{tanproj}
\|\nabla_\theta \nN(\theta)\cdot \Dt +\tau \nabla_\H\cE(\nN(\theta))\|_\H^2 &=
\|\nabla_\theta \nN(\theta)\cdot \Dt +\tau P_{\tanp(\theta)} \nabla_\H\cE(\nN(\theta))\|_\H^2
\nonumber\\
&\qquad\qquad + \|\tau P_{\tanp(\theta)^\perp} \nabla_\H\cE(\nN(\theta))\|_\H^2.
\end{align}
The second summand is independent of $\Dt$ and does not change when minimizing \eqref{ls}
or \eqref{GN}. Hence, it suffices to consider the problem
\be
\label{without}
\min_{\Dt\in\R^D}\Big\{\|\nabla_\theta \nN(\theta)\cdot \Dt +\tau P_{\tanp(\theta)} \nabla_\H\cE(\nN(\theta))\|_\H^2 + \lambda^2\|\Dt\|^2\Big\}.
\ee

\begin{lemma}
\label{lem:svd}
Let
\be
\label{SVD}
G_{\nN,\Theta}(\theta) = \bV \bS^2\bV^\top, 
\ee
be the SVD of the flow matrix, i.e., $\bV$ is an orthogonal matrix and 
$\bS= {\rm diag}\,(s_1,\ldots,s_D)$ is a diagonal matrix with nonnegative entries $s_i$
in decreasing order of magnitude. Moreover, let $D':= \max\{i\le D: s_i>0\}$ and set 
\be
\label{Phi}
\quad \Phi^\top = \nabla_\theta \nN(\theta)^\top\bV\bS^{-1},
\ee 
where $\Phi=\{\phi_1,\ldots,\phi_{D'}\}\subset \H$ is a collection of functions in $\H$. Whenever convenient we view $\Phi$ as a 
column vector $\Phi\in \H^{D'}$ whose entries are the functions $\phi_i$. 

Then, defining 
$$
 \bE:=(E_1,\ldots,E_{D'})^\top\in \R^{D'},\,\, E_i:= (\nabla_\H\cE(\nN(\theta)),\phi_i)_\H,\,\, i=1,\ldots,D',
$$
$\Dt^*:= \bV q^*$ solves \eqref{GN} if and only if  
\be
\label{equivS}
q^*:= \argmin_{q\in \R^{D'}}\Big\{\|\bS q + \tau \bE 
 \big\|^2_{\ell_2}
+ \lambda^2\|  q\|^2_{\ell_2}\Big\}.
\ee
\end{lemma}
\noindent
{\bf Proof:} Consider $T: \R^D\to \tanp(\theta)\subset \H$, defined by
$$
T(\Dt):= \nabla_\theta \nN(\theta)\cdot \Dt :=\nabla_\theta \nN(\theta)^\top\Dt = \sum_{i=1}^D \partial_{\theta_i}\nN(\theta)\Dt_i.
$$
 Since $T$ has finite rank it is compact and hence has a Hilbert-Schmidt representation
$$
T(\Dt) = \sum_{i=1}^{D'} s_i (v_i^\top\Dt )\phi_i =: \Phi^\top \bS \bV^\top\Dt,
$$
where $D'={\rm dim}\,\tanp(\theta)$ and $\Phi=\{\phi_i\}_{i=1}^{D'}$ (viewed as a column vector) is an $\H$-orthonormal basis of $\tanp(\theta)$,  
$\bV = (v_1,\ldots,v_{D'})\in \R^{D\times D'}$ is an orthogonal matrix with columns $v_i$,
 $\bS:= {\rm diag}(s_1,\ldots,s_D)$ is the diagonal matrix with diagonal entries 
$s_i$, and $D'\le D$, when $s_i=0$ for $i>D'$. \eqref{SVD} follows then from the definition of the flow matrix. 

As for the remainder of the claim, recall \eqref{without}, and  note that
\be
\label{projE}
P_{\tanp(\theta)} \nabla_\H\cE(\nN(\theta))= \Phi^\top \big(\nabla_\H\cE(\nN(\theta)),\Phi)_\H =: \Phi^\top\bE,  
\ee
By orthonormality of the $\phi_i$, we clearly have $\|\bE\|_{\ell_2}= \|\nabla_\H\cE(\nN(\theta))\|_\H$.
 Hence, 
we obtain
\begin{align}
\label{deco1}
\|\nabla_\theta \nN(\theta)\cdot \Dt +\tau P_{\tanp(\theta)} \nabla_\H\cE(\nN(\theta))\|_\H^2
&= \big\|\Phi^\top \bS\bV^\top \Dt +\tau \Phi^\top  (\nabla_\H(\nN(\theta)),\Phi)_\H\big\|^2_\H  \nonumber\\
& = \big\|  \bS\bV^\top \Dt +\tau   (\nabla_\H(\nN(\theta)),\Phi)_\H\big\|^2_{\ell_2}.
\end{align}
(Whenever convenient we view $\bV$ as a $D\times D'$ matrix or as 
a $D\times D$ matrix acting on vectors whose entries $D'+1,\ldots,D$ vanish.)
Substituting $\bV^\top\Dt$ (for any $\Dt$ perpendicular to ${\rm ker}\,G_{\nN,\Theta}(\theta)$) by $q\in \R^{D'}$, so that $\|\bV q\|_{\ell_2}= \|q\|_{\ell_2}$, the 
equivalence of \eqref{GN} and \eqref{equivS} follows.\hfill$\Box$\\

Since $\bS$ is diagonal it is easy to compute $q^*$ explicitly.
\begin{lemma}
\label{lem:calc}
The unique minimizer $q^*$ of \eqref{equivS} is
\be
\label{minimizers}
q^*= 
- \tau (\bS^2+\lambda^2 I)^{-1}\bS\bE. 
\ee
Thus, for $\Dt^*:= \bV q^*$, the corresponding correction in the tangent plane is
\begin{align}
\label{corr}
\nN(\theta) + \nabla_\theta \nN(\theta)\cdot \Dt^*&= \nN(\theta) 
-\tau\Phi^\top \bS(\bS^2+\lambda^2I)^{-1}\bS\bE , 
 \end{align}
 where we recall that $\Phi^\top\bE =   P_{\tanp(\theta)}\nabla_H\cE(\nN(\theta))$
 is the projection of the Hilbert space gradient onto $\tanp(\theta)$.
 
Note that for $\lambda=0$ one obtains 
\be
\label{e=0}
q^*:=- \tau \bS^{-1}\bE,\qquad \nN(\theta) + \nabla_\theta \nN(\theta)\cdot \Dt^*=
\nN(\theta)-\tau \Phi^\top \bE.
\ee
\end{lemma}
\noindent
{\bf Proof:} This follows by elementary calculations using the previously derived
relations and the equivalence of \eqref{GN} and \eqref{equivS}.\hfill$\Box$\\

\begin{remark}
\label{rem:special}
 In particular, when $\lambda=0$, and when the $\nabla_\H\cE(\nN(\theta))$ is perfectly aligned with $\tanp(\theta)$, i.e., $P_{\tanp(\theta)}\nabla_H\cE(\nN(\theta))=  \nabla_H\cE(\nN(\theta))$, one obtains, in view of Proposition \ref{prop:descent}, 
 \eqref{decay}, that for $\tau \le 1/\beta$  
$$
\cE(\nN(\theta) -\tau \Phi^\top\bE)-\cE(u^*)\le c(\tau,\alpha,\beta)\big(\cE(\nN(\theta))-
\cE(u^*)\big), 
$$
where $c(\tau,\alpha,\beta)<1$ depends on $\tau,\alpha,\beta$. Although these assumptions
are very unrealistic, this case is instructive as it indicates in which sense
more realistic scenarios can be viewed as a perturbed gradient descent method.
\end{remark}

Obviously, there can only be some energy decay in the tangent plane if $P_{\tanp(\theta)}\nabla_H\cE(\nN(\theta))$ (hence $\bE$) doesn't vanish and the stronger
the tangent space component of $\nabla_H\cE(\nN(\theta))$ the stronger is the decay in
$\tanp(\theta)$. So, a necessary condition for energy decay is that 
there exists   some constant $c_4>0$ such that
\be
\label{portion}
\|P_{\tanp(\theta)}\nabla_H\cE(f(\nN(\theta))\|_\H \ge c_4\|\nabla_H\cE(f(\nN(\theta))\|_\H.
\ee
\begin{lemma}
\label{lem:gen}
Assume that \eqref{portion} holds and define
\be
\label{frac}
 c_5:=\min_{i=1,\ldots,D'}\frac{s_i^2}{\lambda^2+ s_i^2}.
\ee
Moreover, abbreviate
$$
\bg(\theta):= \Phi^\top \bS(\bS^2+\lambda^2I)^{-1}\bS \bE = \Phi^\top \bS q^*.
$$
Then $\bg(\theta)$ is a descent direction for $\cE(\nN(\theta))$ and
\be
\label{descentgen}
\cE(\nN(\theta) -\tau \bg(\theta))-\cE(u^*)\le c_6 \big(\cE(\nN(\theta))-
\cE(u^*)\big)
\ee
holds for 
\be
\label{c6}
 c_6= 1-\alpha\tau\big(c_5 c_4^2(1-\sqrt{1-c_5 c_4^2})-\tau\beta\big),
\ee
provided that
\be
\label{taucond}
\tau < \frac{c_5c_4^2 \big( 1- \sqrt{1- 2c_5 c_4^2}\big)}{\beta}.
\ee
\end{lemma}
\noindent
{\bf Proof:} We view $\bg(\theta)$ as a perturbation of $\nabla_\H\cE(\nN(\theta))$
and apply Proposition \ref{prop:descent} (b). To verify the conditions \eqref{gbounds} and \eqref{align}, note first that, since $\bS(\bS^2+\lambda^2 I)^{-1}\bS$ is diagonal its factors commute so that $\|\bS(\bS^2+\lambda^2I)^{-1}\bS\|= \max\{s^2_i(s_i^2+\lambda^2)^{-1}:i=1,\ldots,D'\}\le 1$. Thus 
\begin{align}
\label{g1}
\|\bg(\theta)\|_\H &= \|\bS(\bS^2+\lambda^2I)^{-1}\bS \bE\|_{\ell_2}\le \|\bE\|_{\ell_2}
= \|P_{\tanp(\theta)}\nabla_\H \cE(\nN(\theta))\|_\H
\nonumber\\
&\le \| \nabla_\H \cE(\nN(\theta))\|_\H.
\end{align}
This establishes the lower estimate in \eqref{gbounds} for $c_1=1$.

Regarding \eqref{align}, we have
\begin{align}
\label{align2}
(\bg(\theta),\nabla_\H\cE(\nN(\theta)))_\H& = \big(\bg(\theta),P_{\tanp(\theta)}
\nabla_\H\cE(\nN(\theta))\big)_\H\nonumber\\
&= (\bS(\bS^2+\lambda^2I)^{-1}\bS\bE,\bE)_{\ell_2}\ge c_5\|\bE\|_{\ell_2} 
\nonumber\\
&\ge c_5 c_4^2 \|\nabla_\H\cE(\nN(\theta))\|^2_\H,
\end{align}
which confirms \eqref{align} with $c_3= c_5 c_4^2$ (because $\|\nabla_\H\cE(\nN(\theta))\|_\H\ge \|\bg(\theta)\|_\H$).

Now we turn to the upper bound in \eqref{gbounds}. We observe first that
(assuming without loss of generality that $2c_5 c_4^2<1$)
\begin{align*}
\|\nabla_\H\cE(\nN(\theta))-\bg(\theta)\|^2_\H &= \|\nabla_\H\cE(\nN(\theta))\|_\H^2
+ \|\bg(\theta)\|_\H^2 - 2(\bg(\theta),\nabla_\H\cE(\nN(\theta)))_\H\\
&\le \|\nabla_\H\cE(\nN(\theta))\|_\H^2
+ \|\bg(\theta)\|_\H^2 - 2c_5c_4^2\|\nabla_\H\cE(\nN(\theta)) \|^2_\H\\
&= (1- 2c_5 c_4^2) \|\nabla_\H\cE(\nN(\theta))\|_\H^2 + \|\bg(\theta)\|_\H^2.
\end{align*}
Thus
\begin{align*}
\|\nabla_\H\cE(\nN(\theta))\|_\H &\le \|\bg(\theta)\|_\H + 
\|\nabla_\H\cE(\nN(\theta))-\bg(\theta)\|_\H\\
& \le \sqrt{1- 2c_5 c_4^2}\|\nabla_\H\cE(\nN(\theta))\|_\H + 2\|\bg(\theta)\|_\H.
\end{align*}
Hence
\be
\label{upper}
\|\nabla_\H\cE(\nN(\theta))\|_\H\le \frac{2}{1- \sqrt{1- 2c_5 c_4^2}} \|\bg(\theta)\|_\H,
\ee
which verifies the upper bound in \eqref{gbounds} with 
\be
\label{c2}
c_2= \frac{2}{1- \sqrt{1- 2c_5 c_4^2}}.
\ee
The assertion follows now from Proposition \ref{prop:descent}
where \eqref{taucond} is just \eqref{rednew} with $c_1=1$, $c_3=c_5 c_4^2$ and $c_2$ from 
\eqref{c2}.\hfill $\Box$

\newcommand{\tanpp}{{\mathfrak{T}_\circ(\theta)}}

\begin{remark}
\label{rem:still}
All the preceding statements remain valid when the (full) tangent space $\tanp(\theta)$
is replaced by some {\em strict} subspace $\tanpp \subset \tanp(\theta)$, provided 
that the constants $c_4$ and $c_5$ refer to this subspace. Of course, alignment
of $\nabla_\H\cE(\nN(\theta))$ with $\tanpp $ may become worse if  $\tanpp $ is
much smaller.
\end{remark}

The following discussion is to distill from these facts an orientation for
subsequent algorithmic specifications.

First, a key role is played by condition \eqref{portion}. If $c_5 c_4^2$ gets very small the admissible step-size $\tau$  in \eqref{taucond}
becomes so small that an effective energy decay (in the tangent plane) is in essence
prohibited. Moreover, the negative impact of $c_4$ being small is particularly strong.
In fact, a small value of $c_4$ means that the ideal Hilbert space descent direction
is nearly orthogonal to the tangent space, i.e., essentially no room for the energy
to decay in a neighborhood of the current point $\nN(\theta)\in \cM(\Theta)$ (Of course,
a too large step-size with respect to \eqref{taucond} may accidentally cause hopping into another energy valley but this is not predictable). Therefore, a too small $c_4$ 
should kick off a {\em network expansion} which is a major point in what follows
and will be addressed in Sections \ref{ssec:expand} and \ref{ssec:align}.

Unfortunately, the value of $c_4$ is not directly assessible. We briefly indicate
a way to at least estimate it. In this regard, a slightly more favorable case is
$\cE(u^*)=0$, i.e., the minimal energy value is zero which we assume first.
We recall the well known Polyak-Lojasiewicz inequality (PL)
\be
\label{PL}
\alpha(\cE(v)-\cE(u^*))  \le\frac 12 \|\nabla\cE(v)\|_\H^2, \quad v\in\H.
\ee
Thus, when $\cE(u^*)=0$ we get
$$
\|\nabla_\H\cE(\nN(\theta))\|_\H\ge \sqrt{2\alpha \cE(\nN(\theta))},
$$
where the right hand side can be evaluated. In principle, one can compute
also $\|P_{\tanp(\theta)}\nabla_\H\cE(\nN(\theta))\|_\H$. If one finds then
\be
\label{c*}
\|P_{\tanp(\theta)}\nabla_\H\cE(\nN(\theta))\|_\H = c_*\sqrt{2\alpha \cE(\nN(\theta))},
\ee
for some positive $c_*$, one would know that $c_4\le c_*$. So, a possible condition that
triggers a network expansion would be that $c_*$ is smaller than a fixed positive threshold $\eta$. Note that if one chooses to compute the basis $\Phi$ according to \eqref{Phi} in Lemma \ref{lem:svd}, we can quickly evaluate $\|P_{\tanp(\theta)}\nabla_\H\cE(\nN(\theta))\|_\H=\|\bE\|$, simplifying determination of $c_*$ in   \eqref{c*}.

If $\cE(u^*)\neq 0$ (as in \eqref{enmin} or \eqref{Egeneral} one can still use that
$\cE(u^*)\le \cE(\nN(\theta)+ \Dt^*\cdot\nabla_\theta \nN(\theta))$ to conclude from 
\eqref{PL} that 
\be
\label{weaker}
2\alpha\big(\cE(\nN(\theta))- \cE(\nN(\theta)+ \Dt^*\cdot\nabla_\theta \nN(\theta))\big)
\le \|\nabla_\H\cE(\nN(\theta))\|_\H^2,
\ee
to compare in an analogous fashion $\|P_{\tanp(\theta)}\nabla_\H\cE(\nN(\theta))\|^2_\H$,
resp. $\|\bE\|^2$, with  $ \alpha\big(\cE(\nN(\theta))- \cE(\nN(\theta)+ \Dt^*\cdot\nabla_\theta \nN(\theta))\big)$.

A computationally simpler and less expensive way to decide when to enlarge the current network, can be based on {\em monitoring} the decay of $\cE$. Specifically, we take 
``stagnation'' of $\cE$ as such a criterion where we postpone the precise meaning
of ``stagnation'', which again depends on the type of energy, to the next section.

Clearly,   $\nabla_\H\cE(\nN(\theta))$ being nearly $\H$-orthogonal to $\tanp(\theta)$
($c_4$ is very small) is one possible reason for a stagnating energy level. It is, however,
not necessarily the only one. In fact, 

(a) the constant $c_5$ from \eqref{frac} could be too small to permit effective step-sizes in \eqref{taucond}, or 

(b) $\Dt^*$ could be too large to keep
the underlying linearization meaningful.   

Regarding (a), evidently, $c_5$ is the smaller the larger
$\lambda$ is in comparison with the smallest positive singular value of $G_{\nN,\Theta}(\theta)$. So, why not taking $\lambda =0$ which in view of \eqref{e=0} looks best anyway? Unfortunately, this may adversely affect (b). In fact, to
 see then the effect of   a very small $s_{D'}$, let us assume for the moment to have the SVD information from Lemma~\ref{lem:svd} available.  
Recall that $\Dt^*= \bV q^* = -\tau \bV \bS^{-1}\bE$ and suppose that the energy in $\bE$ is concentrated on the small singular values. More precisely, let
$\cI= \{D'-k,\ldots,D'-1,D'\}$ and $\cJ= \{1,\ldots,D'\}\setminus \cI$.
Suppose now that for some $1> \eta>0$ we have 
$$
s_i\le \eta,\quad i\in \cI,\quad \|\bE\|_{\ell_2(\cJ)}\le \eta \|\bE\|_{\ell_2(\cI)},
$$
i.e., when $\eta \ll 1$, the dominating part of the energy of  $\bE$ is concentrated on the components in $\cI$
while the large singular values live on the complement $\cJ$. This means
\be
\label{qest}
\|q^*\|_{\ell_2} \ge \|q^*\|_{\ell_2(\cI)}\ge \tau\eta^{-1}\|\bE\|_{\ell_2(\cI)}
\ge \tau (1+\eta^2)^{-1/2}\eta^{-1}\|\bE\|_{\ell_2}.
\ee
Thus, when $\tau\eta^{-1}$ is large, $\|q^*\|_{\ell_2}=\|\bV q^*\|_{\ell_2}=
\|\Dt^*\|_{\ell_2}$ could be too large to comply with linearization. 

A possible strategy can then be summarized as follows:

(i) replace $\tanp(\theta)$ by a subspace $\tanpp$ spanned by $\Phi':= \{\phi_{k_0},
\phi_{k_0+1},\ldots,\phi_{k_1-1},\phi_{k_1}\}$, where $1\le k_0< k_1\le D'$ are chosen
in such a way that $s_{k_0}/s_{k_1}$ stays below a given cap value;

(ii) apply the alignment strategies described below to the resulting subspace $\tanpp:=
{\rm span}\,\Phi'$.

(ii) stagnation of energy decay is then likely to be caused by a poor alignment of
$\nabla_\H\cE(\nN(\theta))$ with $\tanpp$.\\

Following these lines require the SVD of the flow matrix (see Lemma \ref{lem:svd}).
In subsequent experiments we explore whether we can take up the basic mechanisms at a 
lower cost, in particular, avoiding computing many SVDs. In brief: large singular values
are implied by large diagonal entries $G_{\nN,\Theta}(\theta)_{i,i}= \|\partial_{\theta_i}\nN(\theta)\|^2_\H$. Small changes of corresponding weights $\theta_i$ cause large variations in $\cM(\theta)$. So, it seems advisable to reduce the full tangent space to the subspace
$\tanpp$, obtained by discarding the   partials from the spanning set of tangent directions
for which $\|\partial_{\theta_i}\nN(\theta)\|^2_\H$ exceeds a certain cap value.
Then the adverse effect of very small singular values can be countered by a penalization
parameter $\lambda$ which is chosen in dependence of the remaining largest diagonal entry 
of $G_{\nN,\Theta}(\theta)$. We refer to \eqref{eq:lambda} below. A stagnation criterion, depending
on the type of energy, will be used to trigger network expansion and  will be specified for each application example, as described below in Sections \ref{ssec:expand} and
\ref{ssec:align}.

\subsection{Network Architectures}\label{ssec:expand}
\newcommand{\bc}{\mathbf{c}}
\newcommand{\dt}{\delta(\theta)}

We describe now in more detail the type of neural network hypothesis classes
which we focus on. It is well known that deep neural networks (DNNs) are compositions
of  mappings obtained by componentwise acting activation function $\sigma$ on an affine mapping. Typical examples of activation functions are
$$
\sigma(t):= \max\{0,t\} \quad \mbox{(ReLU)},\quad \sigma(t) := \frac 12 \max\{0,t\}^2\quad \mbox{(Squared ReLU)},
$$
and 
\begin{equation*}
    \sigma(t):=\tanh(t)=\frac{e^{t}-e^{-t}}{e^t+e^{-t}} \quad \mbox{(Tanh)}.
\end{equation*}

We concentrate in what follows on a simple ResNet format for the DNNs. The reason is
that, as detailed later, a given network can be viewed as embedded in a network with respect to an expanded
architecture by simply setting expansion weights to zero.
The role of layers is then played by blocks of the form
\be 
\label{block}
z \mapsto { \sigma (\mw z + \mb)\choose z}\mapsto \mJ z +  \sigma (\mw z + \mb)=:\mR(z;\theta),\quad z\in \R^d,\, \mw \in \R^{\bar d\times d},\, \mb\in \R^{\bar d},
\ee
where again $\theta=(\mw,\mb)$, $\mJ\in \R^{\bar d\times d}$ consists of the first $\bar d$ rows of ${\rm id}_d$ when $\bar d< d$, and
\be
\label{if}
\mJ= {{\rm id}_d\choose 0},  \quad \mbox{if}\quad  \bar d \ge d.
\ee
In particular, when $\bar d=d$ the block takes the form
\be
\label{equal}
z \mapsto z+ \sigma(\mw z +\mb)
\ee
and reduces to the identity when $\mb=0$ and $\mw=0$. More generally, it reduces 
to the canonical injection when $\bar d > d$, $\mb=0$ and $\mw=0$. Thus, such 
blocks are just perturbed identities for $\bar d\ge d$ when the trainable parameters, $\mw$ and $\mb$, are small. 

To describe ways of expanding a given architecture $\Theta$ to a larger one, it is
convenient to employ the following notation.

As noted earlier, the elements of $\Theta$ are denoted
by $\theta$, combining the ``nonlinear weights'' denoted by $\cth$ and the coefficient vector $\zeta$ of the closing layer, i.e.,
$$
 \theta = (\cth |\zeta)=(\theta^1|\ldots|\theta^L|\zeta),\quad 
 \mbox{where}\quad \theta^j= (\mw^j,b^j),\,\, \mw^j\in \R^{d_j\times d_{j-1}},\,b^j\in \R^{d_j}.
$$
and $L$ is the depth of the network. 

We consider  ResNet architectures for $\mN(\cdot;\theta)$ based on blocks of the above
type
\be
\label{Rj}
  \mR_j( z;\theta^j) = 
\mJ_j z + \sigma(\mw^j z + \mb^j),\quad \theta^j=(\mw^j, \mb^j), 
\ee
where $\mJ_j$ is the matrix formed by the first $d_j$ rows of the $d_{j-1}\times d_{j-1}$ identity when $d_{j-1}\ge d_j$, and otherwise the $d_{j-1}\times d_{j-1}$-identity
complemented by zero rows. Thus
\be
\label{repT}
\mN(\cdot;\cth |\zeta)=\mN^L(\cdot;\theta ) = \zeta \cdot \nN^L(\cdot;\cth ) = \sum_{i=1}^{d_L} \zeta_i \nN^L_i(\cdot;\cth),
\ee
where $\zeta\in \R^{d_L}$. That is, $\nN^L$ is the vector-valued composition of
all preceding layers
\be
\label{comp}
\nN^L(x;\theta) = \big(\mR_L(\cdot;\theta^L)\circ \cdots\circ \mR_1(\cdot;\theta^1)
\big)(x).
\ee
Roughly speaking
the idea is to enlarge the current budget $\Theta$ in such a way that $\nabla_\H\cE(
\mN(\cdot;\theta_*))$ is ``better seen'' by the new
tangent space.

There are in essence two ways of expanding the current network architecture
$\theta_{(L)} = (\theta^1,\ldots,\theta^L|\zeta)$, namely either adding one or more layers
or extending the current widths, given by the dimensions $d_j$ of the range of the
$j$th layer. \\

\noindent
{\bf Adding a layer:} $\theta_{(L)}=\theta = (\theta^1,\ldots,\theta^L|\zeta)\to 
 \theta_{(L+1)}=(\theta^1,\ldots,\theta^L|\theta^{L+1}|\zeta^L+\xi)$. Here we collect the new trainable weights
in the array 
$$
\psi= (\theta^{L+1}|\xi)= (\mw^{L+1},\mb^{L+1}|\xi)
$$
to obtain 
\be
\label{addl}
\mN^{L+1}(\cdot;\theta|\psi|\xi)= \mN(\cdot;\theta^1|\cdots|\theta^L|\psi|\zeta+\xi).
\ee  
Since we can write
$$
\mN^{L}(\cdot;\theta )=\mN^{L}(\cdot;\cth |\zeta) =   \mN^{L+1}(\cdot;\theta|0|\zeta +0)=\mN^{L+1}(\cdot;\theta^1|\cdots|\theta^L|0|\zeta+0),
$$
  the current network can be written as the expanded network with the 
``new weights''  set to zero.\\

\noindent
{\bf Increasing the width:} To simplify the exposition, we expand only the 
width of the last current layer $R_L(\cdot;\theta_L)$
$$
R_L(z;\theta_L)\mapsto {R_L(z;\theta_L)\choose \sigma(\mw z + b)}
= \mJ_L z+ \sigma(\hat\mw^L + \hat b^j),
$$
where
$$
\hat\mw^L = {\mw^L\choose\mw},\quad \hat b^L = {b^L\choose b}
$$
Again we collect in $\psi=(\mw,b)$ the new ``nonlinear degrees of freedom". 
 They give rise to a network component that can be viewed as a shallow network
 $$
 \mS(z;\psi|\xi)= \xi\cdot \sigma (\mw z + b)
 $$
 appended to the current last ResNet block.
 Note that when adding a layer, formally the original linear closing coefficients $\zeta$ are added to the
new ones $\xi$. When increasing width, the closing vector is {\em extended}
   to $(\zeta,\xi)\in \R^{d_L+ d_e}$ where $d_e$ is the number of added 
 neurons, i.e.,
 \be
 \label{addw}
 \mN(\cdot ;\hat\theta)= \mN(\cdot;\theta)+\xi\cdot \sigma\big(\mw \nN^{L-1}(\cdot;\theta_1|\cdots|\theta^{L-1})
 +b\big).
 \ee
 In this case we have $\mN(\cdot ; \theta)= \mN(\cdot;\cth|0|(\zeta,0))$, i.e., the initial network is still embedded in the expanded one.
 
 In order to have a common representation for the two expansion variants, let
 $$
[\zeta,\xi]:= \left\{
\begin{array}{ll}
\zeta +\xi,& \mbox{in case of \eqref{addl}},\\[1.7mm]
(\zeta,\xi),& \mbox{in case of \eqref{addw}}.
\end{array}\right.
$$
 In what follows we use $\|\cdot\|$ then to denote 
corresponding Euclidean norms.\\[1.6mm]
\subsection{Alignment strategies}
\label{ssec:align}
 Whenever the various criteria, discussed earlier, indicate
 that $\nabla_\H\cE(\nN(\theta))$ is poorly aligned with $\tanp(\theta)$ (the constant $c_4$ is too small), the network should be expanded. We focus on adding a new layer because 
 it is more effective and comes with some subtleties. 
 
 Adhering to 
 the previous conventions, we denote by  $ \theta_*:=(\cth_*|\zeta_*)$ the current network weights. We wish to 
 expand the current network architecture $\Theta$ to a larger one $(\Theta|\Psi)$
 in such a way that the tangent space to $\cM(\Theta|\Psi)$ at the point
 $\mN(\cdot;\theta_*)= \mN(\cdot;\cth_*|0|[\zeta_*,0])\in \cM(\Theta|\Psi)$
 contains a direction that is well-aligned with $\nabla_\H\cE(\mN(\cdot;\theta_*))$.
Here, we collect the ``expansion weights'' in $(\psi|\xi)\in \Psi$.
 
The new tangent space $\T( \cth_*|0 |[\zeta_*,0])$ is spanned by the
partial derivatives of $\mN(\cdot;\cth|\psi|[\zeta,\xi])$ evaluated at 
$(\cth_*|0|[\zeta_*,0])$ which we denote for $\psi=(\mw,b)$ by  
\be
\label{denotepartials}
\begin{array}{ll}
\Gamma_i(\theta_*):= \partial_{\theta_i}\mN(\cdot;\cth|\psi|[\zeta,\xi])\big|_{(\cth_*|0|[\zeta_*,0])},& i=1,\ldots,\#\Theta,\\[1.7mm]
\Gamma_i(\mw):= \partial_{w_{i,j}}\mN(\cdot;\cth|\psi|[\zeta,\xi])\big|_{(\cth_*|0|[\zeta_*,0])},& i,j= 1,\ldots, d_e, d_{L-1},\\[1.7mm]
\Gamma_i(b):= \partial_{b_i}\mN(\cdot;\cth|\psi|[\zeta,\xi])\big|_{(\cth_*|0|[\zeta_*,0])},& i  =1,\ldots, d_e ,\\[1.7mm]
\Gamma_i(\xi):= \partial_{\xi_i}\mN(\cdot;\cth|\psi|[\zeta,\xi])\big|_{(\cth_*|0|[\zeta_*,0])},& i = 1,\ldots, d_e.
\end{array}
\ee
Accordingly, $\Gamma(\psi)$ denotes the gradient portion comprised of the $\Gamma_i(\mw),
\Gamma_i(b)$, $i=1,\ldots,d_e$, (recall $d_e=d_L$ when adding a layer). Likewise
$\Gamma(\psi,\xi)$ is obtained by including also the $\Gamma_i(\xi)$, $i=1,\ldots,d_e$.
 Keep in mind that $\Gamma(\psi,\xi)$ is still a vector whose entries are functions
 of spatial variables in $\H$. In summary, we can group the entries in the gradient
 $$
 \nabla_{\theta,\psi,\xi}\mN(\cdot;\cth|\psi|[\zeta,\xi])\quad \mbox{evaluated at
 $(\cth_*|0|[\zeta_*,0])$}
 $$
  as follows
 $$
 \nabla_{\theta,\psi,\xi}\mN(\cdot;\cth|\psi|[\zeta,\xi])\big|_{(\cth_*|0|[\zeta_*,0])}=
 (\Gamma (\theta_*),\Gamma(\psi,\xi) ).
 $$
In these terms we wish to ideally solve 
\be
\label{ideal}
\min_{\substack{\psi,\zeta\\
\bc\in \R^D, \|\bc\|=1}}\frac{\bc^\top\langle \nabla \cE(\nN(\cdot;\theta_*)),\Gamma(\psi,\zeta)\rangle_\H}{\|\bc^\top\Gamma(\psi,\zeta)\|_\H},
\ee
i.e., we wish to minimize the angle between $-\nabla_\H\cE(\nN(\theta))$ and some linear
combination of the partials of $\nN$ with respect to the new weights in $\bc,(\psi,\zeta)$.

In principle, one can apply SGD to approximately solve \eqref{ideal}. A corresponding
optimization success is of course again uncertain. In favor of mitigating computational cost, we split the task by solving first
\be
\label{split1}
(\hat\psi ,\hat\zeta) \in \argmax_{(\psi,\zeta)\in \Psi} \frac{\big\|( \nabla_\H\cE(\nN(\theta_*)),
\Gamma(\psi,\zeta))_\H\big\|}{\|\Gamma(\psi,\zeta)\|_{\H^{\#(\psi,\zeta)}}},
\ee
which still requires an auxiliary optimization step but of somewhat simpler structure
with fewer variables. Here it is understood that
\begin{align*}
\|\Gamma(\psi,\zeta)\|^2_{\H^{\#(\psi,\zeta)}}&:=\sum_{i=1}^{\#(\psi,\zeta)}
\|\Gamma_i(\psi,\zeta)\|^2_\H,\\
\big\|( \nabla_\H\cE(\nN(\theta_*)),
\Gamma(\psi,\zeta))_\H\big\|^2&:= \sum_{i=1}^{\#(\psi,\zeta)} |(\nabla_\H\cE(\nN(\theta_*)),\Gamma_i(\psi,\zeta))_\H|^2.
\end{align*}
Moreover, the maximal value $Q^*$ of the quotient on the right hand side of \eqref{split1} is one (when $\nabla_\H\cE(\nN(\theta_*))\in \tanp(\cth_*|\psi^*|\zeta^*)$ and otherwise
gives a lower bound for $\|P_{\tanp(\cth_*|\psi^*|\zeta^*)}\nabla_\H\cE(\nN(\theta_*))\|_\H$ (resp. exactly this value when solved exactly). 

Below, we explore an even simpler variant, which is more efficient while retaining positive effects.
For instance, one may randomly choose 
  $R$   samples $(\psi^j,\zeta^j)$, $j=1,\ldots,R$, and take
\be
\label{split1r}
(\hat\psi ,\hat\zeta ) \in \argmax_{1\le j\le R} \frac{\big\|( \nabla_\H\cE(\nN(\theta_*)),
\Gamma(\psi^j,\zeta^j))_\H\big\|}{\|\Gamma(\psi^j,\zeta^j)\|_{\H^{\#(\psi^j,\zeta^j)}}},
\ee
to continue then with \eqref{split2}.

Suppose the extension passes the substance-test (i.e., $Q^*\eqsim \sqrt{2\alpha \cE(\nN(\theta))}$).
Then
\be
\label{split2}
\dt:= \frac{(\Gamma(\hat\psi ,\hat\zeta ),\nabla_\H\cE(\nN(\theta_*)))_\H}{\|(\Gamma(\hat\psi ,\hat\zeta ),\nabla_\H\cE(\nN(\theta_*)))_\H)\|_{\H^{\#(\psi,\zeta)}}}
\ee
has norm one and $\dt^\top \Gamma(\hat\psi ,\hat\zeta )$ maximizes alignment with $\nabla_\H\cE(\nN(\theta))$ 
over all linear combinations 
  $\delta^\top \Gamma(\hat\psi ,\hat\zeta )$ with $\|\delta\|=1$. Note that, on account
  of the chain rule \eqref{grads}, the above expression can be simplified through
 $$
(\Gamma(\hat\psi ,\hat\zeta ),\nabla_\H\cE(\nN(\theta_*)))_\H= \nabla_{\psi,\zeta}
E( \cth_*|0 |[\zeta_*,0]).
$$ 

We need to choose next new weights in the expansion part $\Psi$ of the expanded budget $(\Theta|\Psi)$ that initialize subsequent NGF steps in $(\Theta|\Psi)$. A first
option is to find for some step-size guess $\alpha>0$ 
\be
\label{wt1}
(\wt\psi,\wt\zeta)\in \argmin_{ (\psi|\zeta)\in \Psi 
}\big\|\nN(\cth_*|\psi|[\zeta_*,\zeta]) -\big(\nN(\theta_*)-\alpha \dt^\top \Gamma(\hat\psi ,\hat\zeta )\big)\big\|_\H
\ee
using that $-\dt^\top \Gamma(\hat\psi ,\hat\zeta )$ is a descent direction. Again, we can simplify this task by randomly generating a sequence of parameters $\wt\psi^{(k)}$,    $k=1, \ldots, K$, from
the unit ball in $\Psi$, and take
 \be
\label{wt2}
(\wt\psi|\wt\zeta)\in \argmin_{k=1,\ldots,K}\cE(\nN(\cth_*|\wt\psi^{(k)}|[\zeta_*,\wt\zeta^{(k)}]).
\ee
reducing the minimization in \eqref{wt1} to the finite set
$(\wt\psi^{(k)}|\wt\zeta^{(k)})$, $k=1,\ldots,K$. A disadvantage is that $\alpha$ could
be a poor guess. In our subsequent tests we resort to an  alternative that  exploits the fact that for quadratic energies we can evaluate
inner products of   $-\nabla_\H\cE(\mN(\cdot;\theta_*))$ with {\em any} 
(computable) element in $\H$ (see \eqref{J0}). Specifically, we can choose an orthonormal basis  $\Phi$ for the span of
$\Gamma(\hat\psi,\hat\zeta)$  
and   find 
 $\psi,\xi$ such that the NN update direction is possibly well aligned with the negative Hilbert space gradient 
$-\nabla_\H\cE(\mN(\cdot;\cth_*|0|[\zeta_*,0])) = -\nabla_\H\cE(\mN(\cdot;\theta_*))$. 
Specifically, choosing $\alpha_0>0$, and 
 randomly generating a sequence of parameters $\wt\psi^{(k)}$, for $k=1, \ldots, K$,
 as above, we find corresponding $\wt\zeta^{(k)}, \alpha^k\ge \alpha_0$ solving  
 the simpler problem
\be
\label{alignment}
\min_{\wt\zeta,\alpha\ge \alpha_0}\big\|( \mN^{L+1}(\cdot;\cth_*|\wt\psi^{(k)}|\wt\zeta)-\mN^L(\cdot;\theta_*), \phi_i )_\H 
+ \alpha ( \nabla_\H\cE(\mN(\cdot;\theta_*)), \phi_i   )_\H  \big\|^2.
\ee
We then determine $\wt\psi,\wt\zeta$ by \eqref{wt2}. 

These are only two examples of adaptation strategies. Of course, the effectiveness will
depend on the number of training samples, a larger number increases the probability
of better alignment. Complementing the random search by further SGD steps is another option.
More generally, there are many further conceivable ways of (nearly) solving auxiliary
optimization problems \eqref{alignment} over some ball containing $\psi,\zeta,\alpha$ with $\alpha\ge \alpha_0>0$. An interesting (yet more expensive option)  is to use the random
data $\wt\psi^{(k)}, \wt\zeta^{(k)}$ as (initial) agent positions in a {\em swarm based}
gradient descent method \cite{tadmor}.

\subsection{Structure of the Algorithm}

We now mark the basic structure of the algorithms in terms of two pseudo-codes, namely 
the natural gradient flow step in Alg. \ref{alg:ng}, and the expansion 
 in Alg. \ref{alg:enn}. Several hyper-parameters and specifications will be given later
 in the context of specific applications.

 In order to keep strategies scalable the natural gradient flow steps 
will be applied only to parts of the current network, keeping the remaining trainable weights frozen. To be able to update then all weights, such NG steps are interlaced with 
ADAM steps.

\begin{algorithm}
	\begin{algorithmic}[1]
		\caption{Natural Gradient Flows (NGF)}\label{alg:ng}
		\STATE{\textbf{Input:} Initial neural network with learnable parameters $\theta$, training sample set of $N_b$ mini-batches, initial learning rate $\gamma_0 >0$, maximum number of iterations $M_\text{epoch}$, prescribed criteria {\color{black}(early termination condition and/or loss saturation condition)}.}
		\STATE{\textbf{Output:} {\color{black}Trained model $\nN(\cdot, \theta)$, final loss value $E_k$, stopping flag.}}
		 
		\STATE
		\FOR{$k=1,2,\ldots, M_\text{epoch}$}
		\STATE Initialize loss $E_k = 0$.
		\FOR{$i = 1, \ldots, N_b$}
		\STATE Evaluate loss $E_k^{(i)} := E(\theta)$ on batch $i$.
		\STATE Back-propagate $\nabla_\theta E(\theta)$.
		\STATE Compute neural flow matrix $G_{f, \Theta}(\theta)$.
		\STATE Determine $\lambda$ based on \eqref{eq:lambda}.
		\STATE Set learning rate $\gamma$ by Armijo-based backtracking line search starting from $\gamma_0$.
		\STATE Update learnable parameters $\theta =\theta -\gamma (G_{\nN, \Theta}(\theta)+\lambda \mathbb{I})^{-1}\nabla_\theta E(\theta)$.
		\STATE Accumulate loss: $E_k = E_k + E_k^{(i)}$.
		\ENDFOR
		\STATE Normalize loss: $E_k = E_k/N_b$.
		{\color{black}
		\IF{early termination criterion met}
		\STATE Flag = `early terminated' and break.
		\ELSIF{ loss saturation criterion met}
		\STATE Flag = `saturated' and break.
		\ENDIF 
		}
		\ENDFOR
	\end{algorithmic}
\end{algorithm}

\begin{algorithm}
	\begin{algorithmic}[1]
		\caption{Expansive Neural Network Training}\label{alg:enn}
		\STATE{\textbf{Input:} Initial neural network, index of layer expansion $n_\text{exp} = 1$, maximum number of layer expansions $M_\text{exp}$, prescribed criteria {\color{black}(loss saturation, absolute loss-based early termination and/or relative loss-based early termination )}.}
		\STATE{\textbf{Output:} {\color{black}Trained model $\nN(\cdot, \theta)$.}}
		 
		\STATE
		\STATE Train initial model using Alg.~\ref{alg:ng}, with {\color{black} appropriate stopping criteria (absolute loss-based early termination and/or saturation criteria), store the final loss value to $E^{(0)}$, and output flag. }
		\IF{Flag = `early terminated'}
			\STATE Stop training.
		\ELSIF{Flag = `saturated'} 
			\STATE Expand network by adding a new hidden layer (width $m$) before output layer.
		\ENDIF	
		\WHILE{$n_\text{exp} \leq M_\text{exp}$}
		\STATE Train expanded model (parameters in new layer and output layer are active while the rest is fixed) using Alg.~\ref{alg:ng}, with given criteria {\color{black}(absolute loss-based early termination and/or saturation criteria)} and output flag.
		{\color{black}
		\IF{Flag = `early terminated'}
			\STATE Break. 
		\ELSIF{Flag = `saturated'} 
			\STATE Activate all layers.
			\STATE Train using Adam optimizer with given criteria, output flag, and store the final loss value to $E^{(n_\text{exp})}$.
			\IF{Flag = `early terminated'}
				\STATE Break. 
			\ELSIF{Flag = `saturated'} 
				\STATE {Check whether relative loss-based early termination condition, if given, is satisfied.}
				\IF {relative loss-based early termination met}
				 \STATE Break.
				\ELSE 
				\STATE Expand network by adding a new hidden layer (width $m$) before output layer.
				\ENDIF
			\ENDIF	
		\ENDIF
		}
		\STATE $n_\text{exp} \gets n_\text{exp} + 1$. 
		\ENDWHILE
	\end{algorithmic}
\end{algorithm}

We conclude with specifying ``early termination criteria''. To that end, we distinguish
the cases:
\be
\label{al}
\begin{array}{ll}
\cE(u^*)=0 & \mbox{for a given tolerance $\tau_{t_a}$, terminate the iteration if}\\[1.5mm] 
|E(\theta)| \leq \tau_{t_a}&  { \text{absolute loss-based early termination condition} }
\end{array} 
\ee
and
\be
\label{sl}
\begin{array}{ll}
\cE(u^*)\neq 0 & \mbox{for a given tolerance $\tau_{s_r}$, terminate the iteration if}\\[1.5mm] 
|E_k-E_{\widehat{k}}|<\tau_{s_a} \text{ or }&\\
 |E_k-E_{\widehat{k}}|/|E_{\widehat{k}}|<\tau_{s_r}&  { \text{loss saturation condition}, }
\end{array} 
\ee
where $\widehat{k}$ represents an earlier iteration (typically $\widehat{k} = k - 5$ in our experiments).

For expansive NGF we need in addition a criterion that triggers an expandion step. In the
light of the preceding discussion, loss stagnation or saturation seems to be reasonable 
for either situation $\cE(u^*)=0$ and $\cE(u^*)\neq 0$, i.e., we invoke a satuation condition of the  type \eqref{sl}.

Scenario \eqref{al} is encountered, for instance, in Section \ref{ssec:supervised} as well as for losses of the type  \eqref{uy}, \eqref{var} in Section \ref{ssec:modelred}.
In all those examples the loss is a sharp error bound with respect to the underlying 
model compliant norm. 
Scenario \eqref{sl}, in turn, occurs for instance in Section \ref{ssec:ellip} where the 
loss does not estimate the achieved accuracy and one has to resort to its variation
(in absolute or relative terms, depending on the situation at hand).

\section{Experiments}\label{experiments}

In this section, we evaluate the performance of the proposed optimization method in finding NN models for three tasks (see Section \ref{sec:energy}): function approximation by supervised learning (\textsection~\ref{ssec:supervised}), numerical solution to PDEs by unsupervised learning (\textsection~\ref{ssec:ellip}), and model order reduction by supervised learning (\textsection~\ref{ssec:modelred}). Two types of comparisons will be made: 
\begin{enumerate}
\item First, we train fully-connected neural networks (FCNNs) of depth $L$ (and architecture described in Section \ref{ssec:expand}), with the activation function $\mathsf{Tanh(\cdot)}$, using two optimization methods: Adam and NGF.  
For the former, we set the initial learning rate to $\tau_0 = 5\times 10^{-3}$ and decrease it according to the rule: $\tau_i = \frac{\tau_{i-1}}{1 + \text{decaying rate}\times i}$, where $i$ denotes the index of optimization iterations. The maximum number of iterations is constrained to be  $10^4$. 
For the latter, we set the initial learning rate $\gamma_0 = 10$ and the maximum number of iterations to $10^3$. At each iteration, the learning rate is adapted following the backtracking line search to ensure that the Armijo condition is satisfied. That is, 
\begin{equation}
E(\theta_i - \gamma \Delta \theta_i) \leq E(\theta_i) - c \gamma \|\Delta \theta_i\|^2,
\label{eq:armijo}
\end{equation} 
where $c=2\times 10^{-4}$. The step size $\gamma$ is initialized to $\gamma_0$ and iteratively halved until the condition is met. 
 If the optimal loss is zero (see \eqref{al}, an (absolute loss-based) early termination can be triggered for either optimizer when the absolute loss value drops below a user-defined tolerance $\tau_t$ when $\cE(u^*)=0$ or when a saturation condition is met for a user-defined tolerance.  
\item Second, we investigate the performance of expansive NN models with initial depth $L_0$. 
The NN models are first trained by NGF. An expansdion occurs when the loss saturates,
according to a loss saturation condition of the form \eqref{sl}. 
Specifically, once such a condition is met   
 we expand the NN model by inserting a new hidden layer before the output layer (\textsection~\ref{ssec:expand}). 
For newly introduced parameters after each expansion, we compare two initialization methods: (i) random initialization, where parameters are randomly generated; and (ii) gradient-aligned initialization, described in the previous section. To show the basic effects we
use the simplest version from \eqref{alignment} with  
{ $K = 20$} random instances of learnable parameters for the new hidden layer. The corresponding output layer parameters are then determined  such that the change in the NN model aligns with the (ideal Hilbert space-)energy gradient computed for the pre-expansion model at the current iteration stage. The instance yielding the least loss value is selected as the initial parameters for the new layer and output layer, see Section \ref{ssec:expand} for details. Obviously, the effect of such a strategy is increased by increasing $K$ or employing more elaborate minimizers (see Section \ref{ssec:expand}).

Each expansion increases network complexity and hence the computational cost of assembling the flow matrix. In favor of scalability   we only train the last two layers after each expansion by using the NGF optimizer. Once the loss gets saturated again, all the layers are trained using Adam to relax all learnable parameters. 
{ In particular, if the optimal loss value is unknown a priori (see \eqref{sl}), such as when solving numerical solutions to PDEs via energy minimization (see Section \ref{ssec:ellip}), when the Adam optimization saturates, we terminate the training process if the loss saturation condition   
\begin{equation}
|E^{(n_\text{exp})}-E^{(n_\text{exp}-1)}|\leq \tau_{t_d} 
\text{ or } 
|E^{(n_\text{exp})}-E^{(n_\text{exp}-1)}|/ |E^{(n_\text{exp}-1)}|\leq \tau_{t_r}. 
\end{equation}
is met.
If neither criterion is met, we }
expand the NN model and repeat the process until either the respective early termination condition is met or  the maximum number of iterations (set to $3\times 10^3$) is reached. 
Recall that when $\cE(u^*)=0$ the ``generalized loss'' at termination is a sharp error bound.
\end{enumerate}
One overarching message is to comply with the metric imposed by the underlying energy model.
Rather than working with a more convenient $L_2$-norm, we assemble the flow matrix with respect to the inner product of the ambient Hilbert space in order to be able to 
assess the achieved accuracy. One of the objectives of the experiments is to highlight related effects which are often ignored in the literature.
Again, since we are primarily interested in better understanding the principal mechansisms 
we are contend here with treating relatively simple toy models.

\subsection{Supervised Learning} 
We first consider a function approximation task with the target function 
$$y(x) = e^{\sin(k\pi x)}+x^3-x-1$$ 
defined on the interval $x\in [0,1]$, where $k$ is a parameter controlling the oscillatory behavior of the function. As $k$ increases, the function exhibits higher oscillations, making the approximation task more challenging. 
  Although existing expressivity results ensure that accurate approximations exist,    it is known that
capturing highly oscillatory features via optimization is difficult. 
Specifically, we consider two test cases: (I) $k = 5$ and (II) $k=10$.

In both cases, we use supervised learning to train fully-connected NNs with a fixed number $L$ of layers and a uniform width $N$ in all hidden layers. A dataset consisting of $M_{train}= 201$ randomly selected samples from the domain and range is used for training, making sure that the sampling density exceeds the Nyquist rate. The same set of samples is also used for evaluating integrals involved in the optimization process by Monte-Carlo (MC) method. 
Additionally, a test dataset of the size $M_{test}= 301$ is used to assess the accuracy of the obtained NN models. 
We set a loss tolerance { $\tau_{t_a} = 10^{-5}$} to allow for a possible early termination of the optimization process.

When NGF is applied, we solve $\left(G_{\nN,\Theta}+\lambda I\right)\Delta \theta = -\nabla_\theta E(\theta)$ for the gradient vector. The regularization parameter $\lambda$ is selected based on the maximum diagonal entry $G_{mm}$ of $G_{\nN,\Theta}$. In particular, $\lambda$ is a piecewise function of the form (see the discussion in Section \ref{ssec:issues})
\begin{equation}
\lambda(G_{mm}) = \lambda_1\cdot 1_{G_{mm}<x_1} + \sum_{j=2}^6 \lambda_j \cdot 1_{x_{j-1}\leq G_{mm}<x_{j}} + \lambda_7\cdot 1_{G_{mm}>x_6}
\label{eq:lambda}
\end{equation}
where $1_A$ is the indicator function (1 if condition $A$ holds, 0 otherwise), $x_j = 10^{j-1}$ and $\lambda_j = 5\times 10^{j-6}$, for $j=1,\ldots, 6$. 

\paragraph{Case I: $k = 5$} 

To compare natural gradient flow with Adam as optimizers, we vary the depth of the NN model from $L=2$ to $L=4$, while keeping the width of all the hidden layers fixed at $N=15$.

The history of training loss when Adam is used is shown in Figure~\ref{fig:fak5_fixedL_adam}(a), from which we observe that the loss decreases to the tolerance for $L=3$ and $L=4$, while it fails to reach the tolerance after $10^4$ iterations when $L=2$. 
The trained NNs are evaluated on the test set, which are shown in Figure~\ref{fig:fak5_fixedL_adam}(b). We observe that the NN approximations converge to the target function as $L$ increases from $L=2$ to $L=3$ and $4$.
\begin{figure}[!htb]
  \centering
  	\begin{subfigure}[htbp]{0.475\textwidth}
        \begin{tikzpicture}
        \node[inner sep=0pt] at (0,0)
        {\includegraphics[width=\textwidth]{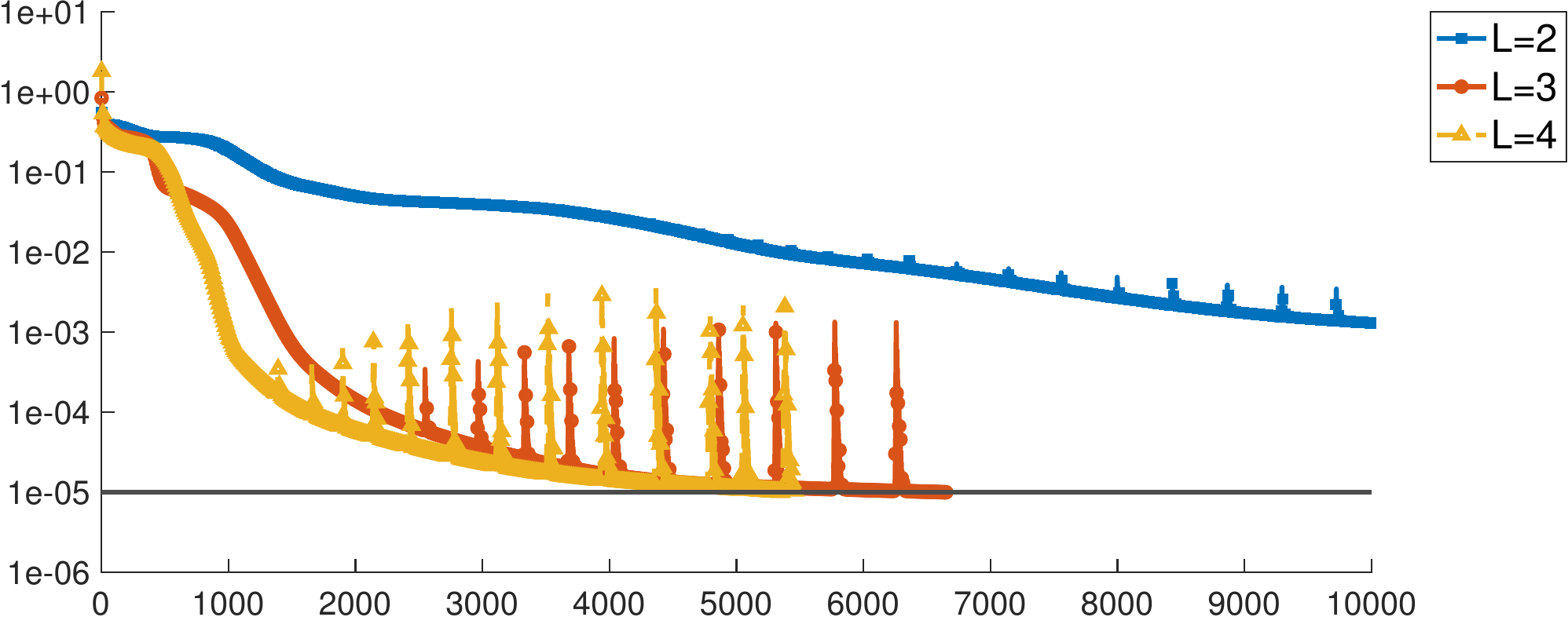}};
        \node[align=center,fill=white] at (3.3,-1.0) {(a)};
        \end{tikzpicture}
    \end{subfigure}
    \quad 
    \begin{subfigure}[htbp]{0.475\textwidth}
        \begin{tikzpicture}
        \node[inner sep=0pt] at (0,0)
        {\includegraphics[width=\textwidth]{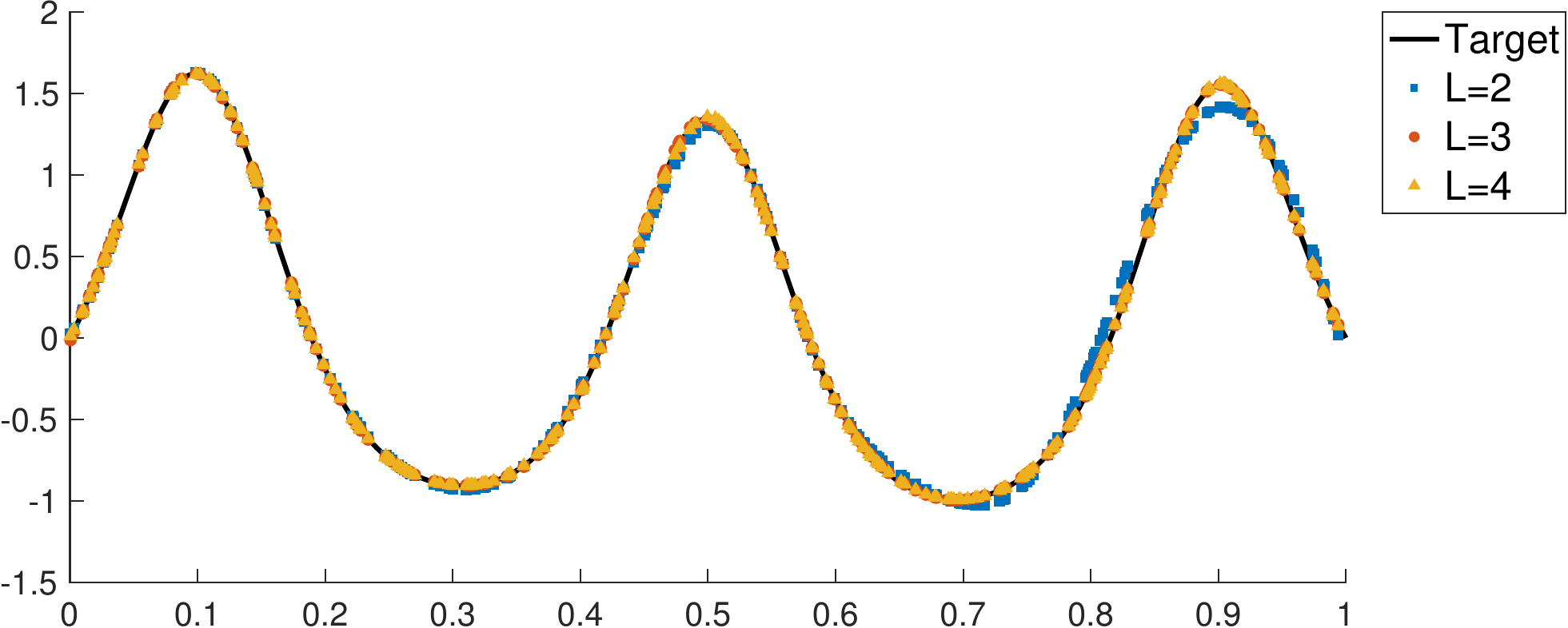}};
        \node[align=center,fill=white] at (3.3,-1.0) {(b)};
        \end{tikzpicture}
    \end{subfigure}
	\caption{SL-Case I ($k=5$): Neural network of fixed depth $L$ trained using Adam optimizer. (a) History of training loss, where the dark line indicates the prescribed loss tolerance; (b) Final approximation on the test set.}
\label{fig:fak5_fixedL_adam}
\end{figure}

When switching to NGF for optimization, the training process converges and reaches the tolerance for all tested depths, as shown in Figure~\ref{fig:fak5_fixedL_ngf}(a). It is observed that NGF converges much faster than Adam, with the loss decaying to the tolerance level within several dozen iterations. Thus, the optimization method strongly impacts the ability of actually exploit  expressivity in a regime of {\em under-parametrization} which is highly desirable in a model reduction context.

The resulting NN approximations are shown in Figure~\ref{fig:fak5_fixedL_ngf}(b), which closely match the target function on the test set. 
\begin{figure}[!htb]
  \centering
  	\begin{subfigure}[htbp]{0.475\textwidth}
        \begin{tikzpicture}
        \node[inner sep=0pt] at (0,0)
        {\includegraphics[width=\textwidth]{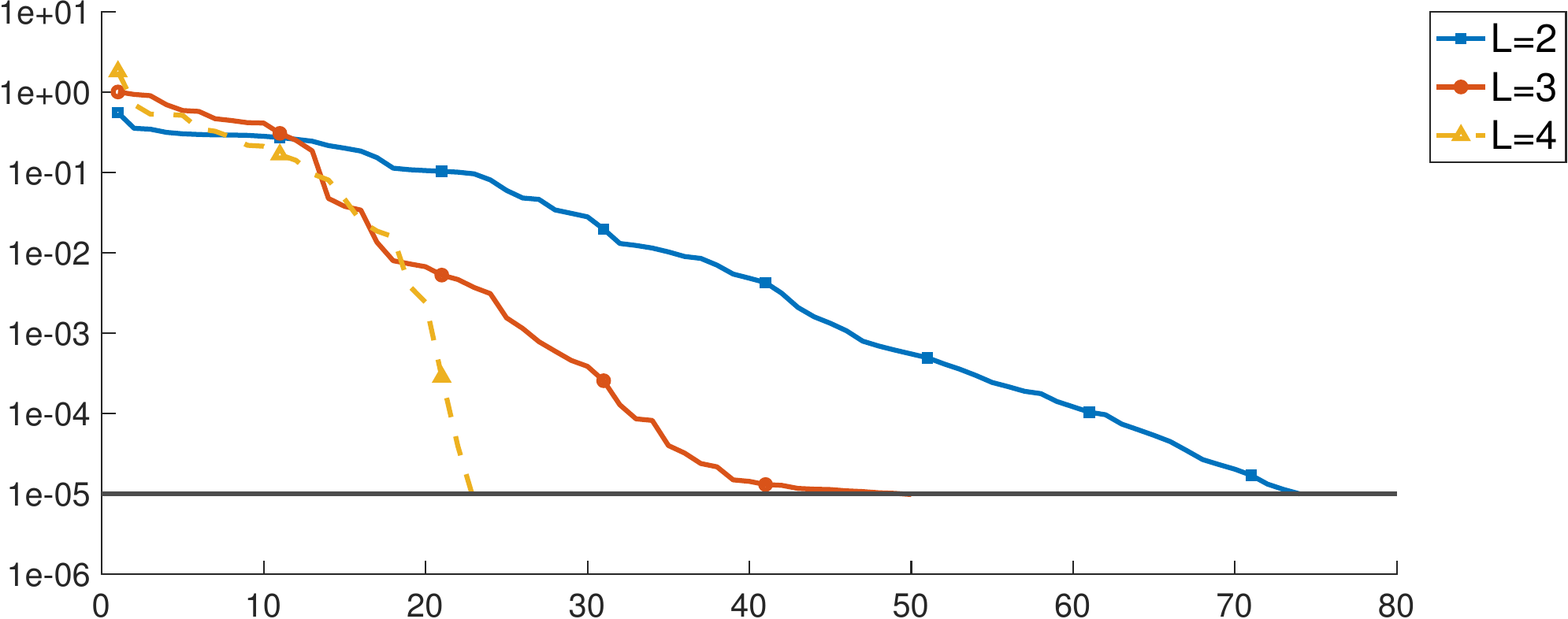}};
        \node[align=center,fill=white] at (3.3,-1.0) {(a)};
        \end{tikzpicture}
    \end{subfigure}
    \quad 
    \begin{subfigure}[htbp]{0.475\textwidth}
        \begin{tikzpicture}
        \node[inner sep=0pt] at (0,0)
        {\includegraphics[width=\textwidth]{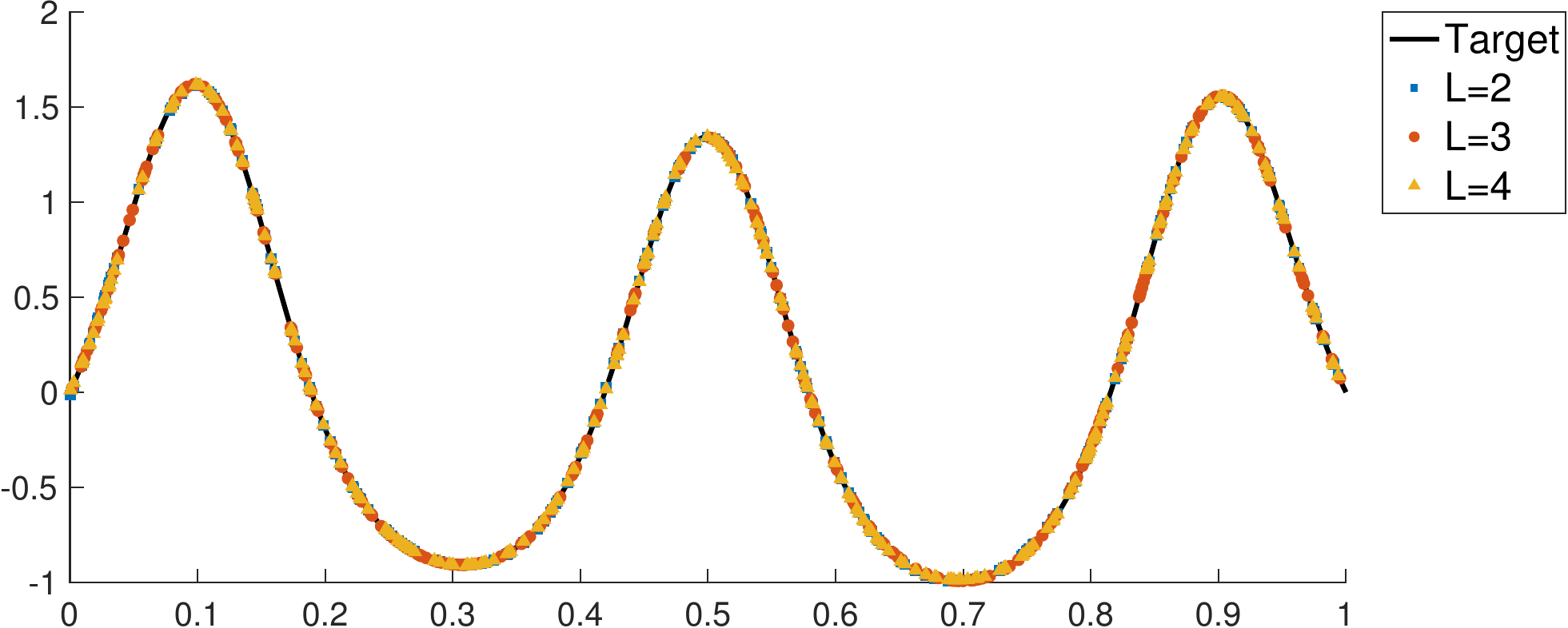}};
        \node[align=center,fill=white] at (3.3,-1.0) {(b)};
        \end{tikzpicture}
    \end{subfigure}
	\caption{SL-Case I ($k=5$): Neural network of fixed depth $L$ trained using NGF optimizer. (a) History of training loss, where the dark line indicates the prescribed loss tolerance; (b) Final approximation on the test set.}
\label{fig:fak5_fixedL_ngf}
\end{figure}

We summarize the total number of iterations, final training loss values, and the $L_2$ approximation errors on the test set in Table~\ref{tab:fak5_fixedL_adam_ngf}. Overall, once the training loss decays to the tolerance level, the obtained NN models  achieve comparable accuracy. 
\begin{table}[!ht]
\centering
\begin{tabular}{c|ccc|ccc}
\hline
\multirow{2}{*}{$L$} & \multicolumn{3}{c|}{Adam} & \multicolumn{3}{c}{NGF} \\ \cline{2-7}
{}    & $N_\text{ite}$ & Loss & $L_2$ error  & $N_\text{ite}$ & Loss & $L_2$ error\\ \hline
$2$ & $10^4$ & $1.30\times 10^{-3}$ & $5.58\times 10^{-2}$ & $74$ & $9.99\times 10^{-6}$ & $3.91\times 10^{-3}$\\
$3$ & $6652$ & $1.00\times 10^{-5}$ & $4.84\times 10^{-3}$ & $50$  & $9.81\times 10^{-6}$ & $5.00\times 10^{-3}$\\
$4$ & $5468$ & $9.96\times 10^{-6}$ & $4.80\times 10^{-3}$ & $23$  & $8.10\times 10^{-6}$ & $3.32\times 10^{-3}$\\
\hline
\end{tabular}
\caption{NN performance for SL-Case I ($k=5$): Comparison of training convergence ($N_\text{ite}$ iterations), training loss, and approximation accuracy ($L_2$) across depths $L$ using Adam/NGF.}
\label{tab:fak5_fixedL_adam_ngf}

\end{table}

Finally, we employ NGF for optimization to evaluate the performance of the {\em expansive} NN with an initial depth of $L_0=2$. However, in this case, the loss fails to saturate before reaching the prescribed tolerance threshold, preventing the expansion from being triggered. So, corresponding tests are only relevant for higher oscillation frequency
considered below.

\paragraph{Case II: $k = 10$}

For the target function with $k=10$, we repeat the same tests. A larger $k$ implies the higher oscillation frequency in the target function, which makes training an NN model to reach the  $10^{-5}$ loss tolerance more difficult. 


When Adam is used, the loss history for training the NN model of depth $L$ is shown in Figure~\ref{fig:fak10_fixedL_adam}(a). We observe that the loss  decays slowly when $L=2$ - although it decreases faster for $L=3$ and $L=4$  - yet fails to meet the tolerance threshold even after $10^4$ iterations. The corresponding NN approximations on the test set are presented in Figure~\ref{fig:fak10_fixedL_adam}(b), indicating accuracy is poor for $L=2$,   improves as 
$L$ increases, approaching the target function more closely.
\begin{figure}[!htb]
  \centering
  	\begin{subfigure}[htbp]{0.475\textwidth}
        \begin{tikzpicture}
        \node[inner sep=0pt] at (0,0)
        {\includegraphics[width=\textwidth]{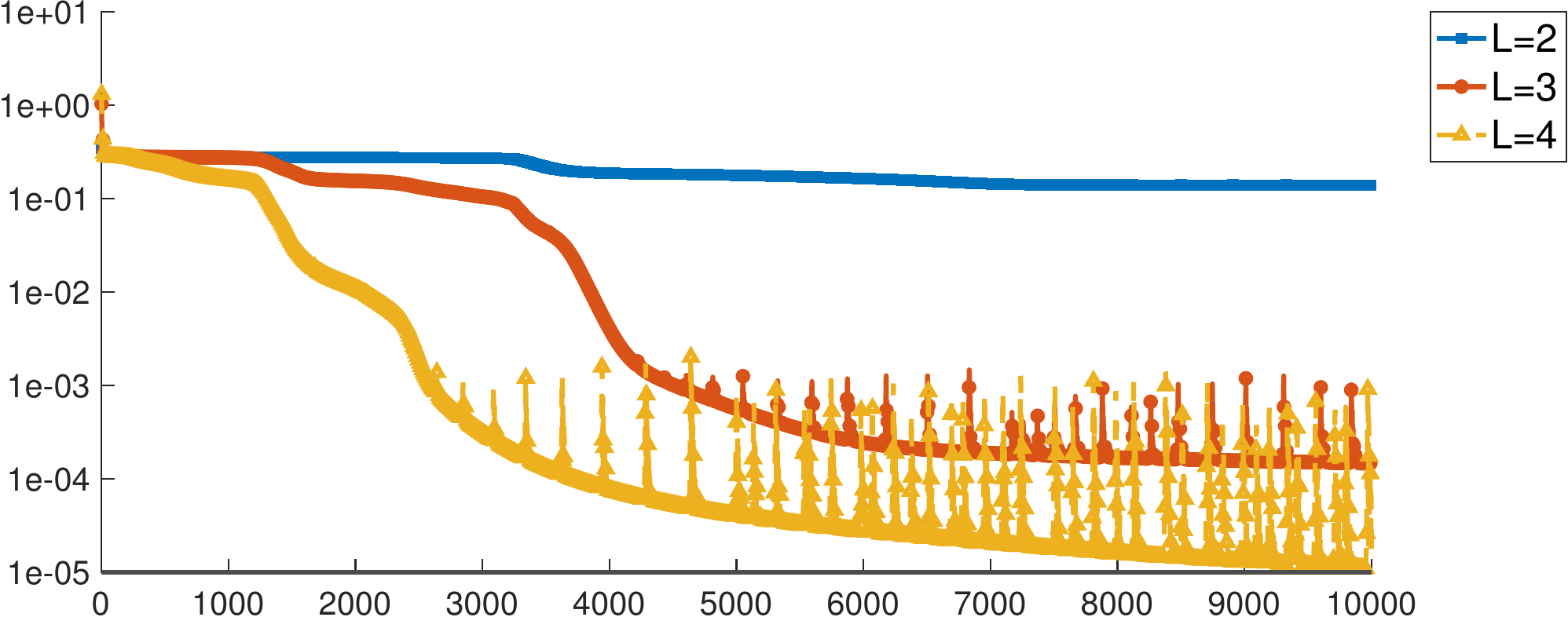}};
        \node[align=center,fill=white] at (3.3,-1.0) {(a)};
        \end{tikzpicture}
    \end{subfigure}
    \quad
    \begin{subfigure}[htbp]{0.475\textwidth}
        \begin{tikzpicture}
        \node[inner sep=0pt] at (0,0)
        {\includegraphics[width=\textwidth]{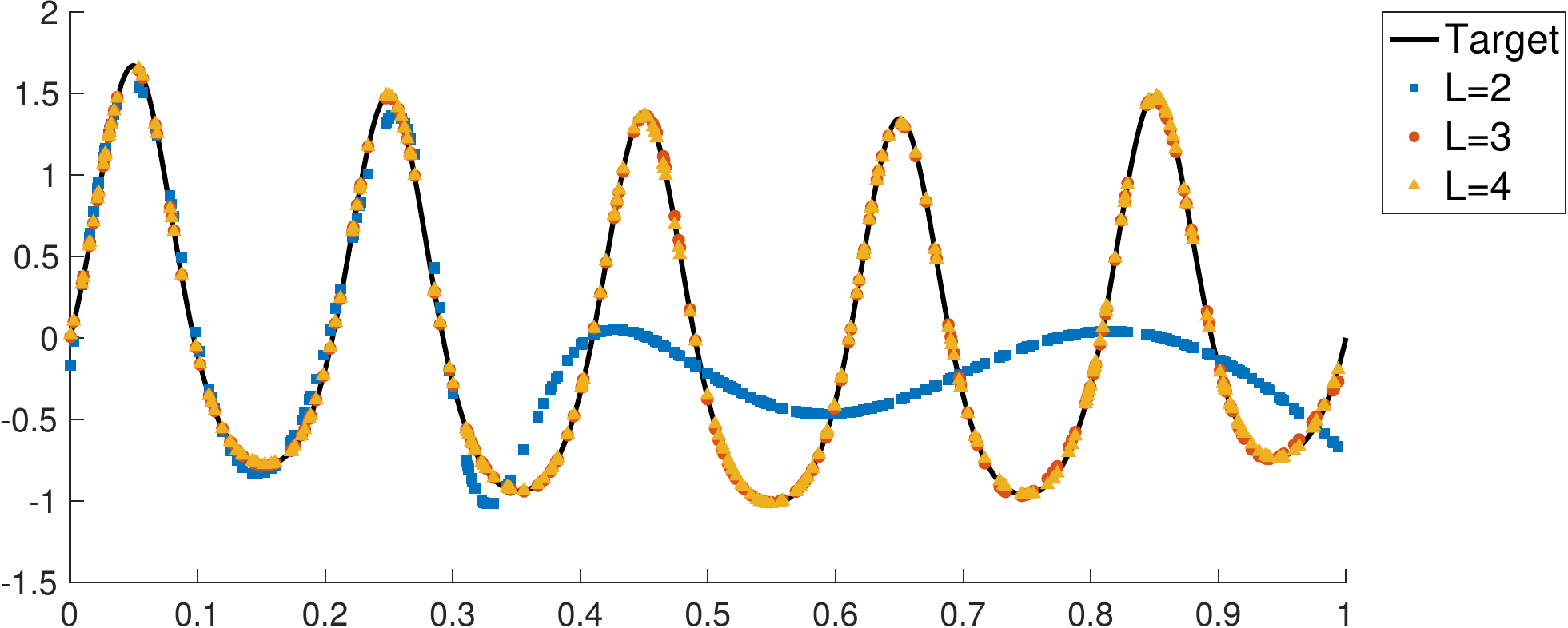}};
        \node[align=center,fill=white] at (3.3,-1.0) {(b)};
        \end{tikzpicture}
    \end{subfigure}
	\caption{SL-Case II ($k=10$): Neural network of fixed depth $L$ trained using Adam optimizer. (a) History of training loss, where the dark line indicates the prescribed loss tolerance; (b) Final approximation on the test set.}
\label{fig:fak10_fixedL_adam}
\end{figure}


When the NGF optimizer is applied, the training process converges more rapidly for both $L=3$ and $L=4$. Figure~\ref{fig:fak10_fixedL_ngf}(a) shows the corresponding loss history, while Figure~\ref{fig:fak10_fixedL_ngf}(b) displays the NN approximations over the test set. Specifically, the NGF optimizer successfully meets the target loss tolerance for both $L=3$ and $L=4$ within the maximum number of iterations, with the resulting approximations closely matching the target function. 
\begin{figure}[!htb]
  \centering
  	\begin{subfigure}[htbp]{0.475\textwidth}
        \begin{tikzpicture}
        \node[inner sep=0pt] at (0,0)
        {\includegraphics[width=\textwidth]{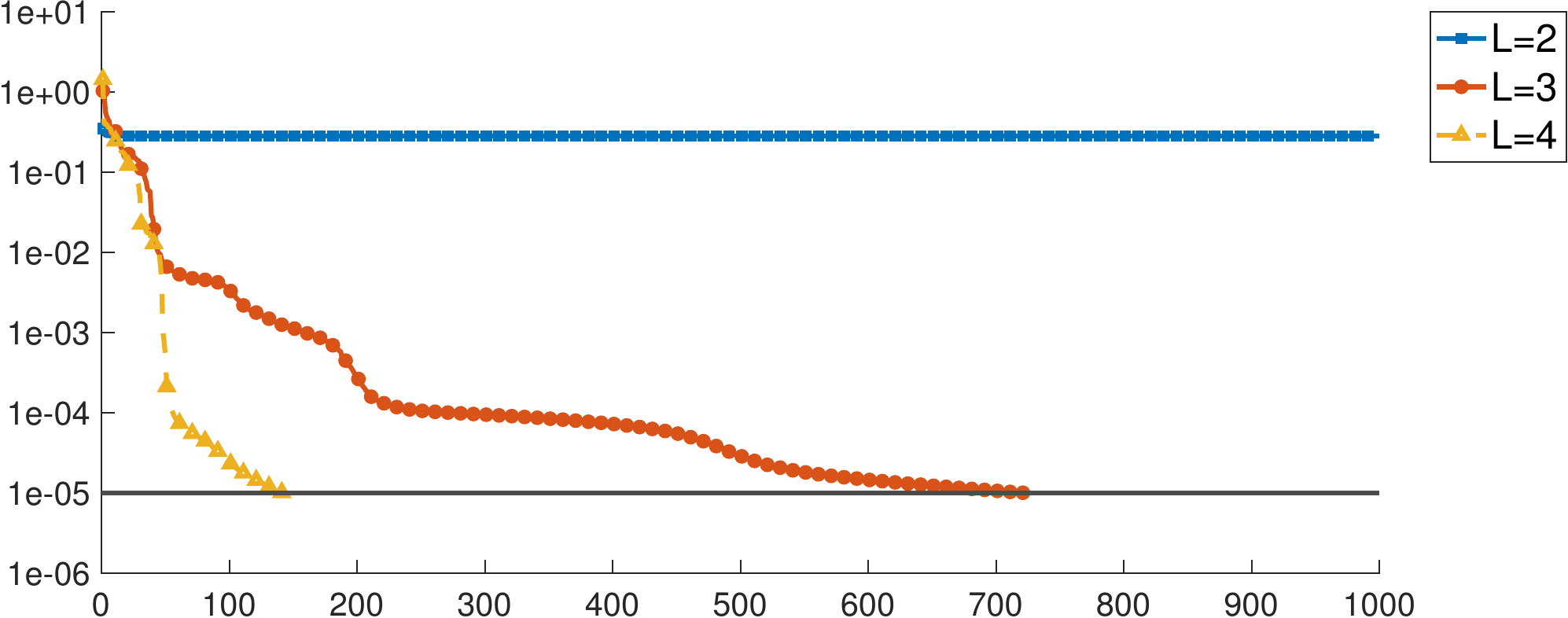}};
        \node[align=center,fill=white] at (3.3,-1.0) {(a)};
        \end{tikzpicture}
    \end{subfigure}
    \quad
    \begin{subfigure}[htbp]{0.475\textwidth}
        \begin{tikzpicture}
        \node[inner sep=0pt] at (0,0)
        {\includegraphics[width=\textwidth]{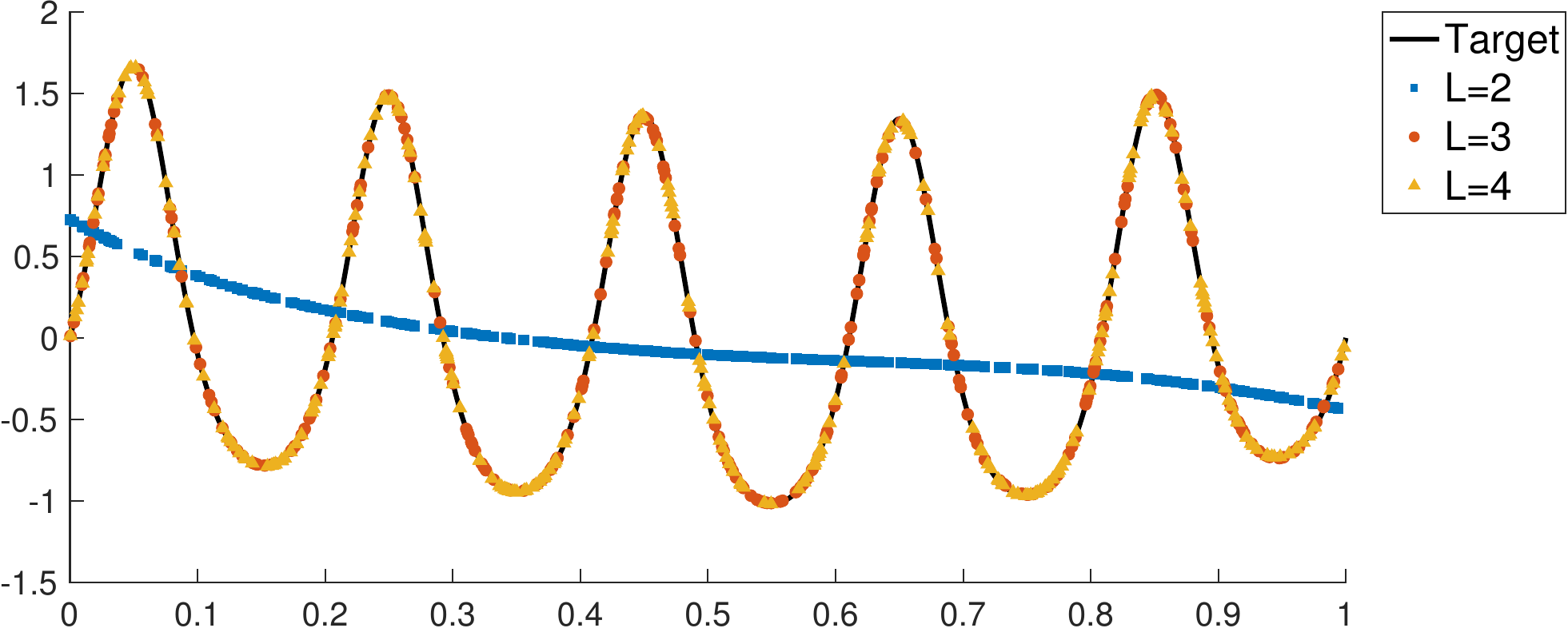}};
        \node[align=center,fill=white] at (3.3,-1.0) {(b)};
        \end{tikzpicture}
    \end{subfigure}
	\caption{SL-Case II ($k=10$): Neural network of fixed depth $L$ trained using NGF optimizer. (a) History of training loss, where the dark line indicates the prescribed loss tolerance; (b) Final approximation on the test set. }
\label{fig:fak10_fixedL_ngf}
\end{figure}


The final loss values and approximation errors are compared in Table~\ref{tab:fak10_fixedL_adam_ngf}.
\begin{table}[!ht]
\centering
\begin{tabular}{c|ccc|ccc}
\hline
\multirow{2}{*}{$L$} & \multicolumn{3}{c|}{Adam} & \multicolumn{3}{c}{NGF} \\ \cline{2-7}
{}    & $N_\text{ite}$ & Loss & $L_2$ error  & $N_\text{ite}$ & Loss & $L_2$ error\\ \hline
$2$ & $10^4$ & $1.38\times 10^{-1}$ & $6.18\times 10^{-1}$ & $10^3$ & $2.80\times 10^{-1}$ & $8.09\times 10^{-1}$\\
$3$ & $10^4$ & $1.47\times 10^{-4}$ & $1.69\times 10^{-2}$ & $725$  & $9.97\times 10^{-6}$ & $6.37\times 10^{-3}$\\
$4$ & $10^4$ & $1.28\times 10^{-5}$ & $6.03\times 10^{-3}$ & $142$  & $9.88\times 10^{-6}$ & $5.37\times 10^{-3}$\\
\hline
\end{tabular}
\caption{NN performance for SL-Case I ($k=10$): Comparison of training convergence ($N_\text{ite}$ iterations), training loss, and approximation accuracy ($L_2$ norm) across depths $L$ using Adam/NGF.}
\label{tab:fak10_fixedL_adam_ngf}
\end{table}

Next, we employ NGF for optimization and assess the performance of the expansive NN with an initial depth of $L_0=2$. One new aspect of interest is the initialization after an expansion step. 
Specifically, after each expansion, we use either random initialization or gradient-alignment initialization for new parameters. As explained earlier, to keep scalability for larger scale problems feasible, we train only the last two layers by using the NGF optimizer. Once the loss gets saturated again, all the layers will be trained using Adam to relax all learnable parameters. After it gets saturated again, we will expand the NN model and repeat the process till either the loss tolerance is met or the maximum number of expansions is reached. 
{ In our loss saturation criterion, $\tau_{s_a} = 10^{-7}$ and $\tau_{s_r} = 5\times 10^{-3}$ are used for NGF optimizer and $\tau_{s_a} = 10^{-8}$ and $\tau_{s_r} = 5\times 10^{-4}$ for Adam optimizer.}

Figure \ref{fig:fak10_expNGF_lastTwo} shows the training loss histories for two expansive NN models, each expanded once at iteration 13 under the two initialization methods. While using random initialization narrowly misses the convergence tolerance by the maximum iteration limit, the gradient-aligned initialization successfully meets the prescribed threshold. Table~\ref{tab:fak10_expNGF_lastTwo} summarizes the iteration counts, final training loss values, and corresponding $L_2$ errors on the test set. The progressive improvements of the NN approximations is further illustrated in Figure~\ref{fig:fak10_expNGF_lastTwo_align}, which displays the model outputs over the test set both before the expansion and at training termination. 
\begin{figure}[!htb]
  \centering
  	\begin{subfigure}[htbp]{0.475\textwidth}
        \begin{tikzpicture}
        \node[inner sep=0pt] at (0,0)
        {\includegraphics[width=\textwidth]{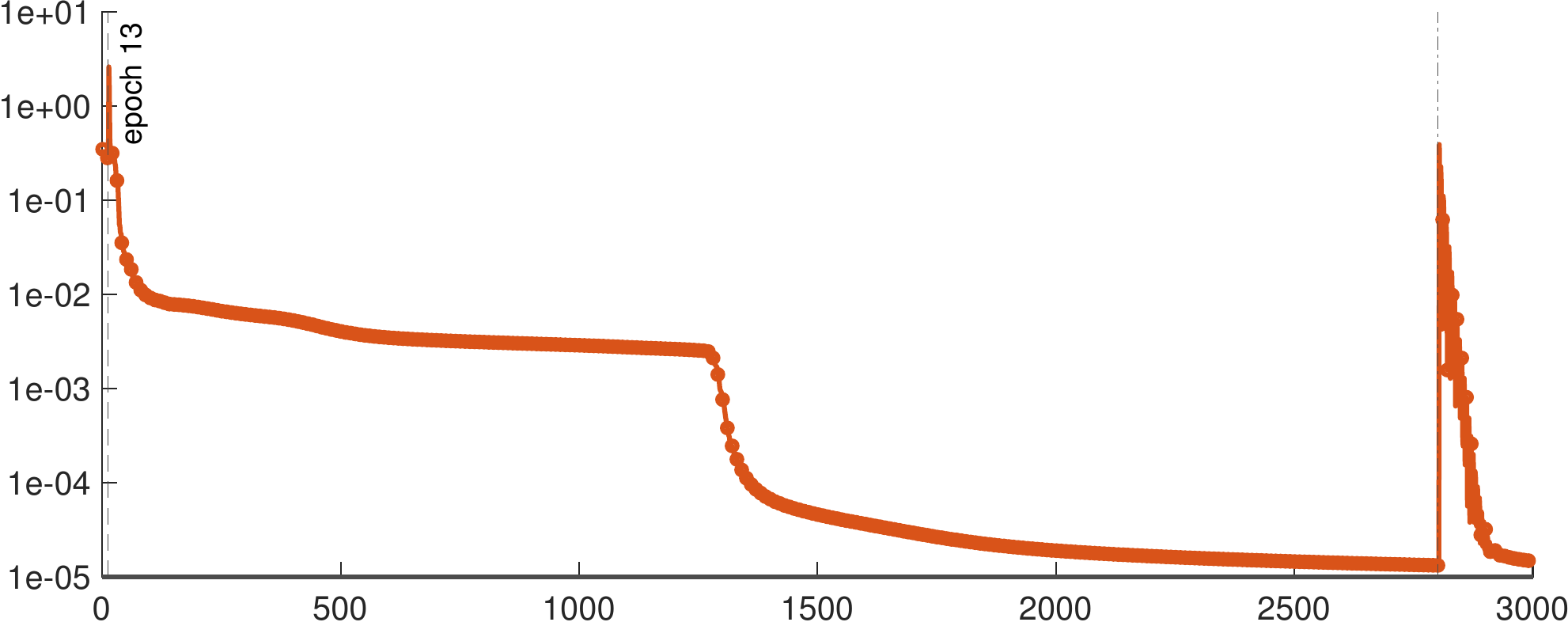}};
        \node[align=center,fill=white] at (2.3,1.1) {(a)};
        \end{tikzpicture}
    \end{subfigure}
    \quad
    \begin{subfigure}[htbp]{0.475\textwidth}
        \begin{tikzpicture}
        \node[inner sep=0pt] at (0,0)
        {\includegraphics[width=\textwidth]{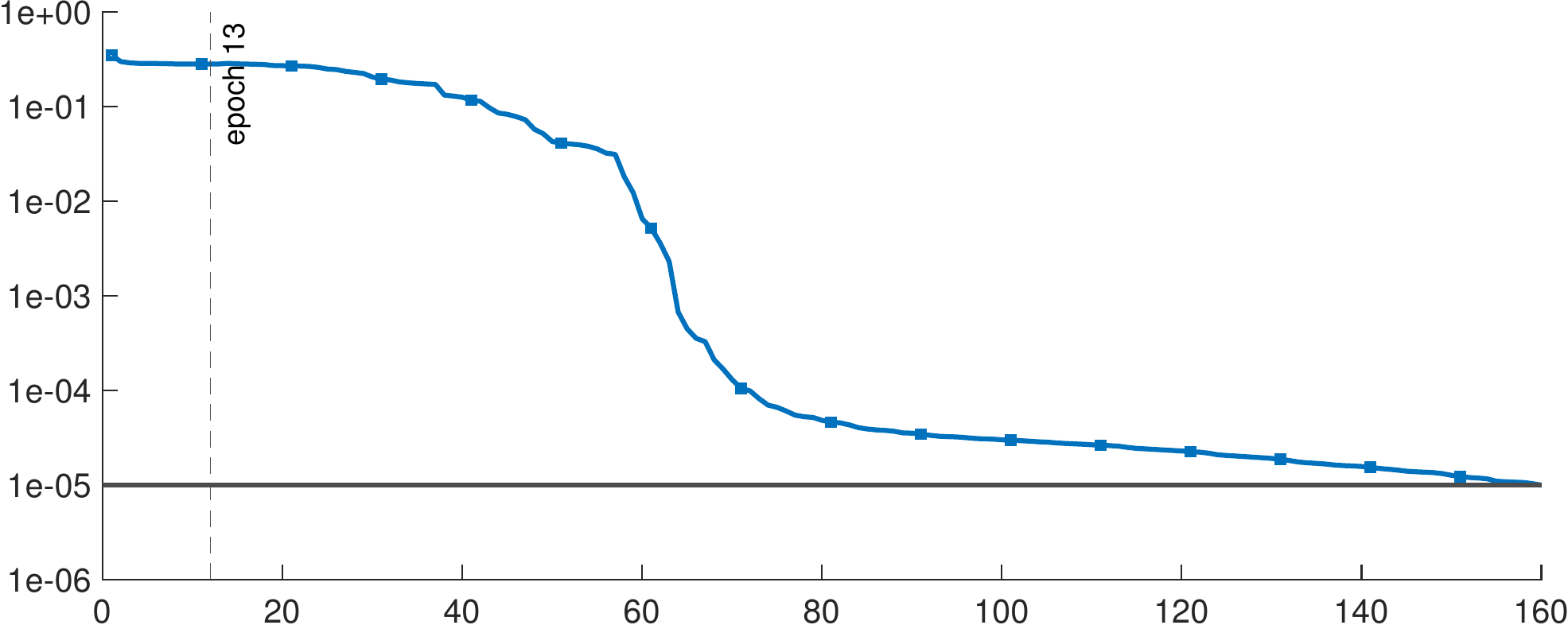}};
        \node[align=center,fill=white] at (2.3,1.1) {(b)};
        \end{tikzpicture}
    \end{subfigure}
	\caption{SL-Case II ($k=10$): expansive NN with initial depth $L_0=2$, comparing (a) random vs. (b) gradient-aligned initialization. Vertical lines indicate expansion events (dashed) and NGF-to-Adam transitions (dash-dot). The dark horizontal line indicates the prescribed loss tolerance.}
\label{fig:fak10_expNGF_lastTwo}
\end{figure}

\begin{table}[!ht]
\centering
\begin{tabular}{c|cccc}
\hline
\multirow{2}{*}{Initalization} & \multicolumn{4}{c}{Expansive NN with $L_0=2$}  \\ \cline{2-5}
{}       & $N_\text{ite}$ & $N_\text{exp}$ & Loss & $L_2$ error                    \\ \hline
Random   & $3\times 10^3$ & 1 & $1.47\times 10^{-5}$ & $6.50\times 10^{-3}$ \\
Gradient-aligned & $160$ & 1 & $9.99\times 10^{-6}$ & $7.56\times 10^{-3}$ \\
\hline
\end{tabular}
\caption{Expansive NN performance for SL-Case II ($k=10$): Comparison of training convergence ($N_\text{ite}$ iterations), number of expansions ($N_\text{exp}$), training loss, and approximation accuracy ($L_2$ norm).}
\label{tab:fak10_expNGF_lastTwo}
\end{table}

\begin{figure}[!htb]
  \centering
  	\begin{subfigure}[htbp]{0.425\textwidth}
        \begin{tikzpicture}
        \node[inner sep=0pt] at (0,0)
        {\includegraphics[width=\textwidth]{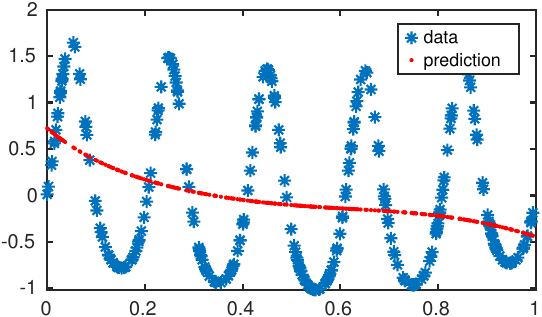}};
        \node[align=center,fill=white] at (3.6,-1.3) {(a)};
        \end{tikzpicture}
    \end{subfigure}
    \qquad
    \begin{subfigure}[htbp]{0.425\textwidth}
        \begin{tikzpicture}
        \node[inner sep=0pt] at (0,0)
        {\includegraphics[width=\textwidth]{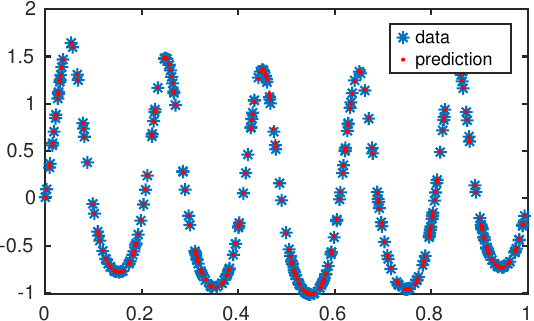}};
        \node[align=center,fill=white] at (3.6,-1.3) {(b)};
        \end{tikzpicture}
    \end{subfigure}
	\caption{SL-Case II ($k=10$): (a) Approximation of the expansive NN model on the test set before the first expansion; (b) Final approximation on the test set.}
\label{fig:fak10_expNGF_lastTwo_align}
\end{figure}

\subsection{Unsupervised Learning - PDE Solution} 

To demonstrate the feasibility of the proposed concept, we consider the following one-dimensional diffusion problem: 
\begin{equation}
-u_{xx} = g \text{ in } \Omega\coloneqq [0,1] \text{ and } u(0)=u(1) = 0.
\label{eq:pde}
\end{equation}
Although the problem is simple, when approximating $u(x)$ by a NN model, $\nN(x; \theta)$, the seek of learnable parameters in the NN model via optimization shares the same difficulty as high-dimensional problems. Moreover, controlling accuracy in the energy norm, here the $H^1$-norm (i.e. controlling fluxes) is typically much more important than
approximating in $L_2$. It is known that its solution minimizes the Ritz energy 
$
\mathcal{E}(v):= \frac 12 \|v\|^2_\H- g(v)
$
with $\H= H_0^1(\Omega)$, and the corresponding gradient of $\mathcal{E}(v)$ satisfies 
$
\langle \nabla \mathcal{E}(v),h\rangle_\H= (v_x, h_x) - g(h)
$, $\forall h\in \H$, see Section \ref{ssec:modelred}, \eqref{enlift}. 
In the following test, we set $g$ in \eqref{eq:pde} such that the exact solution of the problem is 
$$u = e^{\sin(k\pi x)}+x^3-x-1,$$
to compare with previous experiments.
Same as the previous supervised learning task, we test two cases: (I) $k = 5$ and (II) $k=10$. However, in this experiment, the loss function is replaced with the Ritz energy, which makes the learning task unsurpervised.

In the NN model, we use a hard constraint to impose the homogeneous Dirichlet boundary condition - a mask function $y=-4(x^2-x)$ is multiplied to the NN output. To calculate the integrals in the loss function, we use a composite trapezoidal rule with $M_\text{train} = 401$ uniformly distributed nodes from the domain. The accuracy of the trained NN model is evaluated on a test set consisting of $M_\text{test} = 301$ uniformly spaced points.  
We set the maximum number of iterations to be $10^4$ for Adam and $10^3$ for NGF. As the optimal loss value is unknown beforehand for general problems, early termination is not considered when training NNs of a fixed depth.  

When NGF is applied, we solve $\left(G_{T,\Theta}+\lambda I\right)\Delta \theta = -\nabla_\theta E(\theta)$ for the weight increment $\Dt$. The regularization parameter $\lambda$ is chosen by the rule \eqref{eq:lambda}.

\paragraph{Case I: $k=5$}

We first compare the performance of  NGF to Adam for training NNs. The NN architectures have depths ranging from $L=2$ to $L=4$, with a fixed width of $N=15$ in each hidden layer. 

When Adam is applied, the training loss converges in all cases as shown in Figure~\ref{fig:poisson1dk5_fixedL_adam}(a), where the dark line represents the exact value of the Ritz energy (value of -110.90). Particularly, NNs of the depth $L=3$ and $L=4$ are able to converge more rapidly than the $L=2$ case. 
The trained NNs are evaluated on the test set, with results displayed in Figure~\ref{fig:poisson1dk5_fixedL_adam}(b). In all cases, the NN approximations are close to the exact solution.
\begin{figure}[!htb]
  \centering
  	\begin{subfigure}[htbp]{0.475\textwidth}
        \begin{tikzpicture}
        \node[inner sep=0pt] at (0,0)
        {\includegraphics[width=\textwidth]{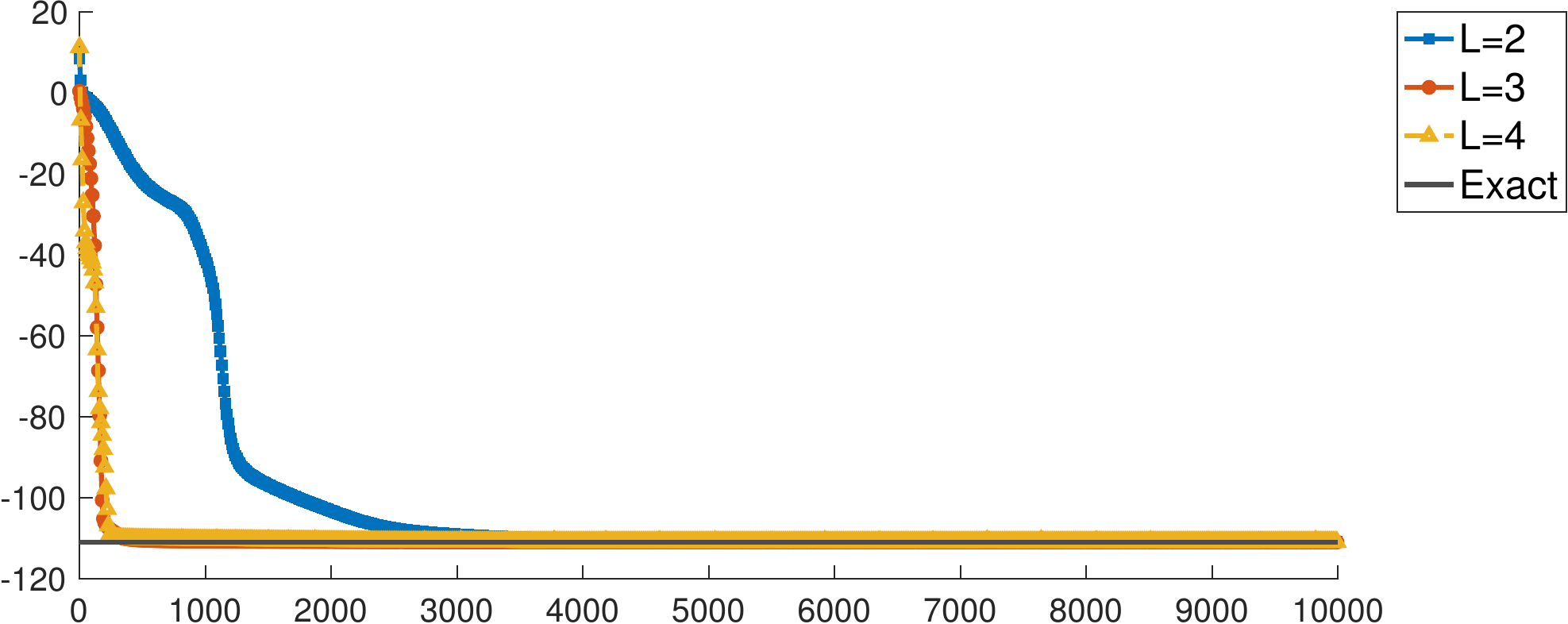}};
        \node[align=center,fill=white] at (3.3,-1.0) {(a)};
        \end{tikzpicture}
    \end{subfigure}
    \quad
    \begin{subfigure}[htbp]{0.475\textwidth}
        \begin{tikzpicture}
        \node[inner sep=0pt] at (0,0)
        {\includegraphics[width=\textwidth]{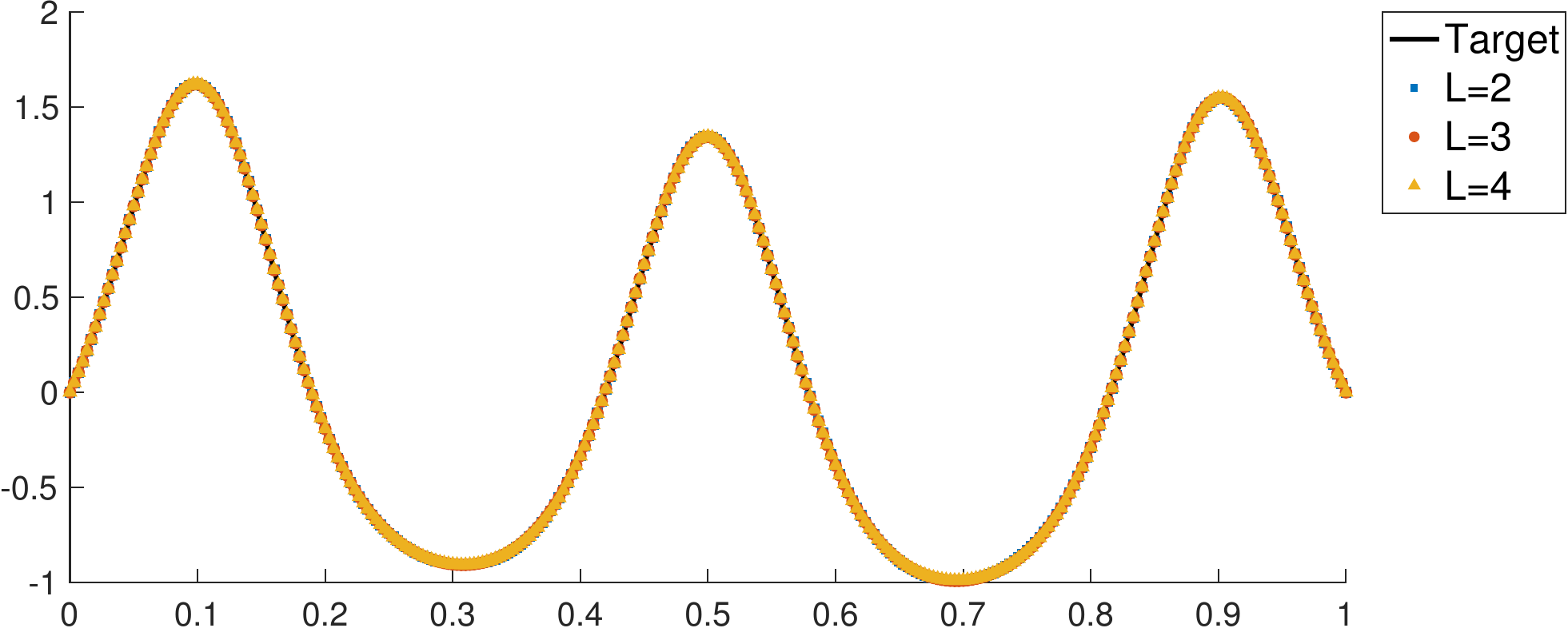}};
        \node[align=center,fill=white] at (3.3,-1.0) {(b)};
        \end{tikzpicture}
    \end{subfigure}
	\caption{PDE-Case I ($k=5$): (a) Training loss history for NNs of fixed depths $L$, along with the exact Ritz energy (dark line); (b) Final NN approximations on the test set. All models were trained using the Adam optimizer.}
\label{fig:poisson1dk5_fixedL_adam}
\end{figure}

When the optimizer is switched to NGF, the training loss is shown in Figure~\ref{fig:poisson1dk5_fixedL_ngf}(a). In all cases, NGF  converges rapidly. The resulting NN approximations on the test set are shown in Figure~\ref{fig:poisson1dk5_fixedL_ngf}(b), which closely match the target function on the test set. 
\begin{figure}[!htb]
  \centering
  	\begin{subfigure}[htbp]{0.475\textwidth}
        \begin{tikzpicture}
        \node[inner sep=0pt] at (0,0)
        {\includegraphics[width=\textwidth]{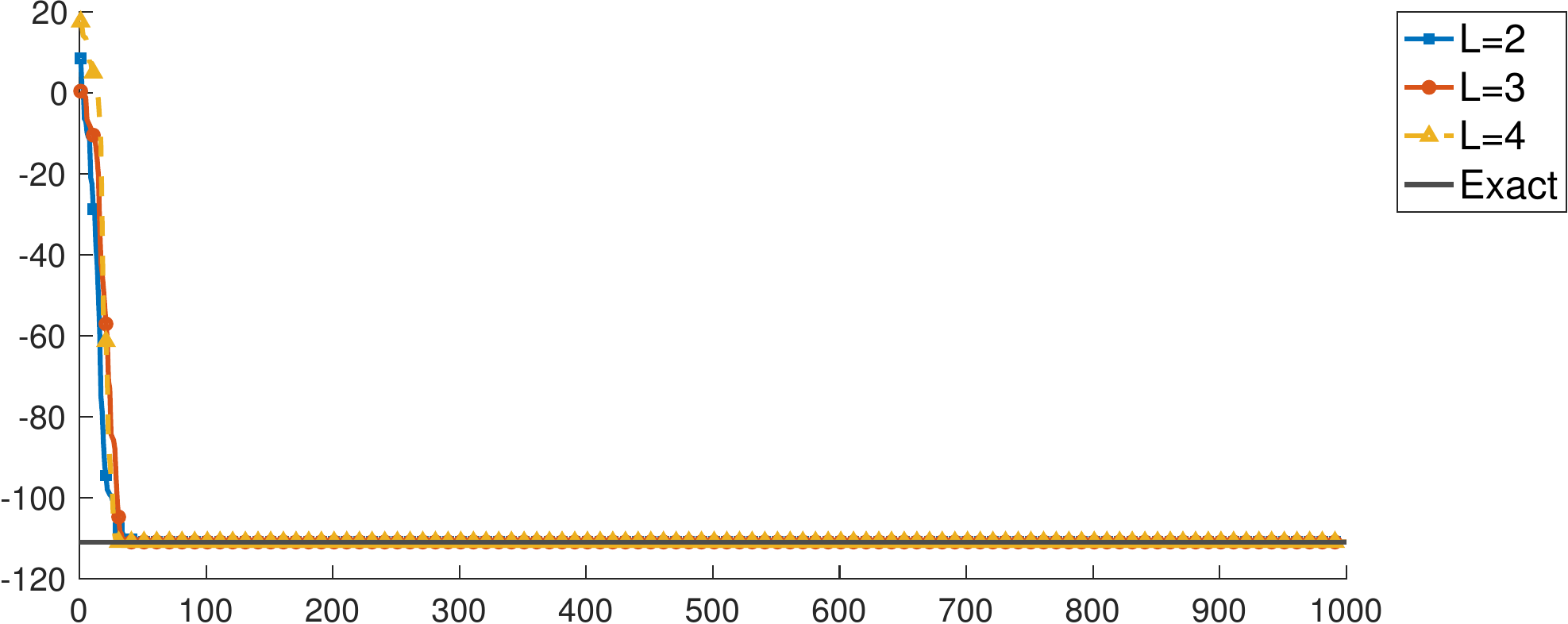}};
        \node[align=center,fill=white] at (3.3,-1.0) {(a)};
        \end{tikzpicture}
    \end{subfigure}
    \quad
    \begin{subfigure}[htbp]{0.475\textwidth}
        \begin{tikzpicture}
        \node[inner sep=0pt] at (0,0)
        {\includegraphics[width=\textwidth]{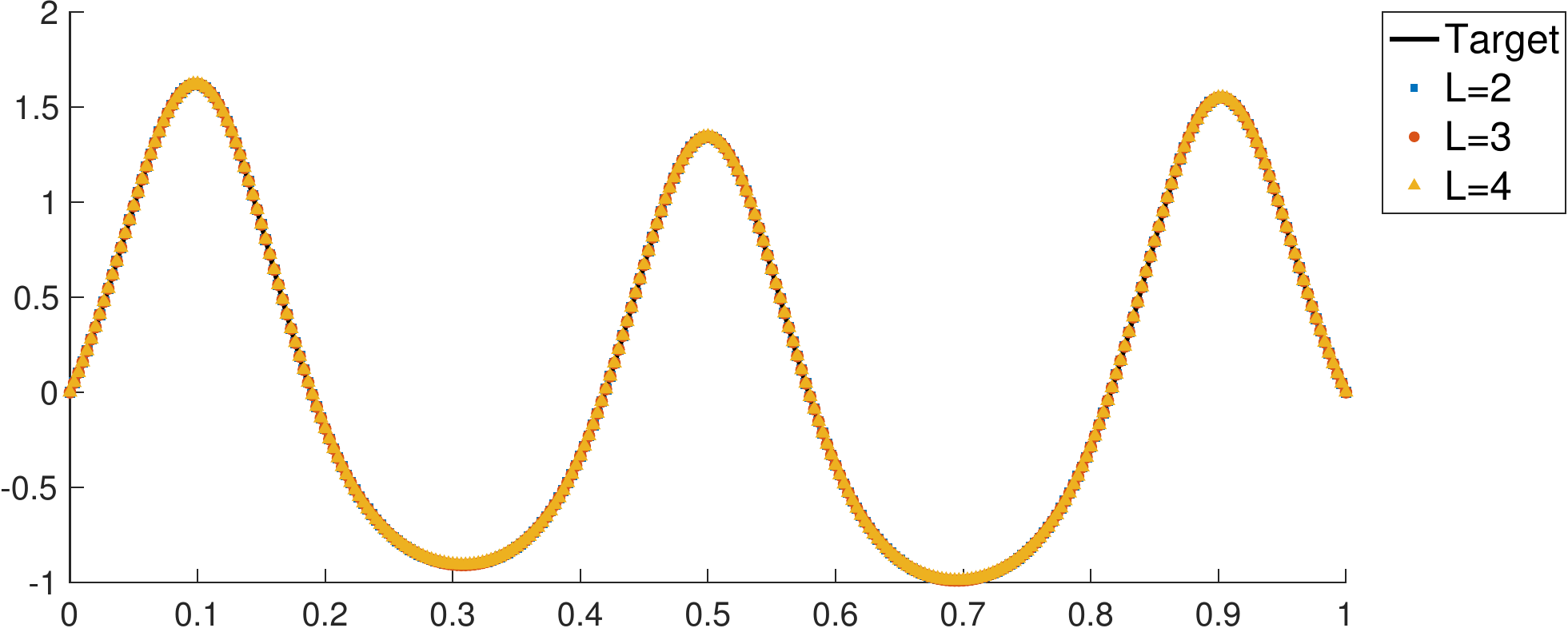}};
        \node[align=center,fill=white] at (3.3,-1.0) {(b)};
        \end{tikzpicture}
    \end{subfigure}
	\caption{PDE-Case I ($k=5$): (a) Training loss history for NNs of fixed depths $L$, along with the exact Ritz energy (dark line); (b) Final NN approximations on the test set. All models were trained using the NGF optimizer.}
\label{fig:poisson1dk5_fixedL_ngf}
\end{figure}

For both optimizers, the number of iterations during training, final loss values, and approximation errors of the trained NN models, measured in $L_2$ and $H_1$ norms, are listed in Table~\ref{tab:poissonk5_ngf_gd}. Here, the relevant accuracy measure is the $H^1$-norm.
In all cases NGF achieves a better accuracy.
\begin{table}[!ht]
\centering
\begin{tabular}{c|cccc|cccc}
\hline
\multirow{2}{*}{$L$} & \multicolumn{4}{c|}{Adam} & \multicolumn{4}{c}{NGF} \\ \cline{2-9}
{}    & $N_\text{ite}$ & Loss & $L_2$ error & $H_1$ error  & $N_\text{ite}$ & Loss & $L_2$ error & $H_1$ error\\ \hline
$2$ & $10^4$ & $-110.90$ & $1.67\times 10^{-3}$ & $9.74\times 10^{-2}$ & $10^3$  & $-110.91$ & $3.61\times 10^{-4}$ & $1.91\times 10^{-2}$\\
$3$ & $10^4$ & $-110.90$ & $1.33\times 10^{-3}$ & $5.26\times 10^{-2}$ & $10^3$  & $-110.91$ & $4.13\times 10^{-4}$ & $2.20\times 10^{-2}$\\
$4$ & $10^4$ & $-110.91$ & $3.63\times 10^{-4}$ & $2.36\times 10^{-2}$ & $10^3$  & $-110.91$ & $2.78\times 10^{-4}$ & $1.72\times 10^{-2}$\\
\hline
\end{tabular}
\caption{NN performance for PDE-Case I ($k=5$): Comparison of training convergence ($N_\text{ite}$ iterations), training loss, and approximation accuracy ($L_2$ and $H_1$ norms) across depths $L$ using Adam/NGF.}
\label{tab:poissonk5_ngf_gd}
\end{table}

Next, we use NGF as the optimizer and train the expansive NN model with an initial depth of $L_0=2$. 
{ The same training strategy is used. In our loss saturation criterion, $\tau_{s_a} = 10^{-8}$ and $\tau_{s_r} = 5\times 10^{-5}$ are used for NGF optimizer and $\tau_{s_a} = 10^{-9}$ and $\tau_{s_r} = 5\times 10^{-6}$ for Adam optimizer. Since the optimal energy is unknown a priori, we use the relative loss-based early termination condition, in which $\tau_{t_d} = 5\times 10^{-3}$ and $\tau_{t_r} = 10^{-6}$. }

Figure~\ref{fig:poisson1dk5_expNGF_lastTwo} compares loss histories between (a) random initialization and (b) gradient-aligned initialization cases. In both cases, the NN model was expanded once after 61 iterations, and the optimization processes were successfully terminated. The results are summarized in Table~\ref{tab:poisson1dk5_exp}. The NN approximate solutions on the test set are shown in Figure~\ref{fig:poisson1dk5_expNGF_lastTwo_align}, where (a) displays the NN solution before the first expansion and (b) displays the final approximation, along with the exact data. 
\begin{figure}[!htb]
  \centering
  	\begin{subfigure}[htbp]{0.475\textwidth}
        \begin{tikzpicture}
        \node[inner sep=0pt] at (0,0)
        {\includegraphics[width=\textwidth]{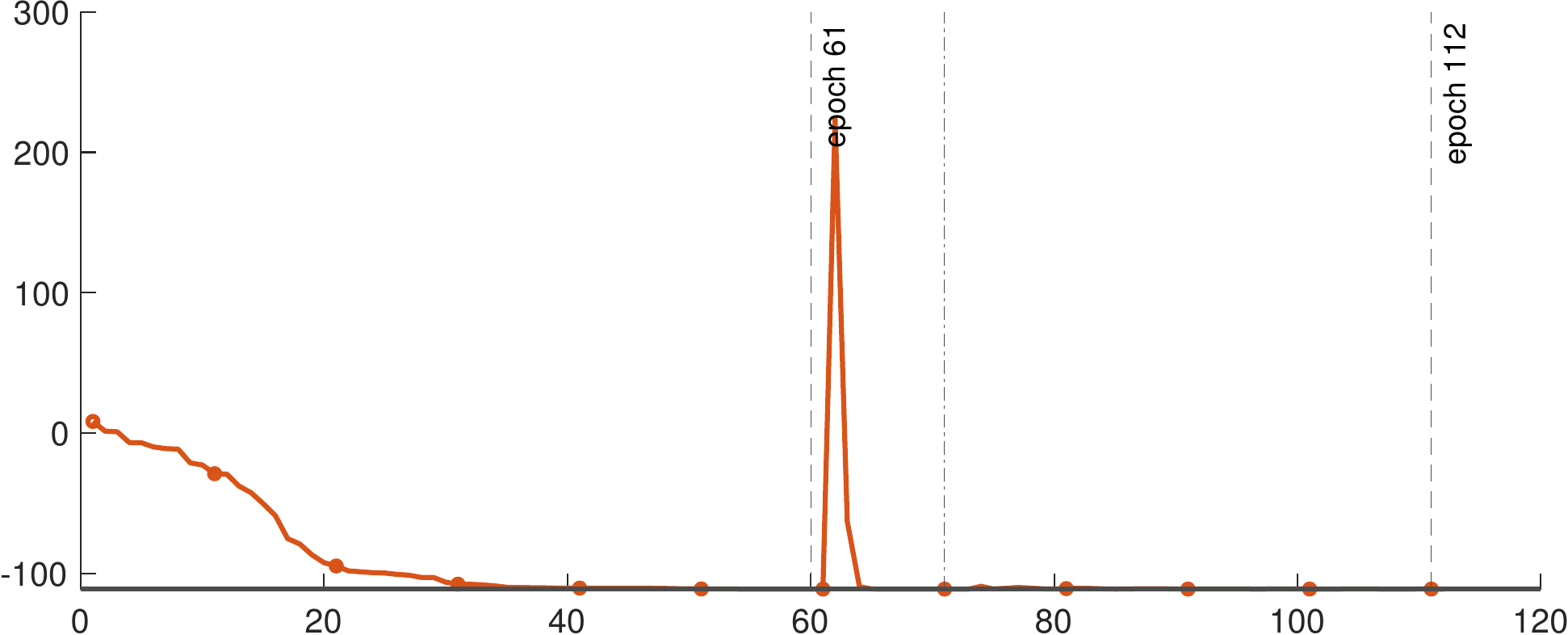}};
        \node[align=center,fill=white] at (2.2,1.1) {(a)};
        \end{tikzpicture}
    \end{subfigure}
    \quad
    \begin{subfigure}[htbp]{0.475\textwidth}
        \begin{tikzpicture}
        \node[inner sep=0pt] at (0,0)
        {\includegraphics[width=\textwidth]{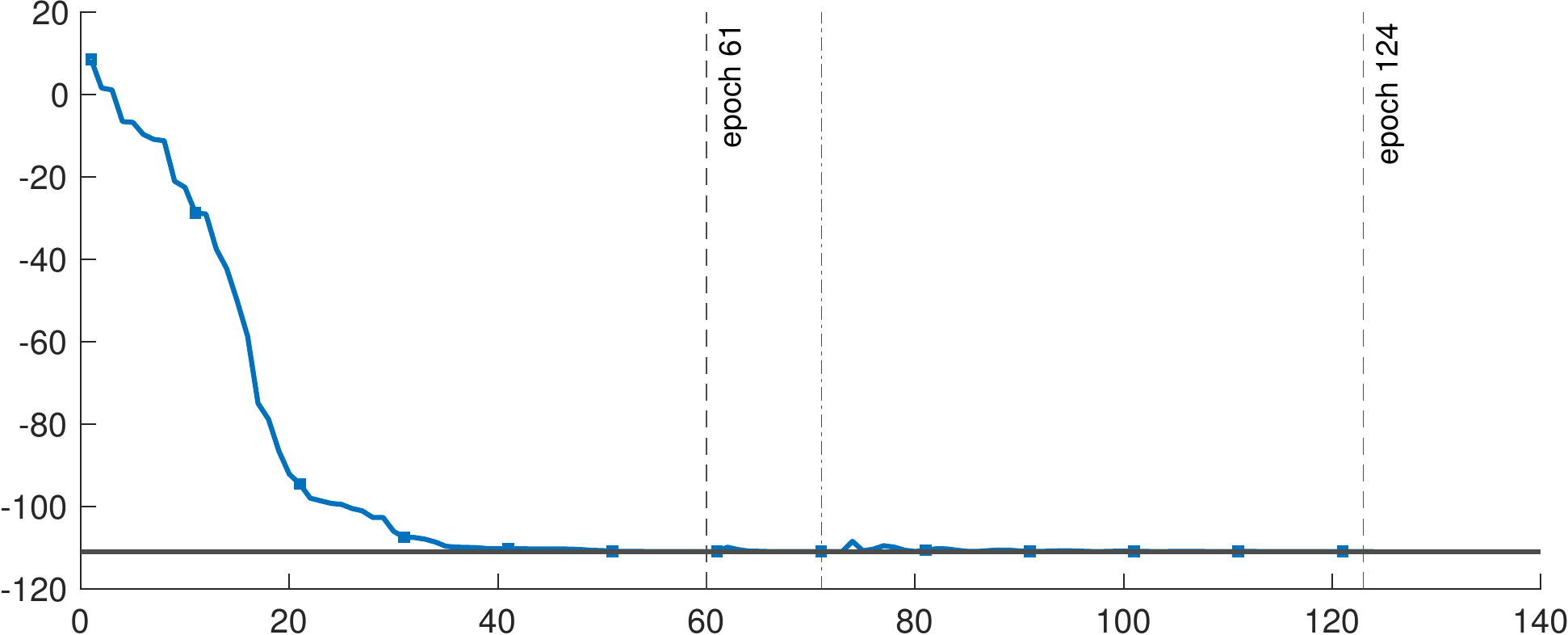}};
        \node[align=center,fill=white] at (2.2,1.1) {(b)};
        \end{tikzpicture}
    \end{subfigure}
	\caption{PDE-Case I ($k=5$): expansive NN with initial depth $L_0=2$, comparing (a) random vs. (b) gradient-aligned initialization. Vertical lines indicate expansion events (dashed) and NGF-to-Adam transitions (dash-dot). The dark horizontal line indicates the exact Ritz energy.}
\label{fig:poisson1dk5_expNGF_lastTwo}
\end{figure}

\begin{table}[!ht]
\centering
\begin{tabular}{c|ccccc}
\hline
\multirow{2}{*}{Initalization} & \multicolumn{5}{c}{Expansive NN with $L_0=2$} \\ \cline{2-6} 
{}    & $N_\text{ite}$ & $N_{\text{exp}}$ & Loss & $L_2$ error & $H_1$ error \\ \hline
Random           & $112$ & 1 & $-110.89$ & $8.60\times 10^{-3}$ & $5.70\times 10^{-2}$ \\
Gradient-aligned & $124$ & 1 &$-110.90$ & $3.01\times 10^{-3}$ & $7.90\times 10^{-2}$ \\
\hline
\end{tabular}
\caption{Expansive NN performance for PDE-Case I ($k=5$): Comparison of training convergence ($N_\text{ite}$ iterations), number of expansions ($N_\text{exp}$), training loss, and approximation accuracy ($L_2$ and $H_1$ norms).}
\label{tab:poisson1dk5_exp}
\end{table}

\begin{figure}[!htb]
  \centering
  	\begin{subfigure}[htbp]{0.425\textwidth}
        \begin{tikzpicture}
        \node[inner sep=0pt] at (0,0)
        {\includegraphics[width=\textwidth]{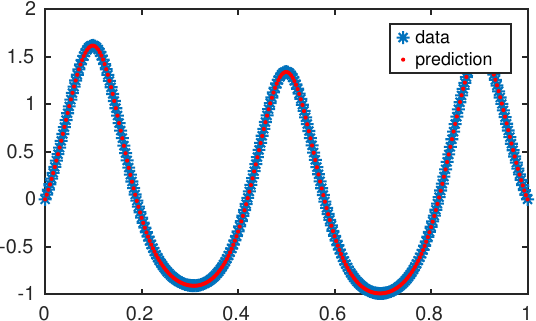}};
        \node[align=center,fill=white] at (3.6,-1.3) {(a)};
        \end{tikzpicture}
    \end{subfigure}
    \qquad
    \begin{subfigure}[htbp]{0.425\textwidth}
        \begin{tikzpicture}
        \node[inner sep=0pt] at (0,0)
        {\includegraphics[width=\textwidth]{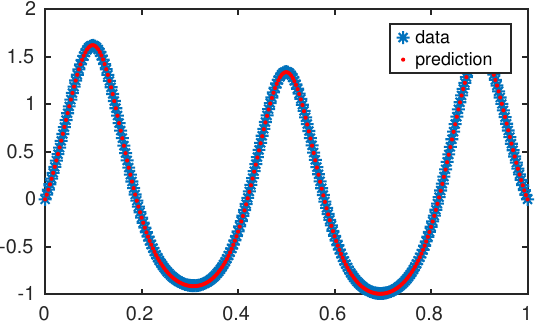}};
        \node[align=center,fill=white] at (3.6,-1.3) {(b)};
        \end{tikzpicture}
    \end{subfigure}
	\caption{PDE-Case I ($k=5$): (a) Approximation of the expansive NN model on the test set before the first expansion; (b) Final approximation on the test set.}
\label{fig:poisson1dk5_expNGF_lastTwo_align}
\end{figure}

\paragraph{Case II: $k = 10$}

For the target function with $k=10$, we repeat the same tests. Note that the exact energy value is -392.29 in this case.

The history of training loss using Adam is shown in Figure~\ref{fig:poisson1dk10_adam_fixedL}(a) along with the exact energy. From which, we observe that the training loss converges to the exact energy value in all the cases. Particularly, NNs of the depth $L=3$ and $L=4$ are able to quickly converge. 
The trained NN models are evaluated on the test set, which are shown in Figure~\ref{fig:poisson1dk10_adam_fixedL}(b). In all the cases, the NN approximations are close to the exact solution.
\begin{figure}[!htb]
  \centering
  	\begin{subfigure}[htbp]{0.475\textwidth}
        \begin{tikzpicture}
        \node[inner sep=0pt] at (0,0)
        {\includegraphics[width=\textwidth]{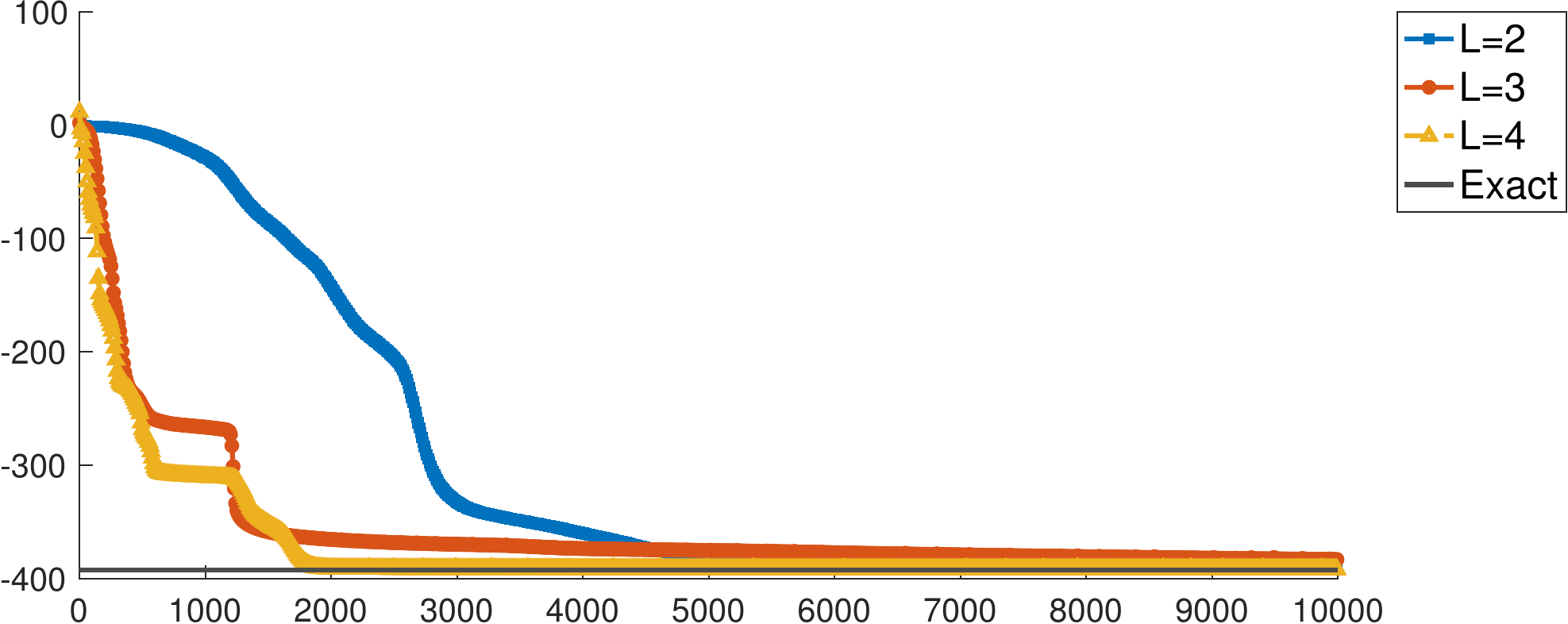}};
        \node[align=center,fill=white] at (3.3,-1.0) {(a)};
        \end{tikzpicture}
    \end{subfigure}
    \quad
    \begin{subfigure}[htbp]{0.475\textwidth}
        \begin{tikzpicture}
        \node[inner sep=0pt] at (0,0)
        {\includegraphics[width=\textwidth]{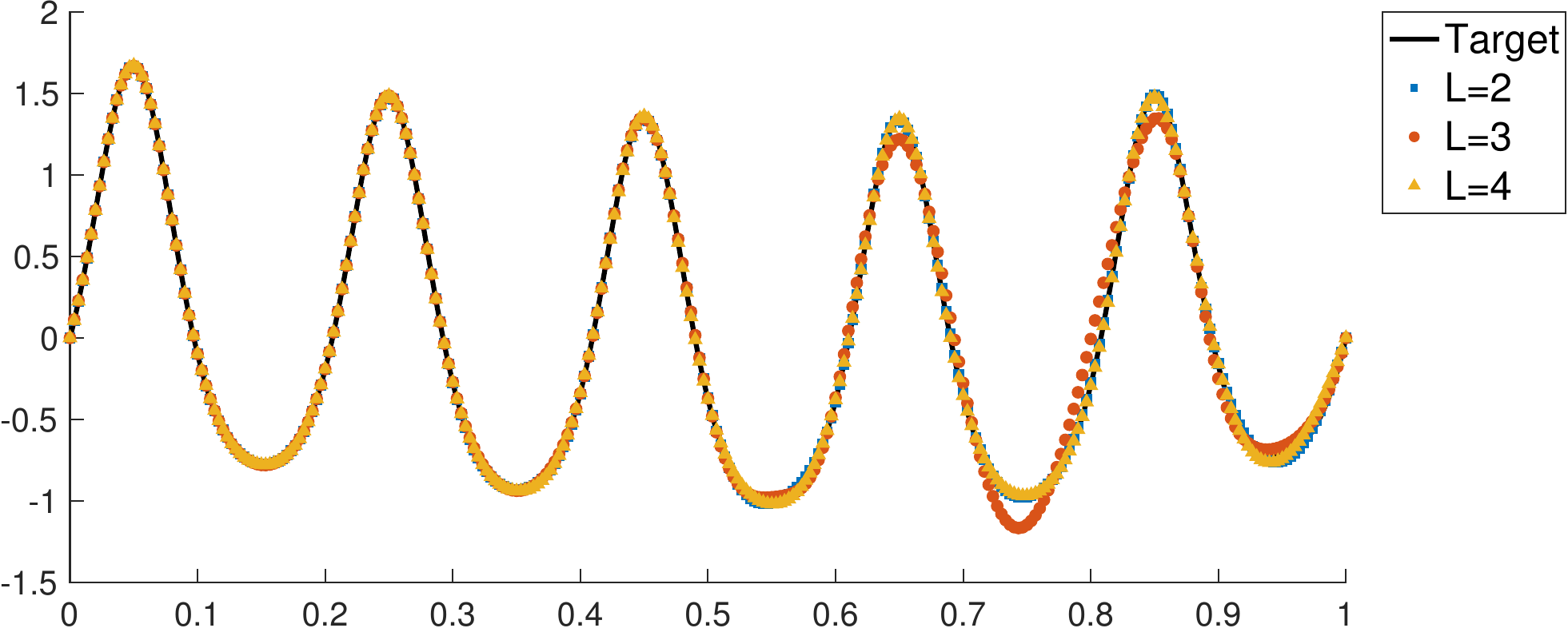}};
        \node[align=center,fill=white] at (3.3,-1.0) {(b)};
        \end{tikzpicture}
    \end{subfigure}
	\caption{PDE-Case II ($k=10$): (a) Training loss history for NNs of fixed depths $L$, along with the exact Ritz energy (dark line); (b) Final NN approximations on the test set. All models were trained using the Adam optimizer.}
\label{fig:poisson1dk10_adam_fixedL}
\end{figure}

Then, we change the optimizer to NGF, the training loss is shown in Figure~\ref{fig:possion1dk10_ngf_fixedL}(a) for NN with different depth. It is observed that, for $L=3$ and $L=4$, NGF rapidly converges to the exact energy. The resulting NN approximations over the test set are shown in Figure~\ref{fig:possion1dk10_ngf_fixedL}(b). It is seen that the NN models with $L=3$ and $L=4$ yield accurate approximations, while $L=2$ results in unsatisfying approximation. 

\begin{figure}[!htb]
  \centering
  	\begin{subfigure}[htbp]{0.475\textwidth}
        \begin{tikzpicture}
        \node[inner sep=0pt] at (0,0)
        {\includegraphics[width=\textwidth]{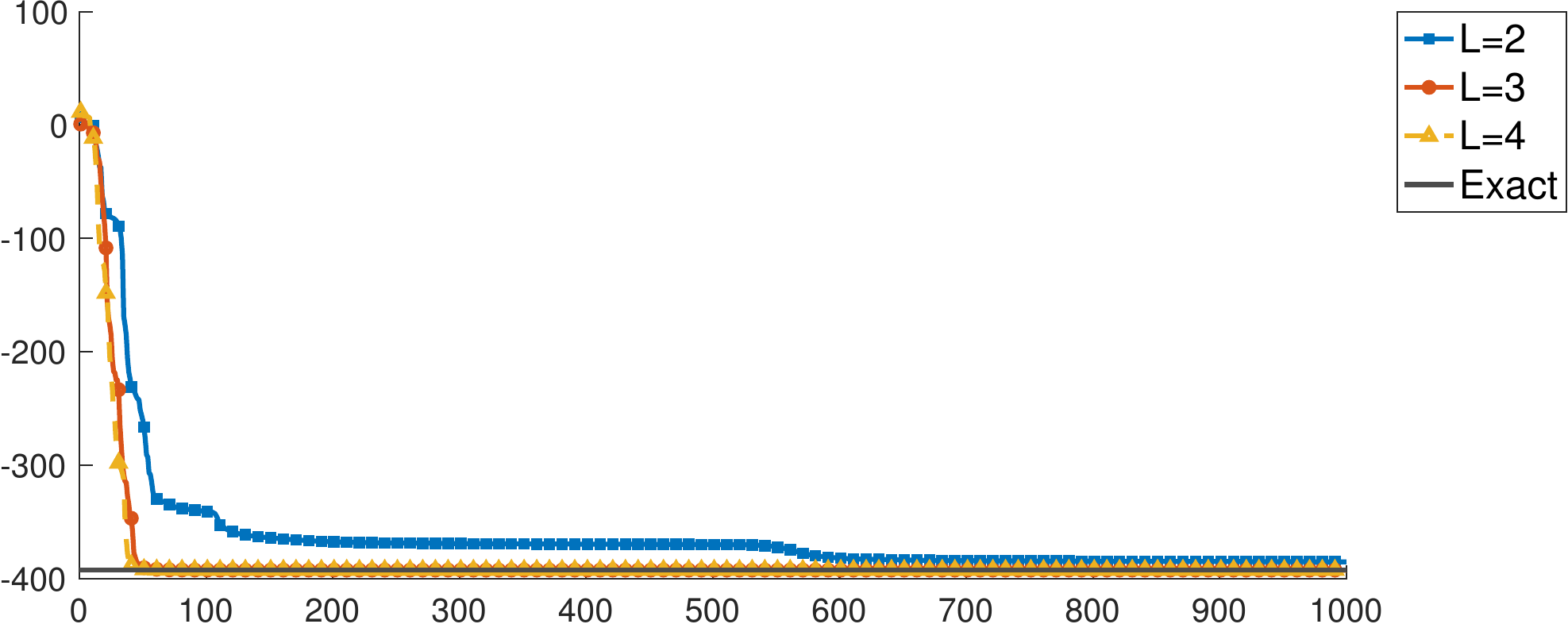}};
        \node[align=center,fill=white] at (3.1,-1.0) {(a)};
        \end{tikzpicture}
    \end{subfigure}
    \quad
    \begin{subfigure}[htbp]{0.475\textwidth}
        \begin{tikzpicture}
        \node[inner sep=0pt] at (0,0)
        {\includegraphics[width=\textwidth]{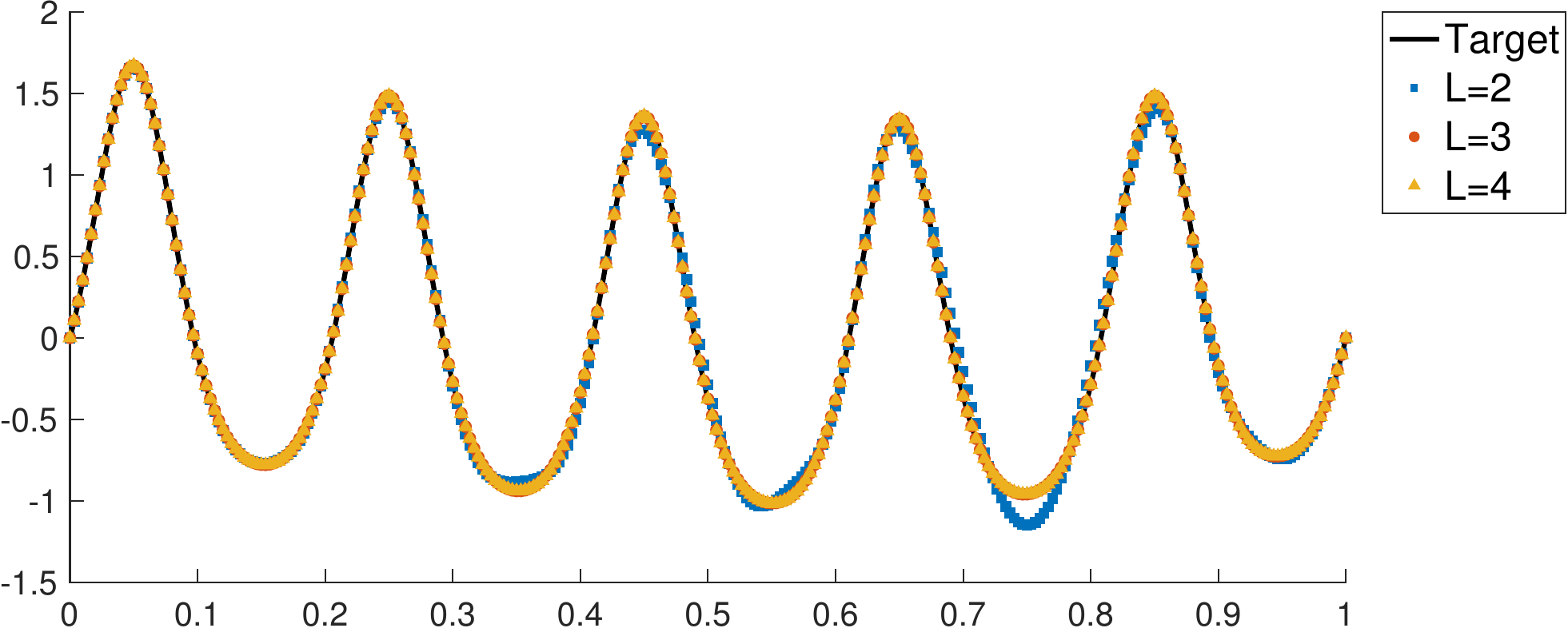}};
        \node[align=center,fill=white] at (3.1,-1.0) {(b)};
        \end{tikzpicture}
    \end{subfigure}
	\caption{PDE-Case II ($k=10$): (a) Training loss history for NNs of fixed depths $L$, along with the exact Ritz energy (dark line); (b) Final NN approximations on the test set. All models were trained using the NGF optimizer.}
\label{fig:possion1dk10_ngf_fixedL}
\end{figure}

\begin{table}[!ht]
\centering
\begin{tabular}{c|cccc|cccc}
\hline
\multirow{2}{*}{$L$} & \multicolumn{4}{c|}{Adam} & \multicolumn{4}{c}{NGF} \\ \cline{2-9}
{}    & $N_\text{ite}$ & Loss & $L_2$ error & $H_1$ error  & $N_\text{ite}$ & Loss & $L_2$ error & $H_1$ error\\ \hline
$2$ & $10^4$ & $-391.76$ & $1.22\times 10^{-2}$ & $1.05$ & $10^3$  & $-385.01$ & $5.61\times 10^{-2}$ & $3.82$\\
$3$ & $10^4$ & $-382.90$ & $7.21\times 10^{-2}$ & $4.33$ & $10^3$  & $-392.28$ & $1.24\times 10^{-3}$ & $1.01\times 10^{-1}$\\
$4$ & $10^4$ & $-391.69$ & $1.51\times 10^{-2}$ & $1.12$ & $10^3$  & $-392.28$ & $1.41\times 10^{-3}$ & $1.44\times 10^{-1}$\\
\hline
\end{tabular}
\caption{NN performance for PDE-Case II ($k=10$): Comparison of training convergence ($N_\text{ite}$ iterations), training loss, and approximation accuracy ($L_2$ and $H_1$ norms) across depths $L$ using Adam/NGF.}
\label{tab:poissonk5_ngf_gd}
\end{table}

Next, we train the expansive NN model with initial depth of $L_0=2$. The same training strategy is applied, and the previously discussed two types of initialization methods are tested. 

\begin{figure}[!htb]
  \centering
  	\begin{subfigure}[htbp]{0.475\textwidth}
        \begin{tikzpicture}
        \node[inner sep=0pt] at (0,0)
        {\includegraphics[width=\textwidth]{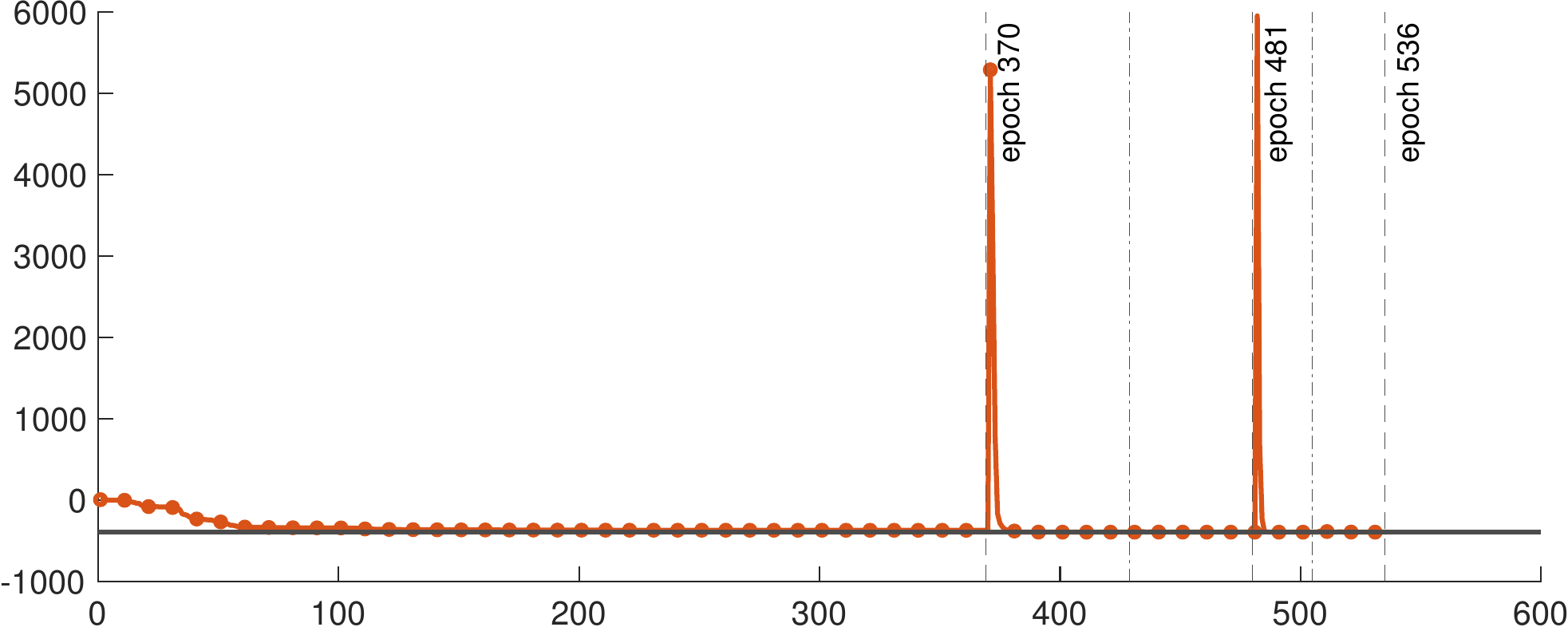}};
        \node[align=center,fill=white] at (0.2,1.1) {(a)};
        \end{tikzpicture}
    \end{subfigure}
    \quad
    \begin{subfigure}[htbp]{0.475\textwidth}
        \begin{tikzpicture}
        \node[inner sep=0pt] at (0,0)
        {\includegraphics[width=\textwidth]{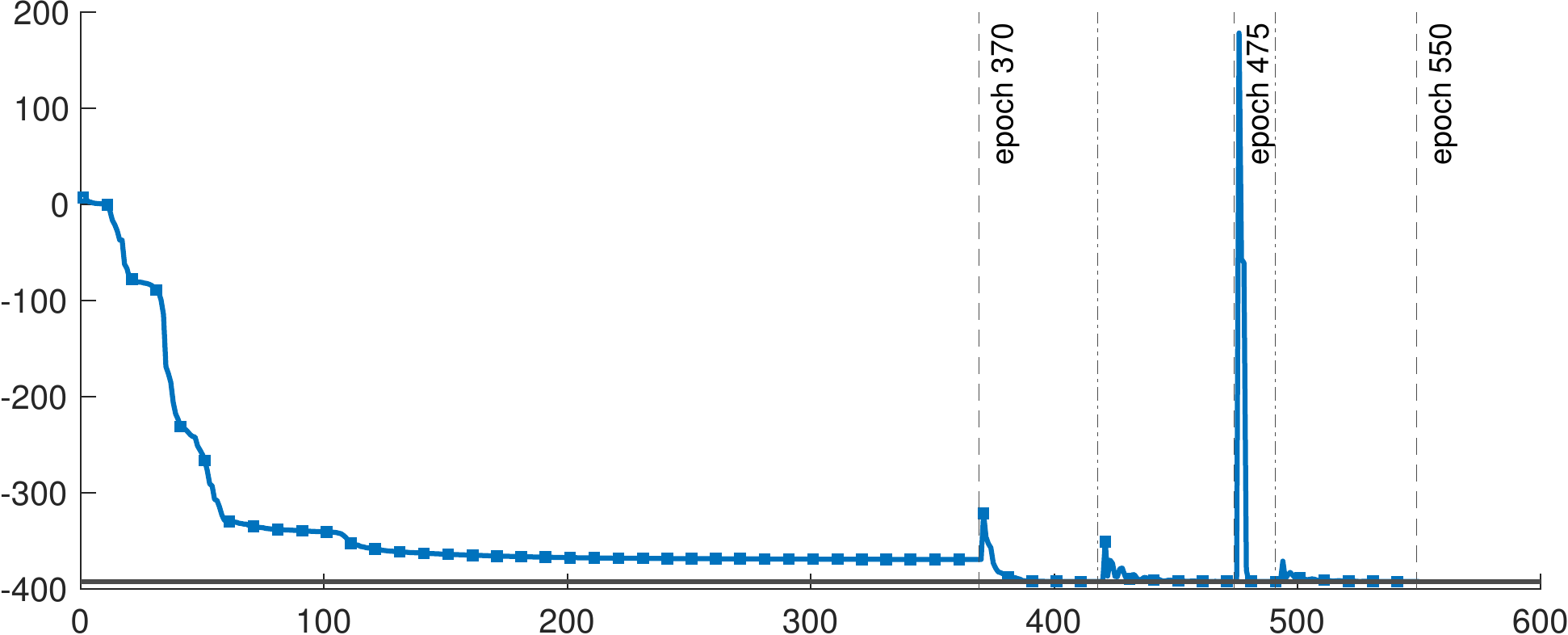}};
        \node[align=center,fill=white] at (0.2,1.1) {(b)};
        \end{tikzpicture}
    \end{subfigure}
	\caption{PDE-Case II ($k=10$): expansive NN with initial depth $L_0=2$, comparing (a) random vs. (b) gradient-aligned initialization. Vertical lines indicate expansion events (dashed) and NGF-to-Adam transitions (dash-dot). The dark horizontal line indicates the exact Ritz energy.}
\label{fig:poisson1dk10_expNGF_lastTwo}
\end{figure}

\begin{table}[!ht]
\centering
\begin{tabular}{c|ccccc}
\hline
\multirow{2}{*}{Initalization} & \multicolumn{5}{c}{Expansive NN with $L_0=2$} \\ \cline{2-6} 
{}    & $N_\text{ite}$ & $N_{\text{exp}}$ & Final loss & $L_2$ error & $H_1$ error \\ \hline
Random           & $536$ & 2 & $-391.97$ & $6.31\times 10^{-2}$ & $7.23\times 10^{-1}$ \\
Gradient-aligned & $550$ & 2 & $-392.17$ & $1.38\times 10^{-2}$ & $5.21\times 10^{-1}$ \\
\hline
\end{tabular}
\caption{Expansive NN performance for PDE-Case II ($k=10$): Comparison of training convergence ($N_\text{ite}$ iterations), number of expansions ($N_\text{exp}$), and approximation accuracy ($L_2$ and $H_1$ norms).}
\label{tab:poisson1dk10_exp}
\end{table}

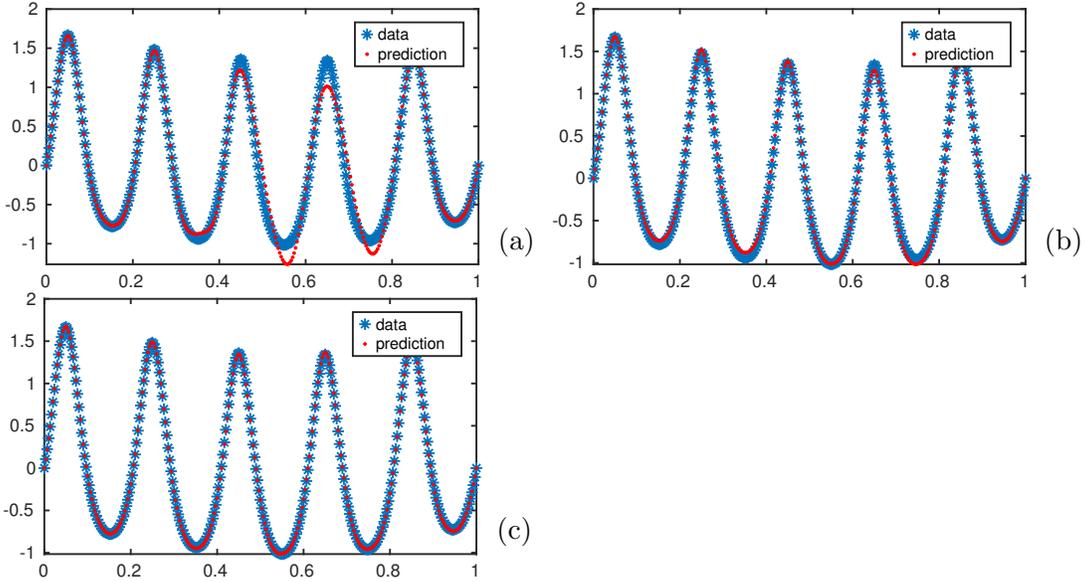
\begin{figure}[!htb]
  \centering
  	\begin{subfigure}[htbp]{0.425\textwidth}
        \begin{tikzpicture}
        \node[inner sep=0pt] at (0,0)
        {\includegraphics[width=\textwidth]{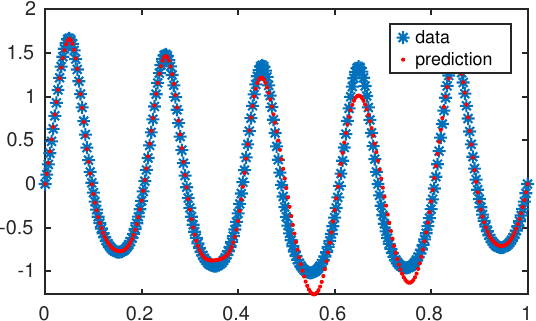}};
        \node[align=center,fill=white] at (3.6,-1.3) {(a)};
        \end{tikzpicture}
    \end{subfigure}
    \qquad
    \begin{subfigure}[htbp]{0.425\textwidth}
        \begin{tikzpicture}
        \node[inner sep=0pt] at (0,0)
        {\includegraphics[width=\textwidth]{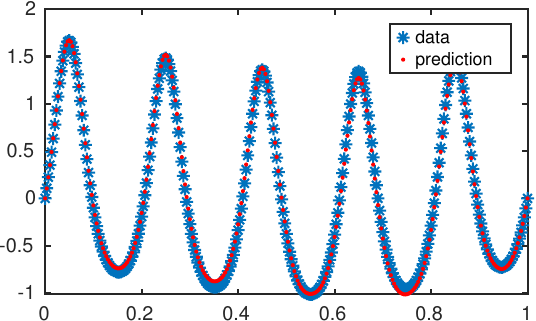}};
        \node[align=center,fill=white] at (3.6,-1.3) {(b)};
        \end{tikzpicture}
    \end{subfigure}
    \\
    \begin{subfigure}[htbp]{0.425\textwidth}
        \begin{tikzpicture}
        \node[inner sep=0pt] at (0,0)
        {\includegraphics[width=\textwidth]{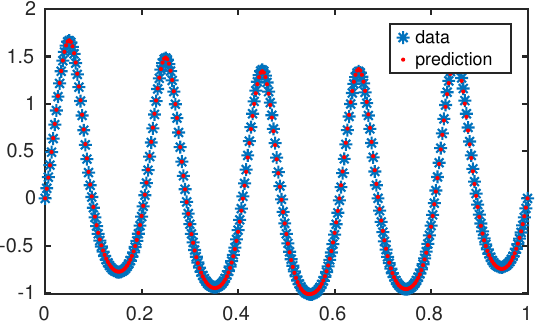}};
        \node[align=center,fill=white] at (3.6,-1.3) {(c)};
        \end{tikzpicture}
    \end{subfigure}
    \hspace{0.48\textwidth}
	\caption{PDE-Case II ($k=10$): (a)-(b) Approximation of the expansive NN model on the test set before the first and second expansion, respectively; (b) Final approximation on the test set.}
\label{fig:poisson1dk10_expNGF_lastTwo_align}
\end{figure}

\newpage
\subsection{MOR for Parametrized Burgers}\label{sec:pburg} 
Next, we consider solving $u(x, t; \mu)$ from a parametrized inviscid Burgers equation:
\begin{equation}
u_t + uu_x = 0.02e^{\mu x}, \quad \forall x\in [0, 1],\, t\in[0, 20],
\label{eq:pburg}
\end{equation}
associated with the initial condition $u(x, 0; \mu) = \left\{\begin{array}{cl} 4.25, & \text{ if } x=0 \\ 1, & \text{otherwise} \end{array}\right.$, boundary condition $u(0, t; \mu)=4.25$ for $t>0$, and the parameter $\mu \in [0.015, 0.030]$. 
The solution's behavior varies with the change of $\mu$. For example, numerical solutions at two parameters $\mu = 0.0174$ and $0.0296$ are presented in Figure \ref{fig:burgers_fd}.

\begin{figure}[!htb]
  \centering
  	\begin{subfigure}[htbp]{0.475\textwidth}
        \begin{tikzpicture}
        \node[inner sep=0pt] at (0,0)
        {\includegraphics[width=\textwidth]{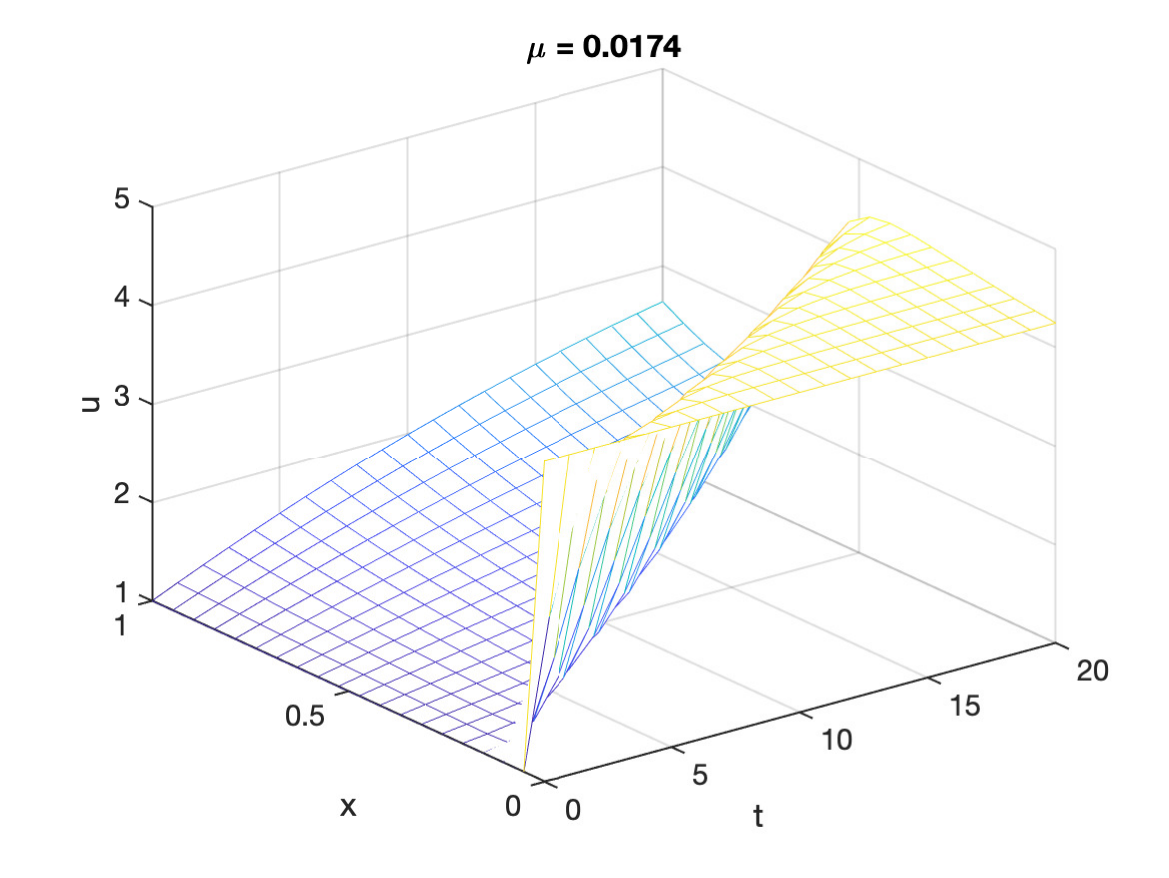}};
        \node[align=center,fill=white] at (3.2,1.1) {(a)};
        \end{tikzpicture}
    \end{subfigure}
    \quad
    \begin{subfigure}[htbp]{0.475\textwidth}
        \begin{tikzpicture}
        \node[inner sep=0pt] at (0,0)
        {\includegraphics[width=\textwidth]{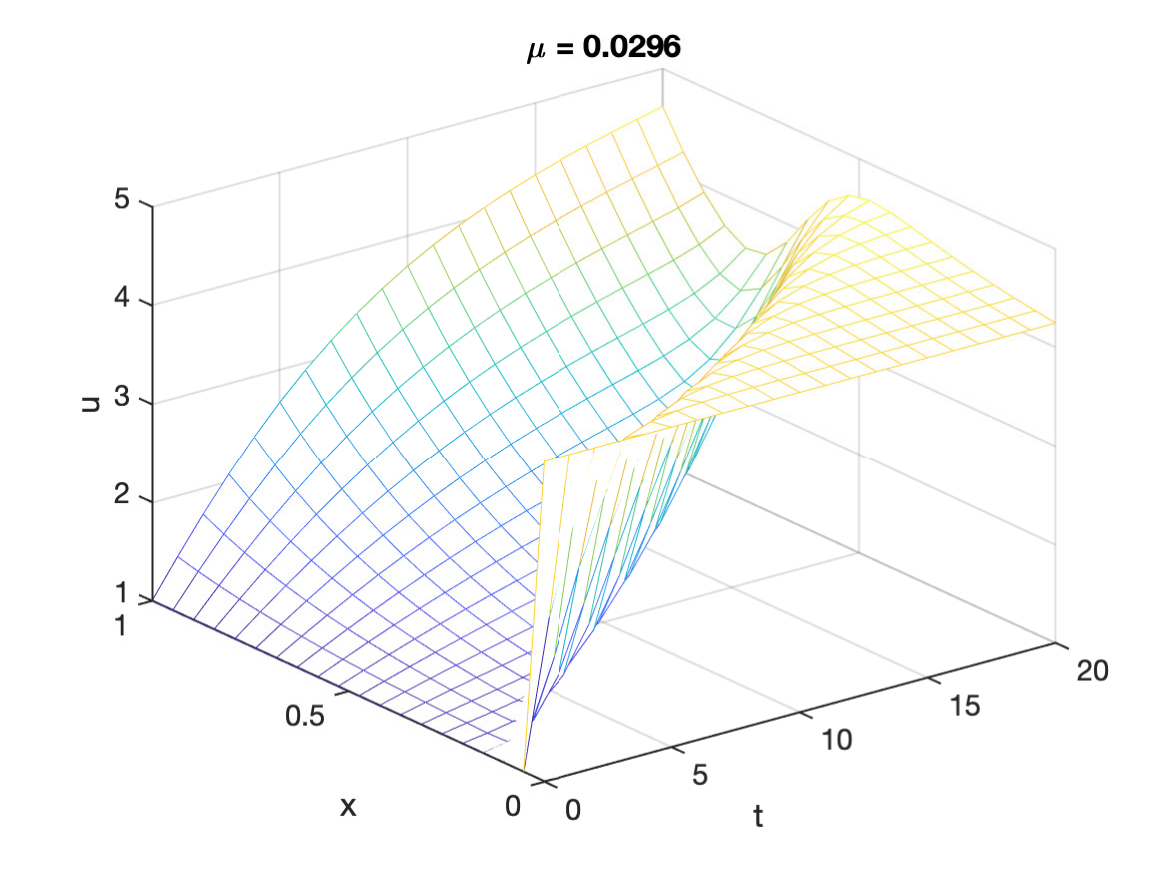}};
        \node[align=center,fill=white] at (3.2,1.1) {(b)};
        \end{tikzpicture}
    \end{subfigure}
	\caption{Full order simulations of the parametrized Burgers equation at (a) $\mu=0.0174$ and (b) $\mu = 0.0296$, respectively.}
\label{fig:burgers_fd}
\end{figure}

Our goal is to approximate the parametric solution by a NN model $f(x, t, \mu; \Theta)$ via supervised learning. 
To generate training data, numerical simulations by finite difference method using a uniform spatial and temporal discretization $\delta x = \delta t = \frac{1}{20}$ are performed at 11 uniformly spaced parameter samples from the prescribed interval. The numerical solutions on the grids are collected once every 20 time steps, which yields $M_{train}=4,620$ samples. 
This set is partitioned into 11 mini-batches for training. The integrals involved in the optimization process are calculated by Monte-Carlo (MC) method.
For testing the performance of trained neural networks, another 6 $\mu$-values from the interval $[0.015, 0.030]$ are randomly selected. Using the same fully-discrete numerical scheme, we generate solutions and form a test set containing $M_{test} = 2,520$ samples.  

The same type of fully connected NN with residual connections is applied, whose input dimension is $3$ and output dimension is $1$. All hidden layers have a uniform width of $N=15$.  
Training is performed using batches of size $N_b=420$, which will be terminated when the loss value reaches 
{\color{blue} $\tau_{t_a} = 10^{-5}$.}

We first compare the performance of Adam and NGF in training the NN with a fixed depth $L$. 
The training loss history of Adam is shown in Figure~\ref{fig:mor_pburg_fixedL_adam}(a). The results indicate that the loss decreases as $L$ increases; however, none of them reach the prescribed loss tolerance. 
In Figure~\ref{fig:mor_pburg_fixedL_adam}(b), we present the NN approximate solutions at the final time corresponding to two test parameter samples, $\mu = 0.0174$ and $\mu = 0.0296$. We observe that the NN approximation improves   as the depth increases from $L=2$ to $L=4$.
\begin{figure}[!htb]
  \centering
  	\begin{subfigure}[htbp]{0.475\textwidth}
        \begin{tikzpicture}
        \node[inner sep=0pt] at (0,0)
        {\includegraphics[width=\textwidth]{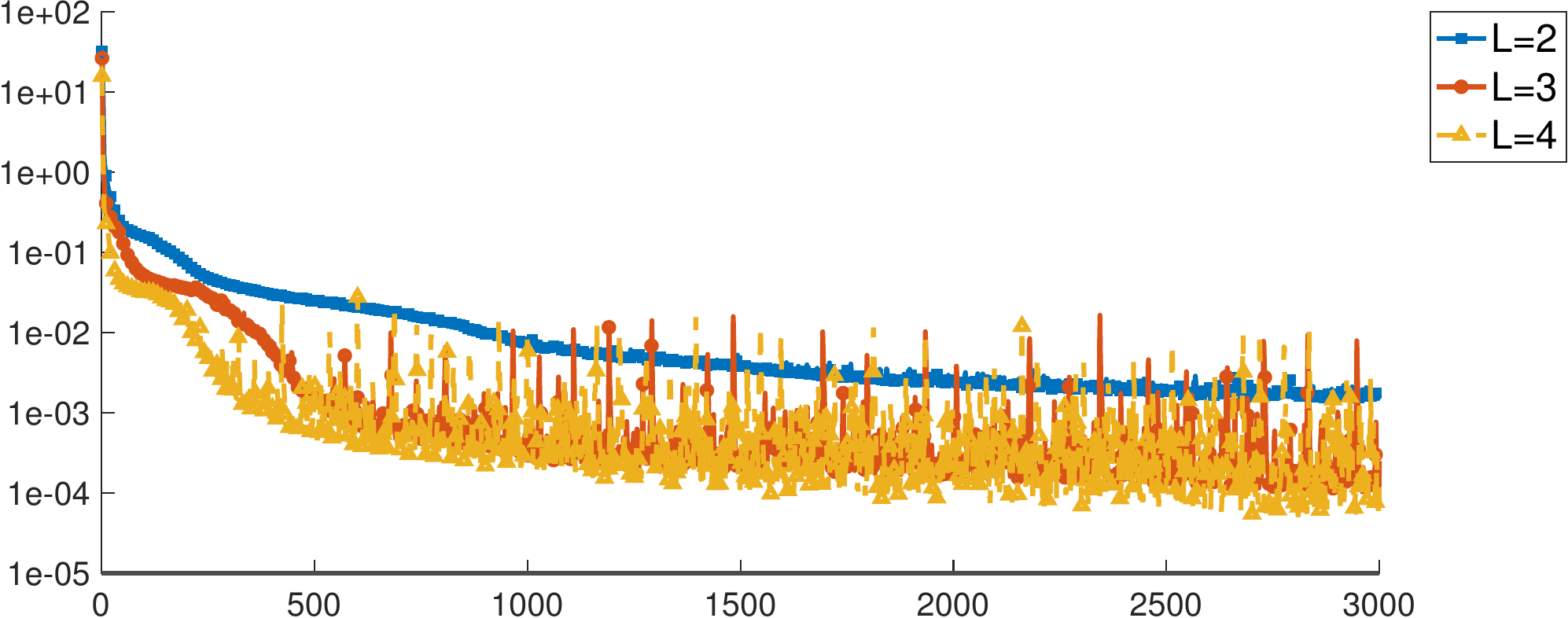}};
        \node[align=center,fill=white] at (3.3,-1.0) {(a)};
        \end{tikzpicture}
    \end{subfigure}
    \quad 
    \begin{subfigure}[htbp]{0.475\textwidth}
        \begin{tikzpicture}
        \node[inner sep=0pt] at (0,0)
        {\includegraphics[width=\textwidth]{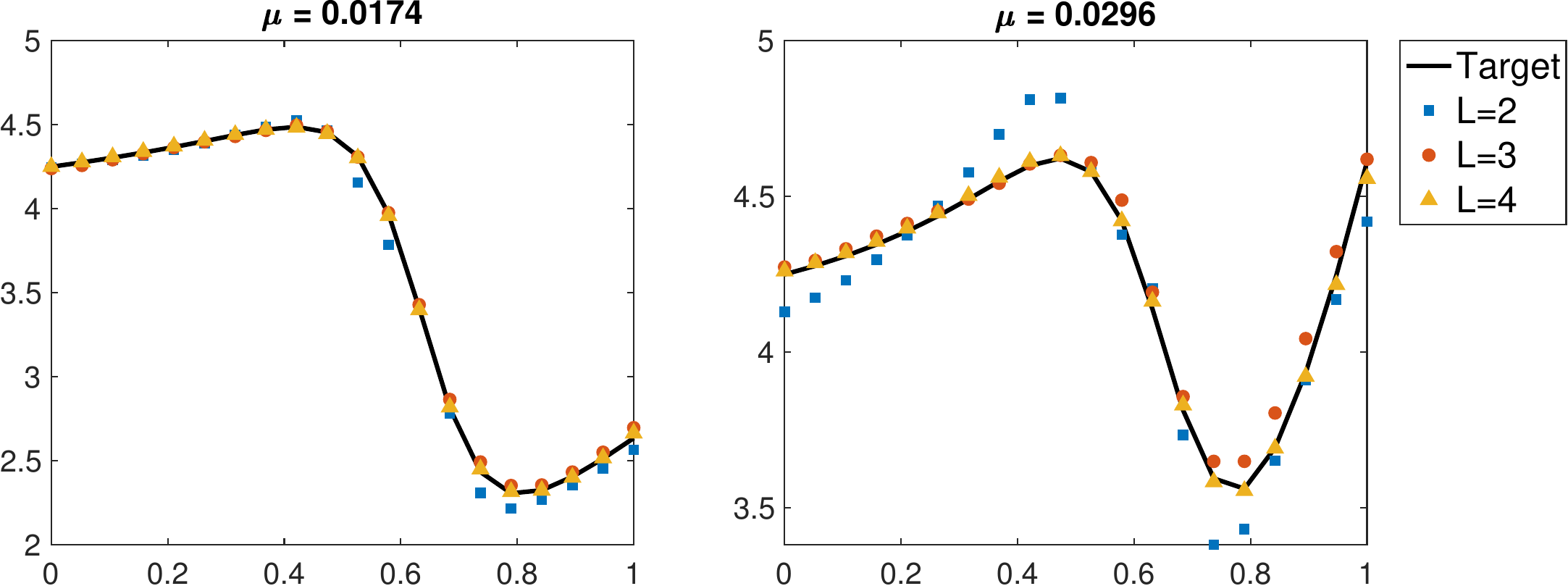}};
        \node[align=center,fill=white] at (3.3,-1.0) {(b)};
        \end{tikzpicture}
    \end{subfigure}
	\caption{MOR: Neural network of fixed depth $L$ trained using Adam optimizer. (a) History of training loss; (b) Approximate solutions at two test parameters.}
\label{fig:mor_pburg_fixedL_adam}
\end{figure}

When NGF is applied, we solve $\left(G_{T,\Theta}+\lambda I\right)\Delta \theta = -\nabla_\theta E(\theta)$ for the gradient update. The regularization parameter $\lambda$ follows the same definition as in \eqref{eq:lambda}, with $\lambda_j = 1\times 10^{j-8}$, for $j=1,\ldots, 6$.
The training process with NGF converges faster than Adam, yet it still fails to attain the prescribed tolerance within the maximum allowed iterations for all tested depths, as illustrated in Figure~\ref{fig:mor_pburg_fixedL_ngf}(a). The corresponding NN  approximations are shown in Figure~\ref{fig:mor_pburg_fixedL_ngf}(b). 
\begin{figure}[!htb]
  \centering
  	\begin{subfigure}[htbp]{0.475\textwidth}
        \begin{tikzpicture}
        \node[inner sep=0pt] at (0,0)
        {\includegraphics[width=\textwidth]{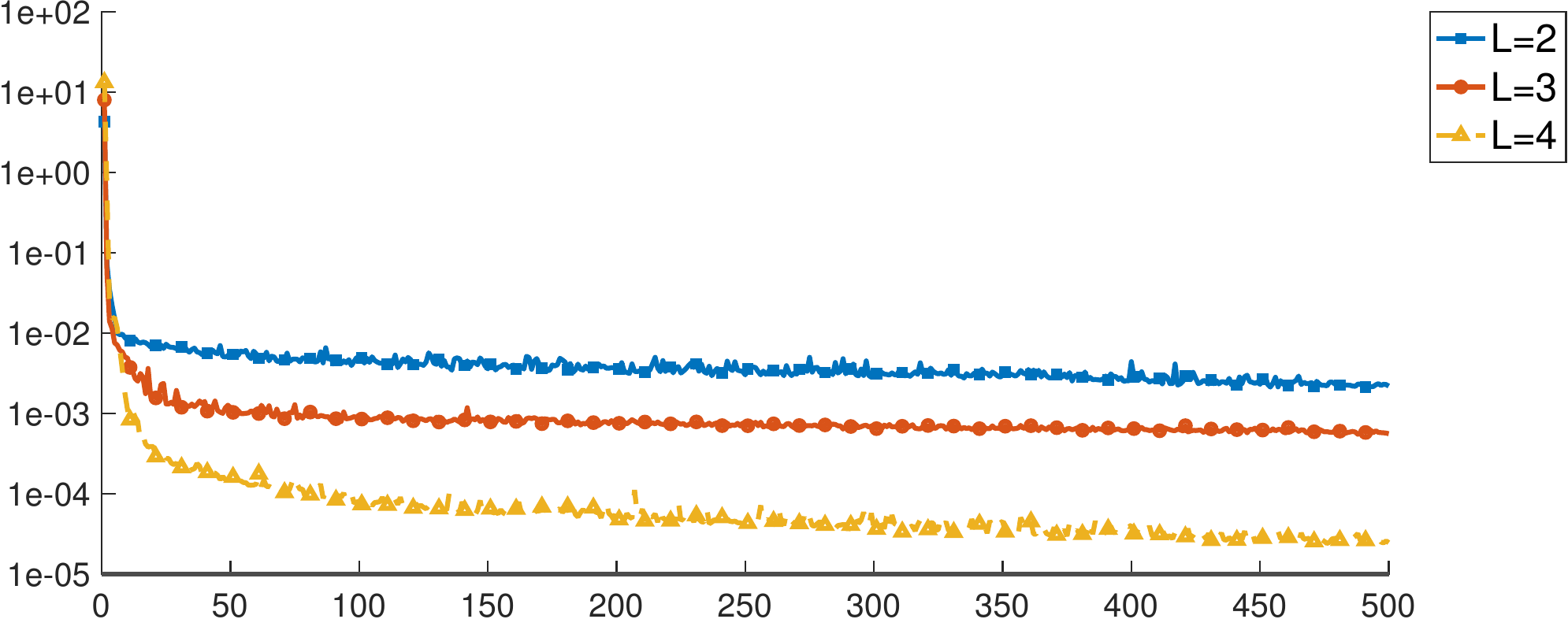}};
        \node[align=center,fill=white] at (3.3,-1.0) {(a)};
        \end{tikzpicture}
    \end{subfigure}
    \quad 
    \begin{subfigure}[htbp]{0.475\textwidth}
        \begin{tikzpicture}
        \node[inner sep=0pt] at (0,0)
        {\includegraphics[width=\textwidth]{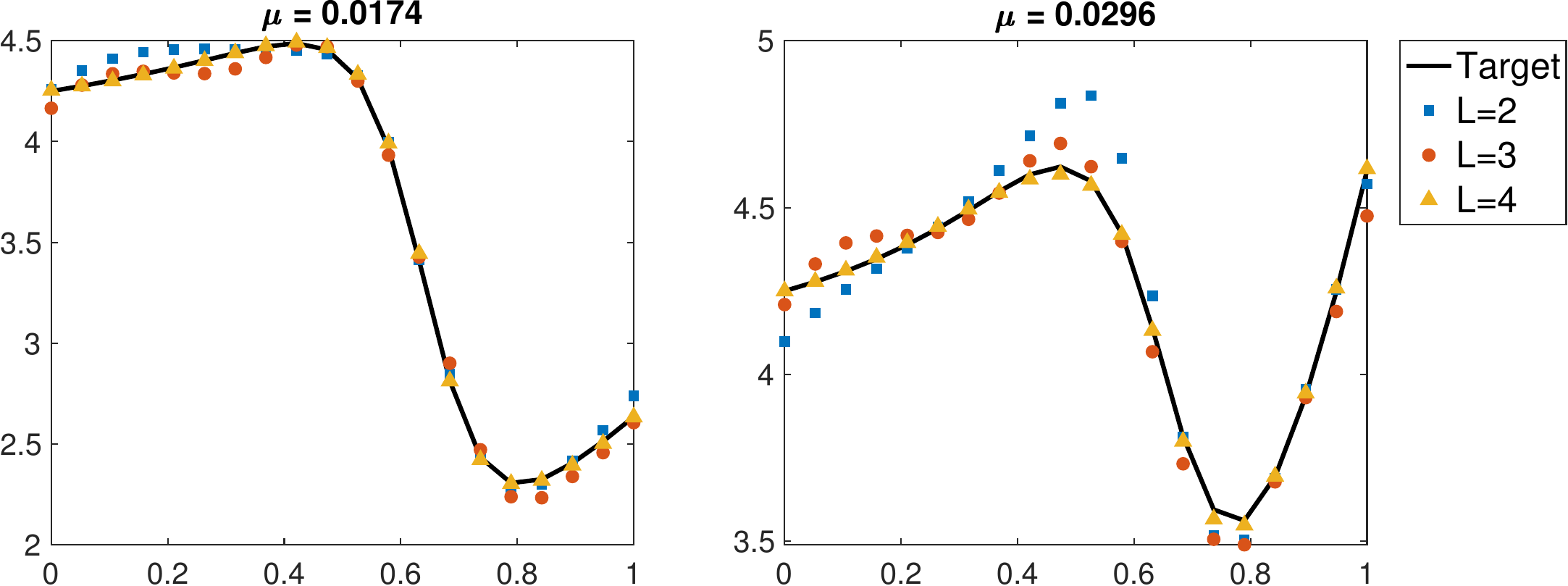}};
        \node[align=center,fill=white] at (3.3,-1.0) {(b)};
        \end{tikzpicture}
    \end{subfigure}
	\caption{MOR: Neural network of fixed depth $L$ trained using NGF optimizer. (a) History of training loss, where the dark line indicates the prescribed loss tolerance; (b) Approximate solutions at two test parameters.}
\label{fig:mor_pburg_fixedL_ngf}
\end{figure}

For both optimizers, the training iterations, final loss values, and approximation errors of the trained NN models, measured in $L_2$ norm, are listed in Table~\ref{tab:mor_pburg_ngf_gd}. 
\begin{table}[!ht]
\centering
\begin{tabular}{c|ccc|ccc}
\hline
\multirow{2}{*}{$L$} & \multicolumn{3}{c|}{Adam} & \multicolumn{3}{c}{NGF} \\ \cline{2-7}
{}    & $N_\text{ite}$ & Loss & $L_2$ error & $N_\text{ite}$ & Loss & $L_2$ error\\ \hline
$2$ & $10^4$ & $1.02\times 10^{-3}$ & $4.25\times 10^{-2}$ & $10^3$  & $1.72\times 10^{-3}$ & $6.55\times 10^{-2}$ \\
$3$ & $10^4$ & $4.30\times 10^{-4}$ & $1.46\times 10^{-2}$ & $10^3$  & $3.16\times 10^{-4}$ & $2.65\times 10^{-2}$ \\
$4$ & $10^4$ & $2.86\times 10^{-5}$ & $5.60\times 10^{-3}$ & $10^3$  & $1.70\times 10^{-5}$ & $5.54\times 10^{-3}$ \\
\hline
\end{tabular}
\caption{NN performance for MOR: Comparison of training convergence ($N_\text{ite}$ iterations), training loss, and approximation accuracy ($L_2$ norm) across depths $L$ using Adam/NGF.}
\label{tab:mor_pburg_ngf_gd}
\end{table}

Next, we train the expansive NN model with an initial depth of $L_0=2$. The maximum number of expansion $M_\text{exp}=6$, the maximum number of iteration $3\times 10^3$, and the loss tolerance $\tau_{t_a}=10^{-5}$ are used in training. { In the loss saturation criterion,    $\tau_{s_a} = 10^{-7}$ and $\tau_{s_r} = 5\times 10^{-3}$ are used for NGF optimizer and $\tau_{s_a} = 10^{-8}$ and $\tau_{s_r} = 5\times 10^{-4}$ for Adam optimizer.} The same expansion strategy as used before is applied, and the previously discussed two types of initialization methods are tested. The corresponding loss histories are shown in Figure~\ref{fig:mor_pburg_expNGF_lastTwo}, where (a) corresponds to the random initialization and (b) the gradient-aligned initialization. We observe   that the gradient-aligned initialization method outperforms the random initialization, achieving the prescribed loss tolerance after 4 times of expansions, while random initialization fails to converge to the target accuracy despite reaching the maximum of 6 expansions and allowed iterations. Table~\ref{tab:mor_pburg_exp} summarizes the number of iterations and expansions, final loss values, and numerical errors. Finally, we present the NN approximations in the gradient-aligned initialization case at the test parameters $\mu = 0.0174$ and $0.0296$ in Figure~\ref{fig:mor_pburg_expNGF_lastTwo_align}. 

\begin{figure}[!htb]
  \centering
  	\begin{subfigure}[htbp]{0.475\textwidth}
        \begin{tikzpicture}
        \node[inner sep=0pt] at (0,0)
        {\includegraphics[width=\textwidth]{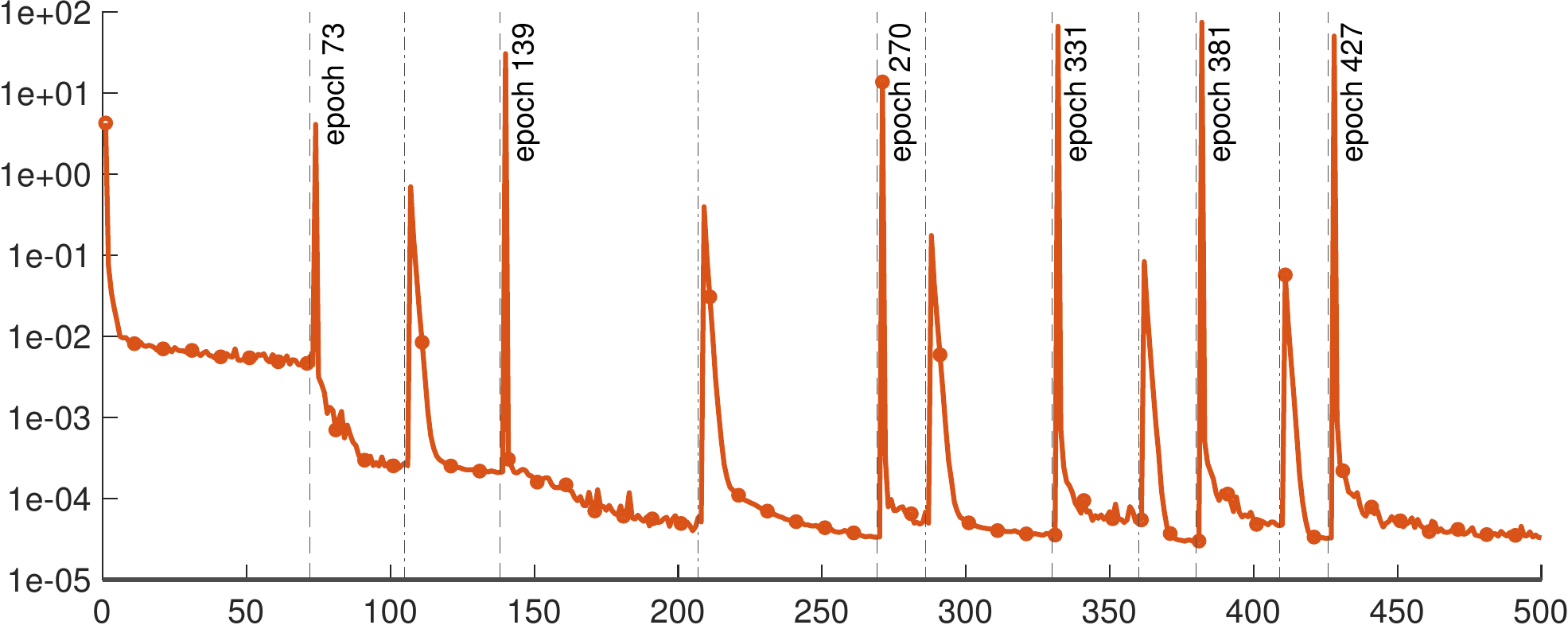}};
        \node[align=center,fill=white] at (0.2,1.1) {(a)};
        \end{tikzpicture}
    \end{subfigure}
    \quad
    \begin{subfigure}[htbp]{0.475\textwidth}
        \begin{tikzpicture}
        \node[inner sep=0pt] at (0,0)
        {\includegraphics[width=\textwidth]{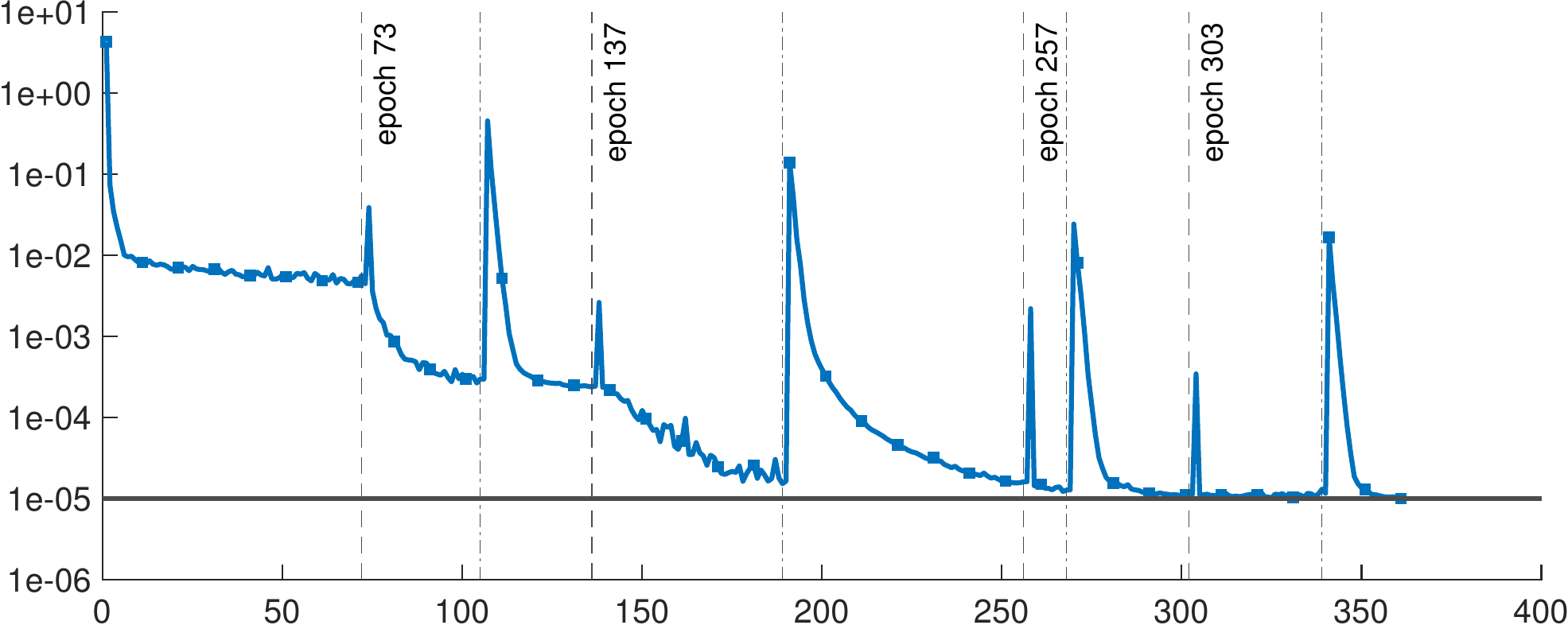}};
        \node[align=center,fill=white] at (0.2,1.1) {(b)};
        \end{tikzpicture}
    \end{subfigure}
	\caption{MOR: Expansive NN with initial depth $L_0=2$, comparing (a) random vs. (b) gradient-aligned initialization. Vertical lines indicate expansion events (dashed) and NGF-to-Adam transitions (dash-dot).}
\label{fig:mor_pburg_expNGF_lastTwo}
\end{figure}

\begin{table}[!ht]
\centering
\begin{tabular}{c|cccc}
\hline
\multirow{2}{*}{Initalization} & \multicolumn{4}{c}{Expansive NN with $L_0=2$} \\ \cline{2-5} 
{}    & $N_\text{ite}$ & $N_{\text{exp}}$ & Final loss & $L_2$ error \\ \hline
Random           & $3\times 10^3$ & 6 & $2.52\times 10^{-5}$ & $7.03\times 10^{-3}$  \\
Gradient-aligned & $362$ & 4 & $9.94\times 10^{-6}$ & $4.68\times 10^{-3}$  \\
\hline
\end{tabular}
\caption{Expansive NN performance for MOR: Comparison of training convergence ($N_\text{ite}$ iterations), number of expansions ($N_\text{exp}$), and approximation accuracy ($L_2$ norm).}
\label{tab:mor_pburg_exp}
\end{table}

\begin{figure}[!htb]
  \centering
  	\begin{subfigure}[htbp]{0.225\textwidth}
        \begin{tikzpicture}
        \node[inner sep=0pt] at (0,0)
        {\includegraphics[width=\textwidth,height=\textwidth]{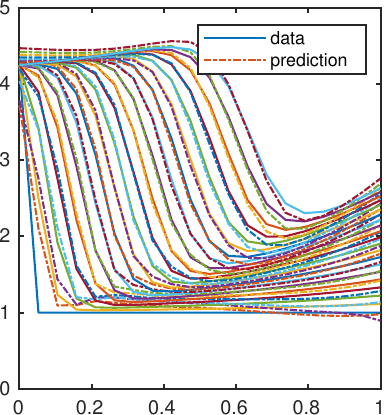}};
        \node[align=center,fill=white] at (2.1,-1.3) {(a)};
        \end{tikzpicture}
    \end{subfigure}
    \qquad
    \begin{subfigure}[htbp]{0.225\textwidth}
        \begin{tikzpicture}
        \node[inner sep=0pt] at (0,0)
        {\includegraphics[width=\textwidth,height=\textwidth]{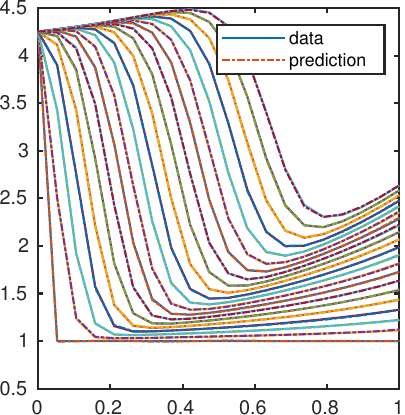}};
        \node[align=center,fill=white] at (2.1,-1.3) {(b)};
        \end{tikzpicture}
    \end{subfigure}
    \qquad
    \begin{subfigure}[htbp]{0.225\textwidth}
        \begin{tikzpicture}
        \node[inner sep=0pt] at (0,0)
        {\includegraphics[width=\textwidth,height=\textwidth]{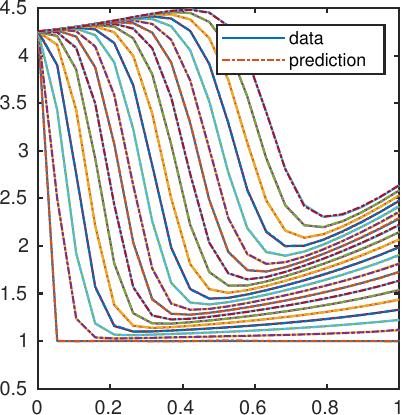}};
        \node[align=center,fill=white] at (2.1,-1.3) {(c)};
        \end{tikzpicture}
    \end{subfigure}
    \\
      	\begin{subfigure}[htbp]{0.225\textwidth}
        \begin{tikzpicture}
        \node[inner sep=0pt] at (0,0)
        {\includegraphics[width=\textwidth,height=\textwidth]{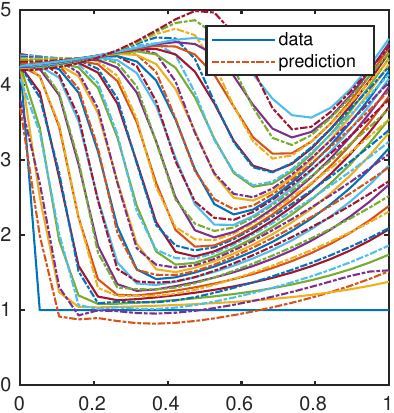}};
        \node[align=center,fill=white] at (2.1,-1.3) {(d)};
        \end{tikzpicture}
    \end{subfigure}
    \qquad
    \begin{subfigure}[htbp]{0.225\textwidth}
        \begin{tikzpicture}
        \node[inner sep=0pt] at (0,0)
        {\includegraphics[width=\textwidth,height=\textwidth]{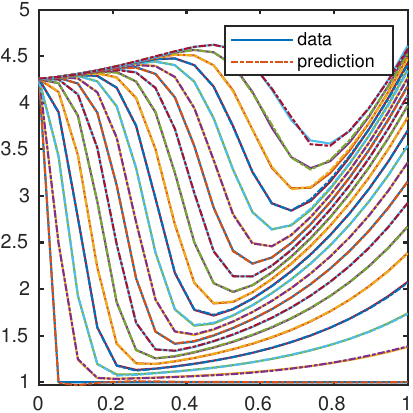}};
        \node[align=center,fill=white] at (2.1,-1.3) {(e)};
        \end{tikzpicture}
    \end{subfigure}
    \qquad
    \begin{subfigure}[htbp]{0.225\textwidth}
        \begin{tikzpicture}
        \node[inner sep=0pt] at (0,0)
        {\includegraphics[width=\textwidth,height=\textwidth]{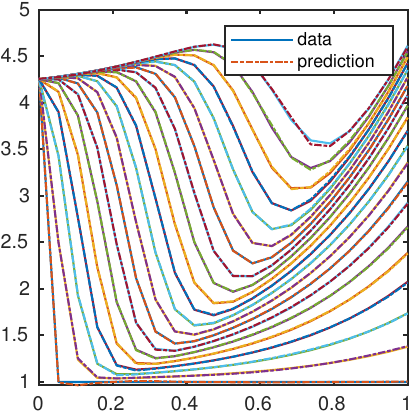}};
        \node[align=center,fill=white] at (2.1,-1.3) {(f)};
        \end{tikzpicture}
    \end{subfigure}
	\caption{MOR: Approximations of the expansive NN model for test parameters $\mu = 0.0174$ (first row) and $\mu = 0.0296$ (second row), where (a) and (d) show pre-first-expansion, (b) and (e) pre-third-expansion, and (c) and (f) final results.}
\label{fig:mor_pburg_expNGF_lastTwo_align}
\end{figure}
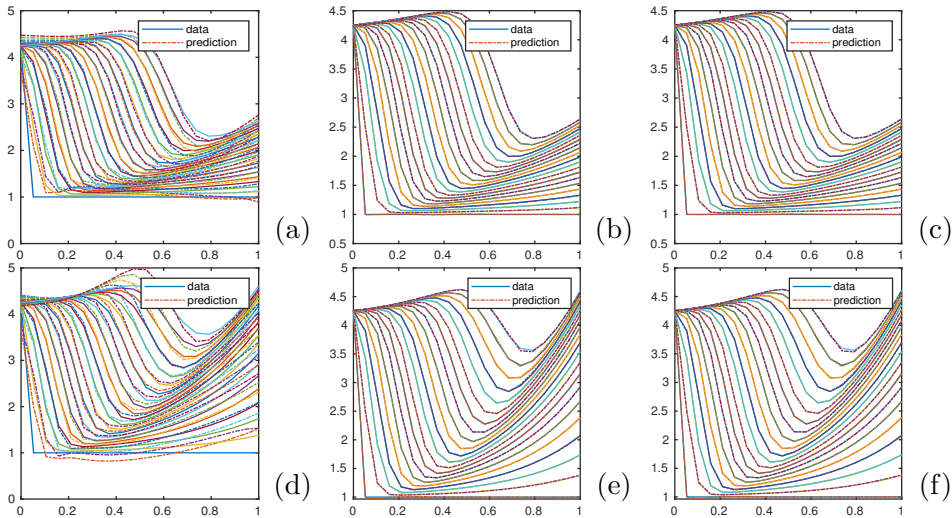

\section{Summary and Conclusions}
We have presented a general setting that accommodates both supervised and
unsupervised learning tasks in the form of minimizing a strictly convex (in fact quadratic)
energy in some ambient (possibly infinite dimensional) Hilbert space. Since this setting primarily applies to physics-informed tasks and thus aims, in one way or another, 
at producing (ideally certified) reduced models, over-parametrization as a means to
ensure optimization accuracy is ruled out. The overarching objective is therefore, on the one hand, to
achieve a given accuracy target at the expense of possibly small budgets of trainable weights, and on the other hand, to see to what extent achievable accuracy is affected
by specific optimization strategies. The proposed approach hinges on two main pillars:
(a) natural gradient flow, (b) adaptive network expansion. Regarding (a), we stress the importance of respecting the correct inner product imposed by the ambient Hilbert space,
not just $L_2$ as a convenient default which, in particular, leads to accuracy certificates in model compliant norms (in the sense of generalization errors) when $\cE(u^*)=0$.  (b) aims at actually realizing a target tolerance by gradually increasing network complexity. Moreover, the role of the ambient Hilbert space
becomes particularly evident in connection with the {\em adaptation criterion}.
As brought out in our discussions the expansion is to better align the tangent space of the
expanded neural manifold with the {\em ideal Hilbert space gradient direction}, leveraging the fact that projections of this gradient (although itself not being computable) can be
evaluated. Hence, this serves two purposes: first, to increase expressive power when needed, and second, to align the evolution with the ideal Hilbert space descent path, thereby enhancing minimization efficacy.

We have presented preliminary tests of these concepts for a collection of (admittedly simple) model problems. Nevertheless, these examples highlight the scope of the setting under consideration and exhibit the challenges that come with respecting the correct norms.
Moreover, we have chosen the simplest versions of the proposed strategies.
Although these relatively simple means show a clear effect, especially alignment strategies leave much room for improvement, and several such options have been already indicated.
We finally emphasize that the proposed concepts are not restricted to the 1D examples because the  expensive part of assembling flow matrices is controlled by
the proposed hybrid strategy of confining NGF steps to two layers, while ADAM serves to
update  all parameters in an interlacing fashion.


\end{document}